\documentclass[10pt]{article}

\usepackage{amsmath,latexsym,amssymb,amsthm}

\numberwithin{equation}{section} \allowdisplaybreaks

\newtheorem{theorem}{\bf\normalsize Theorem}[section]
\newtheorem{corollary}[theorem]{\bf\normalsize Corollary}

\newtheorem{lemma}[theorem]{\bf\normalsize Lemma}
\newtheorem{proposition}[theorem]{\bf\normalsize Proposition}
\newtheorem{remark}[theorem]{\bf\normalsize Remark}

\begin{document}

\title{The  Gaussian Correlation
Inequality\\ for  Symmetric Convex Sets}
\author{QingYang Guan\footnote{
 Institute of Applied Mathematics,  AMSS, CAS.
   }
   }
\date{} \maketitle

\begin{abstract}
The paper is to prove the Gaussian correlation conjecture stating
that, under  the standard    Gaussian measure, the measure of the
intersection of any two  symmetric convex sets  is greater than or
equal to the product of their measures. Characterization of the
equality and some applications are given.
\end{abstract}


{\textbf{\small\noindent AMS(2000) Subject Classification:
  Primary $60\mathrm{G}15$;\  Secondary
$28\mathrm{C}20$, $60\mathrm{E}15$ }

   \tableofcontents

\section{Introduction}
\subsection{\normalsize the main results}     The standard Gaussian
measure $\mu_n$ on $\mathbb{R}^n$ is given  by
\begin{align} d\mu_n(x)=\frac{1}{(2\pi)^{n/2}} {\exp{\{-|x|^2/2\}}}dx,
\ \ \ \ \
\end{align}
where $|x|$ is   the Euclidean norm of $x$. The main subject of the
paper is to prove the   conjecture    that, under
 the standard  Gaussian measure,
  a pair of symmetric convex sets are   positively correlated
  or independent.
  This conjecture  is often  called the Gaussian correlation conjecture
in   literatures.

Some special cases of the conjecture appeared first   in
 the study of    multidimensional confidence regions in
 statistics
 for
Gaussian measure. See, e.g.,    Dunn \cite{DUNN}, Dunnett and Sobel
\cite{D}, Khatri \cite{Kha67} and \v{S}id\'{a}k
\cite{Sid67}\cite{Sid68}. We refer to Das Gupta, Eaton, Olkin,
Perlman, Savage and Sobel \cite{DEOPSS72} and Schechtman,
Schlumprecht and Zinn \cite{SSZ98} for   more  historical background
of the conjecture.

 The statement of the conjecture in
  Theorem 1.1 below is from      Pitt \cite{Pitt77}.  See  \cite{SSZ98} for
  some  other  equivalent forms  of
  the
conjecture. Besides the original conjecture, characterization of the
equality is also given in  Theorem 1.1. Notice that counterexamples
of  a   stronger  conjecture in \cite{DEOPSS72} for general
 elliptically  contoured  distributions
are  stated   in Theorem 1.2 of Lewis and Pritchard \cite{LP03}.

 The following definitions
and notations are adopted in Theorem 1.1 and also    the rest part
of the paper. A subset $A$ of a Banach space  is called a  symmetric
set if its indicator function $I_A$ is an even function.     A
measurable subset $A$ of $ \mathbb{R}^n$ is called degenerate
  if $\mu_n(A)=0$. For $A\subseteq \mathbb{R}^n$, denote  by $\overline{A}$ the closure
of   $A$. Denote by $\mathcal{C}_n$ the class of
   symmetric convex subsets of $\mathbb{R}^n$.
 We call  two   subsets $A$ and $B$ of $\mathbb{R}^n$     unlinked   if either one of $A$ and  $B$
  is equal to
  $\mathbb{R}^n$,
  or there exists some  orthogonal
transformation  $Q$ of $\mathbb{R}^n$ such that
$Q(A)=\widetilde{A}\times \mathbb{R}^k$,
$Q(B)=\mathbb{R}^{n-k}\times\widetilde{B} $ for some
$\widetilde{A}\subseteq \mathbb{R}^{n-k}$, $\widetilde{B}\subseteq
\mathbb{R}^k$ with  $1\leq k\leq n-1$.

\begin{theorem}\label{A}   For every  $n\geq 1$ and  every   $A,B\in \mathcal{C}_n$
  \begin{align}\label{1}
\mu_n(A\cap B)\geq \mu_n(A)\mu_n(B).
\end{align}
Moreover, (\ref{1}) is  an    equality  if and only if  one of $A$
and  $B$ is
 degenerate or $\overline{A}$ and $\overline{B}$  are    unlinked.
\end{theorem}

  Various partial results of the conjecture have
been proved under additional conditions, e.g.,
 one of $A$ and $B$ is a symmetric slap   in \cite{Kha67} and
\cite{Sid67}; the two dimensional case in     \cite{Pitt77};
 both  $A$ and $B$ are contained in  the centered  ball
of radius
 $ 2^{-1/2}\Gamma(1+n/2)^{1/n} $ or they are both centered
ellipsoids in
 \cite{SSZ98}; one of $A$ and $B$  is a
symmetric  ellipsoid in Harg\'{e} \cite{Har99},  and etc. When one
of the symmetric sets is a slap, inequality (\ref{1}) is   usually
called
  \v{S}id\'{a}k-Khatri inequality. Some other proofs of this special case
  are given
  in   \cite{Sid68},  Jogdeo \cite{JOG},    \cite{DEOPSS72}
  with extension to elliptically  contoured distributions and    Szarek  and Werner \cite{SW99}
with  extension to an asymmetric case. See also the one-sided case
in Slepian \cite{SLE62}.  The result of \cite{Pitt77} is extended to
some multidimensional cases in Borell \cite{Bo81}. See also Figalli,
Maggi  and Pratelli \cite{FMP12}. For the case when one of the
symmetric convex set is an
 ellipsoid      in \cite{Har99}, another proof  can be found in         Cordero-Erausquin
\cite{CDE02}. See also  Lim Adrian   and Luo \cite{LS12} for a
special asymmetric case. Moreover, inequality $\mu_n(A\cap B)\geq
\mu_n(\lambda A)\mu_n(\sqrt{1-\lambda^2}B)$ ($0\leq \lambda \leq 1$)
is proved for $\lambda=2^{-1/2}$ in \cite{SSZ98}. The case for
general $\lambda$ is given
  in Theorem 1.1 Li \cite{Li99}.
For  applications of this inequality   on   Gaussian processes, we
refer to
 Li and Shao   \cite{LS01}.

  Since the family of symmetric convex sets  is stable under
linear transformation, inequality (\ref{1}) holds for general
centered Gaussian measures. Dimension free is another remarkable
property of (\ref{1}). In  Theorem \ref{C},
   inequality (\ref{1})
 is extended to   Wiener    space, which
  verifies   the
   conjecture formulated    in
\cite{LS01} and Lata{\l}a \cite{La} for instance.

Besides      the potential   applications of Theorem 1.1 to
Gaussian processes, in Theorem \ref{Ba2} we show
   that the Gaussian correlation inequality
 implies  the same  correlation inequality    for subordinate Brownian motion.
  Moreover, for any     open  sets $A,B\in  \mathcal{C}_n$ such that neither of them
  is equal to $\mathbb{R}^n$,
 we show that $\lambda_{A\cap B}\leq
\lambda_A+\lambda_B$ in Theorem \ref{La}. Here $\lambda_D$ is
denoted for  the spectral gap of Dirichlet Laplacian on a domain
$D$. The same  spectral gap inequality can also  be verified  for
  generators of subordinate Brownian motion with Dirichlet
boundary condition.

%for example,  $\mu_n(A\cap B)\leq \mu_n(\lambda
%A)\mu_n(\sqrt{1-\lambda^2}B)$($\lambda=1/2$ in   \cite{SSZ98} and
%$\lambda\in [0,1]$

\subsection{\normalsize introduction of the proof of Theorem 1.1}

 Next we introduce   the  proof of the conjecture    together with some related  methods
  used
   before. The proof given here is  based on some  previous    results
   about
   log-concave distributions   which will be  clear from the
   introduction in what below.
    Denote by $\langle \cdot,\cdot\rangle$ the
 standard inner product of $\mathbb{R}^n$.
For bounded  measurable functions $u$ and $v$ on  $\mathbb{R}^n$,
define
\begin{align}\label{defi3}\psi_\lambda(u,v )=\int
\int  u(x)v(y)f_{2n}(x,y;\lambda)dxdy,\ \ \ \ \forall\
\lambda\in[0,1),
\end{align}
where
\begin{align}\label{definition}
f_{2n}(x,y;\lambda)
=&\frac{1}{(2\pi)^n(1-\lambda^2)^{n/2}}\exp\{-\frac{|x|^2+|y|^2-2\lambda\langle
x,y\rangle}{2(1-\lambda^2)}\}, \ \ \ \ \ \ \forall\ x,y\in
\mathbb{R}^n .
\end{align}
 Denote
also
\begin{align}
\psi_1(u,v)=\int  uvd\mu_n.\nonumber
\end{align}

Let   $A,B\in \mathcal{C}_n$ in the rest part of this section. In
the study of   (\ref{1}), the following  relations are   used
frequently
\begin{align}\label{3'}
\mu_n(A\cap B)=& \psi_1(I_A,I_B),\ \ \ \ \ \mu_n(A) \mu_n(B)=
\psi_0(I_A,I_B),
\end{align}
which   can be verified   directly.    Therefore, to prove (\ref{1})
it is sufficient to show  that the derivative of $\psi_\lambda(A,B)$
is nonnegative.  This monotone property may be taken  as a finer
version of the  Gaussian correlation conjecture which  has been
verified  for  the case when one of the   symmetric convex sets is a
  slap in \cite{Sid68}\cite{JOG}  and also for   the
two dimensional case in \cite{Pitt77}.

The correlation  parameter $\lambda$ above  is   often used  in
statistical  literature. For function $\psi_\lambda$, another
parameter defined by $-\ln \lambda$ is  introduced  in
\cite{Pitt77}. Denote for every $t\geq 0$
\begin{align}\label{sem} \phi_t(u,v)= \psi_{e^{-t/2}}(u,v),\end{align}
when the right hand side above is well defined.  In what below
  $\psi_\lambda(I_A,I_B)$ and $\phi_t(I_A,I_B)$ are also
denoted by $\psi_\lambda(A,B)$ and $\phi_t(A,B)$, respectively. By
 (\ref{3'}) and (\ref{sem}),
\begin{align}\label{3''}
\mu_n(A\cap B)=&  \phi_0(A,B),\ \ \ \ \ \mu_n(A) \mu_n(B)=\lim_{t
\rightarrow \infty} \phi_t(A,B).
\end{align}

To study the monotone property of $\phi_t$, it is derived in
\cite{Pitt77} that, for   smooth functions   $u$ and $v$ with
  gradients  controlled by some polynomial for instance,
 \begin{align}
  \frac{d}{dt} \phi_t(u,v)=-\frac{1}{2}\int
\langle\nabla P_tu,\nabla v\rangle d\mu_n,\ \ \ \ \ \forall\  t\geq
0, \label{d}
\end{align}
where    $(P_t)$ is the  Ornstein-Uhlenbeck semigroup defined by
\begin{align}
P_t u(x)=&\frac{1}{(2\pi(1-e^{-t}))^{n/2}}\int u(
y)\exp\{-\frac{|y-e^{-t/2}x|^2}{2(1-e^{-t})}\}dy,\ \ \ \forall\
x\in
\mathbb{R}^n,\forall\ t>0;\label{P1} \\
 \ \ \ \ \ \ \  P_0u(x)=&u(x),\ \ \ \ \ \ \ \ \
 \ \ \ \ \ \ \ \ \ \ \ \ \ \ \ \ \ \  \ \ \ \ \ \ \ \ \ \ \ \ \ \ \ \ \ \
  \ \ \ \ \ \ \ \ \ \ \ \ \ \ \   \forall\  x\in
\mathbb{R}^n.\nonumber
\end{align}
One    way to derive  formula (\ref{d}) is  from    the fact that
the infinitesimal generator of $(P_t)$ is
$\frac{1}{2}(\Delta-\langle x,\nabla\rangle)$ and the following
relation
\begin{align}\label{semi} \phi_t(u,v) = \int(P_t
u)vd\mu_n,\ \ \  \  \forall\  t\geq 0.\end{align} The semigroup
point of view is adopted in \cite{Har99} on the conjecture.

It is  given  in \cite{SSZ98} that the following  inequality is
equivalent to   the conjecture:  for any $\varepsilon'>0$, there
exists some integer $N_0\geq1$ such that
 \begin{align}
\mu_n(A\cap B)\geq \exp\{- \varepsilon' n\}\mu_n(A)\mu_n(B),\ \ \ \
\ \forall\  A,B\in \mathcal{C}_n,\ \forall\ n\geq N_0.\label{2}
\end{align}
Instead of proving (\ref{1}) directly, our aim is to verify
(\ref{2}). To this end,  the parameter $\lambda$ and the parameter
$t$ are both  crucial in the proof. Next we give an outline of the
proof of (\ref{2})
     according to the correlation  of the assistant function   is strong, moderate and
     small, respectively.

\emph{Strong  correlation I} \ \  Let $\varepsilon>0$. The
derivative estimate of $\psi_\lambda(A,B)$,     given in Lemma
\ref{lower3},  shows that
 \begin{align}
\mu_n(A\cap B)=\psi_1(A,B)\geq & \exp\{- \varepsilon
n\}\psi_{1-\varepsilon}(A,B) ,\ \ \ \ \ \forall\  A,B\in
\mathcal{C}_n. \label{2'}
\end{align}    The proof of the   estimate (\ref{2'}) is based on
a   functional form  of the special case when one of the symmetric
convex sets is an    ellipsoid     mentioned above. See Lemma
\ref{Ho19}.

\emph{Strong  correlation II} \ \   We have
\begin{align}
 \psi_{1-\varepsilon}(A,B)
=& \psi_1(P_{ \delta}I_{A},P_{\delta}I_{B}) ,\ \ \ \ \ \forall\
A,B\in \mathcal{C}_n, \nonumber
\end{align}
where $\delta=-\ln (1-\varepsilon)$. The action
 of     $P_{ \delta}$    on $I_A$ and $I_B$
 allows  us to prove   the conjecture  under  an extra  assumption
that   both of  the    symmetric convex sets  contain large ball
with radius of order $\sqrt{n}$. See Corollary \ref{comb} for more
details.

 \emph{Small correlation I} \ \ Let $\alpha\in (0,1)$.
 For every
  $A\in \mathcal{C}_n$,  a  symmetric
log-concave function $h_{A,\alpha}(x)=\exp\{-H_{A,\alpha
}(x)\}=\exp\{-n\rho_A(x)-2^{-1}\alpha |x|^2\}$ on $\mathbb{R}^n$ is
defined in (\ref{HA'}). The estimate  for the increments of
$\phi_t(A,B)$   can then  be reduced to that of
$\phi_t(h_{A,\alpha},h_{B,\alpha})$ when $\alpha$ is close to zero.
We show that for some   $t_0>0$
\begin{align}
\frac{d}{dt}\phi_t(h_{A,\alpha},h_{B,\alpha}) <0 ,\ \ \ \ \ \forall\
t\in (t_0,\infty). \label{se}
\end{align}
The proof of (\ref{se}) is proved  by the first derivative estimate
and the second derivative   estimate of $\psi_\lambda$ at
$\lambda=0$ which are given in Lemma \ref{mid,a} and  Lemma
\ref{second} respectively.

\emph{Small  correlation II} \ \ In Lemma \ref{low},  the following
  inequality is given:
\begin{align}
\frac{d^2}{dt^2}\phi_t(h_{A,\alpha},h_{B,\alpha})
>& -\frac{1}{2}
\frac{d}{dt}\phi_t(h_{A,\alpha},h_{B,\alpha}) ,\ \ \ \ \ \forall\
t\in (T(\alpha),\infty), \label{A'B'}
\end{align}
where $T(\alpha)$ is a positive constant depending on $\alpha$. The
relation above  relies on the following uniform estimate
\begin{align}
 C(\alpha) e^{-t}I_n\leq \nabla^2 H_{A,\alpha,t}(x) \leq &2(1\wedge t)^{-1}e^{-t}I_n ,\ \ \ \ \forall\  x\in
 \mathbb{R}^n,\ \forall\  t>0,\label{un}
\end{align}    where $H_{A,\alpha,t}$ is defined by
$P_th_{A,\alpha}=\exp\{-H_{A,\alpha,t}\}$ and
$C(\alpha)=\min(e^{-3}\alpha, 2^{-6}e^{-3})$. The estimate
(\ref{un}) is given in Lemma
  \ref{C1,2,3;}. Combing (\ref{se}) and (\ref{A'B'}), we
get \begin{align} \frac{d}{dt}\phi_t(h_{A,\alpha},h_{B,\alpha}) <0,\
\ \ \ \ \ \ \ \ \
 \forall t\in (T(\alpha),\infty).\nonumber\end{align}
The estimate  (\ref{un}) is  a quantity version of the fact that
$P_t u\in \mathcal{CF}_n$ if $u\in \mathcal{CF}_n$,
 which is a consequence of Theorem 7
in Pr\'{e}kop \cite{Pre73}.

 \emph{Moderate  correlation I} \ \
 As the estimate (\ref{2'}), the increments of function
$\phi_t(h_{A,\alpha},h_{B,\alpha})$ on $ [0,\varepsilon)$ can be
controlled well     in order to verify (\ref{2}).

 \emph{Moderate  correlation II} \ \ The conclusion  in the  strong correlation II  above allows us  further  assume that
 $B_n(\delta \sqrt{n}) \subseteq A\cap B$ for some    $\delta>0$.
 In order to estimate the increments of $ \phi_t(h_{A,\alpha},h_{B,\alpha})$
 for  $t\in [\varepsilon,T(\alpha)]$,  we show in Lemma \ref{upper} that
  for $n$ big enough depending on $\varepsilon,\delta$ and $\alpha$
\begin{align}
 \frac{d}{dt}\phi_t(h_{A,\alpha},h_{B,\alpha})<
 \varepsilon n
 \phi_t(h_{A,\alpha},h_{B,\alpha})
,\ \ \ \ \ \forall\  t\in [\varepsilon, T(\alpha)] . \nonumber
\end{align}
Since the order of $T(\alpha)$   is $-\ln \varepsilon$ when taking
$\alpha=\varepsilon$ in the final proof of Theorem 1.1,  the
inequality (\ref{2}) can be verified from  the estimates introduced
all above.  The proof of Lemma \ref{upper} is based a monotone
estimate for some assistant function with dilation parameter;  see
Lemma \ref{upper,pp}. The main
   tools in the proof of Lemma \ref{upper,pp} and some
other related estimates to  prove   Lemma \ref{upper} are some
concentration inequalities  for   certain log-concave distributions,
including   the Poincar\'{e} inequality given in Brascamp and Lieb
\cite{BL76}
  and the isoperimetric inequality given in Bakry and Ledoux \cite{BL96}. We refer to
   Ledoux \cite{Le99}\cite{Le01}  for    more information of  this subject.

\subsection{\normalsize structure of the paper with   some further comments and notations}

The rest part of the paper is organized as follows. The second
section provides some basic estimates of $
{d\psi_\lambda}/{d\lambda}$ together with a formula for the second
derivative. In Lemma \ref{zz} we show that the monotone property of
$\psi_{\lambda}$ can be obtained when the Ornstein-Uhlenbeck
semigroup is replaced by the semigroup of Brownian motion together
with a change of   reference measure.

The first part of     section three  is to reduce  the condition of
the  conjecture to the case that the symmetric  convex sets
containing  large ball.  Most effort of this section is to prove
Proposition \ref{C1,2,3}. We remark that the result  of Proposition
\ref{C1,2,3}
 is motivated  by  Pr\'{e}kopa's result in \cite{Pre73}, however, the proof
 can be modified by induction without applying this result.
  The
last part of this section is to introduce   some log-concave
functions associated with symmetric convex sets and   prove  the
uniform estimate in (\ref{un}) above.

To study the moderate correlation part   introduced above, in
  section four, we  give  some  basic derivative estimates of
$\phi_t$ for the  associated log-concave functions. Some related
  formulas for the derivatives can be found  in, e.g., \cite{DEOPSS72},
\cite{Hu97}, Houdr\'{e},
 P\'{e}rez-Abreu  and   Surgailis  \cite{HPS98}, \cite{LS01} and
 Harg\'{e}
\cite{Har05}.
  The relation  of  $\phi_t$ for  the symmetric
convex sets and the associated log-concave functions is   given in
section five.  The proof of Theorem \ref{A} is also given in section
five. When the convex sets are bounded,  we note that  the result of
Lemma \ref{second} has been given in Koldobsky
  and   Montgomery-Smith  \cite{KM96} by  the method of  Fourier analysis.
  Therefore, the estimate (\ref{se}) above is essentially given in
  \cite{KM96}. For the unbounded case,  we use  Anderson inequality in Anderson
\cite{An55}  in the proof of  Lemma \ref{second}.
     Some applications
of Theorem 1.1 mentioned above  are given in the last section.

Next we briefly introduce  another two   types  of functional
correlation inequalities for Gaussian measure which are  closely
related to the paper. One is a correlation inequality for two convex
functions under $\mu_n$ given in Hu \cite{Hu97}. The other is a
  correlation inequality for a  convex function and
a log-concave function   under general Gaussian measures   given in
  Harg\'{e} \cite{Har04}. The  relation   of these  two
types of     inequalities can be found in  \cite{Har04}. See also
Remark \ref{com} below for some further comments.

The derivative estimate given in the second section can also   be
proved by Harg\'{e}'s inequality in \cite{Har04} mentioned above.
The proof of Harg\'{e}'s inequality in \cite{Har04} relies  on a
sharp regularity estimate for certain Brenier  map given in
Caffarelli \cite{Caff2000}. Since we only need some special cases of
Harg\'{e}'s inequality,   this approach is not adopted here.
Moreover,   Lemma \ref{lower3}  can also be proved by
\v{S}id\'{a}k-Khatri inequality since the convex function used in
the proof of Lemma \ref{lower3} is   square function.

Throughout the paper, notations $m,n,i,j,k,l$ are always denoted for
    integers  with $n\geq 1$,  notations
   $x,y,z$ are always denoted for
 elements  of $ \mathbb{R}^n$.  The coordinates of
  $x$    are denoted by    $(x_1,\cdots,x_n)$ and  the same  convention
 is  applied also to $y$ and $z$.
 For  $A \subseteq \mathbb{R}^n$, denote  $F(A)=\{F(x):x\in A\}$
for a map $F$ defined on $\mathbb{R}^n$ and  denote by $A^c$ the
complement  of $A$. The notation of the integral
 $\int$
    is over all of $\mathbb{R}^n$ unless
explicitly stated otherwise. Denote by $f_n(\cdot)$   the density
function of the standard Gaussian measure   on $\mathbb{R}^n$.
Denote   $S_{n-1}=\{ x\in \mathbb{R}^n:|x|=1 \}$ and $B_n(r)=\{x\in
\mathbb{R}^n:|x|< r\}$ for every  $r\geq 0$.
 Denote by $m_k$ the $k-$dimensional Hausdorff
measure for every  $k\geq 0$. Denote by $I_n$ the identity matrix on
$\mathbb{R}^n$. For two
  $n$ by $n$ matrixes $Q_1$ and $Q_2$, denote  $Q_1\leq Q_2 $ when
$Q_2-Q_1$ is a nonnegative definite matrix.

 A
nonnegative function $f$ on $\mathbb{R}^n$ is called  log-concave if
$f(\lambda x+(1-\lambda)y)\geq f(x)^\lambda f (y)^{1-\lambda}$ holds
for every  $x,y\in \mathbb{R}^n$ and every  $0<\lambda <1$. Denote
by $\mathcal{CF}_n$ the class of symmetric log-concave functions on
$\mathbb{R}^n$.
 For an unit vector $\mathbf{e}\in
\mathbb{R}^n$, denote by $\partial_{\mathbf{e}}$ the partial
derivative   along    $\mathbf{e}$. For $1\leq i\leq n$, denote by
$\mathbf{e}_i$ the unit vector of $\mathbb{R}^n$  of which the
$i$-th coordinate is equal to one. For $1\leq i\leq n$, denote
$\partial_{ \mathbf{e}_i}$ by $\partial_i$. For twice differentiable
function $f$ on $\mathbb{R}^n$, denote by $\nabla^2f$ the Hessian of
$f$.  For $a>0$ denote by $\lfloor a \rfloor$ the integer part of
$a$.  For $a,b\in \mathbb{R}$, denote $a\vee b=\max\{a,b\}$ and
$a\wedge b=\min\{a,b\}$.  Some other notations will be introduced in
what below when necessary.

\section{Derivative estimates for   correlation parameter
   }
 }

\subsection{\normalsize  Harg\'{e}'s   correlation inequality}
 The main aim of this section is to give some derivative estimate of
 $\psi_\lambda(A,B)$ which is helpful   in particular when $\lambda$ is close to one or zero
 in the proof of  Theorem 1.1.

In what below we say that  a set $A\subseteq \mathbb{R}^n$ is a
centered ellipsoid if $A=\{x: |\langle \Sigma x,x\rangle|\leq 1\}$
for some symmetric nonnegative definite   matrix $\Sigma$. We say
that the lower   level sets of
  a nonnegative function $f$ on $\mathbb{R}^n$ are
centered ellipsoids if    $\{x:f(x)\leq r\}$ is either  a centered
ellipsoid or a degenerate set for every   $r\geq 0$.

\begin{lemma}\label{B} \emph{[Corollary 3 in \cite{Har99}]}\ Let $A\in \mathcal{C}_n$
and $\gamma_n$  be a centered Gaussian measure on $\mathbb{R}^n$.
Then for every centered ellipsoid $B \subseteq\mathbb{R}^n$
  \begin{align}
\gamma_n(A\cap B)\geq \gamma_n(A)\gamma_n(B).\nonumber
\end{align}
\end{lemma}

For function $u$ on $\mathbb{R}^n$, denote by $\mathrm{Supp}(u)$ the
closure of the set $\{x:u(x)>0\}$ in $\mathbb{R}^n$. The following
inequality is motivated by  Harg\'{e}'s   correlation inequality in
\cite{Har04}  and it is   a functional form of the result above. See
also Theorem 2 in \cite{Har99} for another form. Notice that
(\ref{21}) below still holds when the left hand side is infinity.
 \begin{lemma} \label{Ho19} Let   $\gamma_n$ be a
centered  Gaussian measure on $\mathbb{R}^n$. Let $u\in
\mathcal{CF}_n$ and    $f$ be a  nonnegative function of which the
lower  level sets are centered ellipsoids. Then
   \begin{align}
\int  fu d\gamma_n\leq  \int u  d\gamma_n \int fd\gamma_n,
\label{21}
\end{align}
provided that both sides above are well defined finite  integrals.
Let $M\in(0,\infty)$. If further assuming that $0\leq f(x)\leq M$
for every  $x\in \mathrm{Supp}(u)$, then
     \begin{align}\label{qu}
  \int  fu d\gamma_n\leq  \int u  d\gamma_n \int  (f \wedge M) d\gamma_n.
\end{align}
\end{lemma}
\noindent{\bf Proof}\  We have by Fubini theorem,  the assumptions
of $u,f$ and Lemma \ref{B}
\begin{align}
 \int f u  d\gamma_n
 =& \int \big(\int_0^\infty I_{f> r}dr\big)  u d\gamma_n\nonumber\\
  =& \int_0^\infty\big(\int u d\gamma_n- \int  uI_{f\leq  r}
  d\gamma_n\big)dr\nonumber\\
 \leq &\int_0^{\infty}\big(\int u  d\gamma_n- \int u  d\gamma_n\int
I_{f\leq r}  d\gamma_n
      \big)dr\nonumber\\
      =& \int u  d\gamma_n \int fd\gamma_n,\nonumber
  \end{align}
which gives  (\ref{21}).  Suppose further that $0\leq f(x)\leq M$
for all $x\in \mathrm{Supp}(u)$. Then  we have $\int fu d\gamma_n=
\int (f\wedge M)u d\gamma_n$. Noticing that the lower level sets of
$f\wedge M$ are centered ellipsoids, we get the second conclusion by
(\ref{21}). \qed
\medskip

\subsection{\normalsize some derivative estimates}

Recall that    $f_{2n}(x,y;\lambda)$ is defined by
(\ref{definition}) and $f_n$ is the density function of $\mu_n$. We
refer to \cite{An58} for some basic properties of Gaussian measure.
When $\lambda=0$,
\begin{align}\label{inde}
f_{2n}(x,y;0)=&f_n(x)f_n(y), \ \ \ \forall\  x,y\in \mathbb{R}^n.
\end{align}
For every  $\lambda\in[0,1)$,
\begin{align}
 \int\int \langle x,y\rangle
f_{2n}(x,y;\lambda)dxdy =&\lambda n,\label{e2q}\\\label{e1q}
 \int\int|x|^2
f_{2n}(x,y;\lambda)dxdy=&  n.\end{align}
  \begin{lemma}\label{22}
Let $u$ and $v $ be two bounded measurable functions on
$\mathbb{R}^n$. Then for every   $\lambda\in [0,1)$
  \begin{align}\label{de'}
&\frac{d \psi_{\lambda}(u,v)}{ d\lambda}=\int\int
h_\lambda(x,y)u(x)v(y)f_{2n}(x,y;\lambda)dxdy,
 \end{align}where
 \begin{align}  h_\lambda(x,y)= \frac{
-\lambda(|x|^2+|y|^2)+(1+\lambda^2)\langle x,
y\rangle+n\lambda(1-\lambda^2)}{(1-\lambda^2)^2}.\label{hl}
\end{align}

\end{lemma}\noindent\textbf{Proof}\
  For every $x,y\in \mathbb{R}^n$ and every
$\lambda\in [0,1)$, we have by (\ref{definition})
\begin{align}
\frac{\partial f_{2n}(x,y;\lambda)}{ \partial \lambda}
=&\frac{\partial }{
\partial \lambda }\Big(\frac{1}{ (2\pi)^n(1-\lambda^2)^{n/2}}
 {\exp{\{-  \frac{ |x|^2+|y|^2-2\lambda\langle x,y\rangle}{2(1-\lambda^2)} \}}}\Big)\nonumber\\
=& \frac{ -\lambda(|x|^2+|y|^2)+(1+\lambda^2)\langle x,
y\rangle+n\lambda(1-\lambda^2)}{ (1-\lambda^2)^2}
f_{2n}(x,y;\lambda) .\nonumber
 \end{align}
Then we get  (\ref{de'})   by definition (\ref{defi3}).\qed\medskip

\begin{lemma}\label{mid,a}
Let $u$ and $v $ be two bounded measurable  functions on
$\mathbb{R}^n$  and assume that  $v$ is symmetric. Then
    \begin{align}
\big(\frac{d}{d\lambda}\psi_{\lambda}(u,v)\big)_{\lambda=0}
=0.\nonumber
    \end{align}

\end{lemma}\noindent\textbf{Proof}\
By (\ref{inde}), (\ref{de'}) and Fubini theorem,
    \begin{align}
\big(\frac{d}{d\lambda}\psi_{\lambda}(u,v)\big)_{\lambda=0}=&\int
\int \langle x,y\rangle u(x) v(y)
d\mu_{n}(x)d\mu_{n}(y)\nonumber\\=&\sum_{i=1}^n\int x_iu(x)d\mu_n(x)
\int y_iv(y) d\mu_{n}(y)=0,\nonumber
    \end{align}
where we use     assumption   $v(y)=v(-y)$ in the last equality
above.
\medskip\qed

\begin{lemma}\label{mix}
Let $u$ and $v $ be two  measurable functions on $\mathbb{R}^n$
which are both bounded and  nonnegative. Suppose further that $v$ is
symmetric. Then for every  $\lambda\in [0,1)$
    \begin{align}
\int\int  \langle x, y\rangle u(x)v(y)f_{2n}(x,y;\lambda) dxdy\geq
0.\nonumber
    \end{align}

\end{lemma}\noindent\textbf{Proof}\
Let $\lambda\in [0,1)$.  To prove the lemma, by Fubini theorem and
the assumption $u\geq 0$, it is sufficient to verify that   for
every $x\in \mathbb{R}^n$
    \begin{align}\label{first po2}
   \langle x,\int  y  v(y)f_{2n}(x,y;\lambda)dy \rangle \geq 0.
    \end{align}
From the   assumption of $v$,  function  $v(Q(\cdot))$ is
  symmetric   for any orthogonal transformation
$Q$ of $\mathbb{R}^n$.  We also have  that
   $ f_{2n}(\cdot,\cdot;\lambda)=f_n(Q(\cdot),Q(\cdot);\lambda)$  for any orthogonal transformation
$Q$ of $\mathbb{R}^n$. Therefore,
     to prove (\ref{first po2}) we can assume   in what below
     that $x=a\textbf{e}_1$
    for some  $a\geq 0$.  Here $\textbf{e}_1=(1,0,\cdots, 0)$. In other
    words,
    to prove (\ref{first po2}) it is sufficient to show that
     \begin{align}\label{first po3}
  \int  y_1v(y) f_{2n}(a\textbf{e}_1,y;\lambda)dy  \geq 0.
    \end{align}
The assumption  $a\geq 0$ implies that $|y-a\lambda\textbf{e}_1|\leq
|y+a\lambda\textbf{e}_1|$
    when  $y_1\geq 0$. Therefore,  by the symmetric and nonnegative  assumptions of $v$, we
    get
         \begin{align}
 & \int  y_1 v(y)f_{2n}(a\textbf{e}_1,y;\lambda)dy\nonumber\\
  =&\frac{1}{(2\pi)^n(1-\lambda^2)^{n/2}}\int y_1v(y)
  \exp\{-\frac{|y-\lambda a\textbf{e}_1|^2}{2(1-\lambda^2)}\}dy\nonumber\\
   =&\frac{1}{(2\pi)^n(1-\lambda^2)^{n/2}}\int_{y_1\geq 0} y_1
   v(y)
  \Big(\exp\{-\frac{|y-\lambda a\textbf{e}_1|^2}{2(1-\lambda^2)}\}-\exp\{-\frac{|y+\lambda a\textbf{e}_1|^2}
  {2(1-\lambda^2)}\}\Big)dy\nonumber\\
  \geq&0,\nonumber
    \end{align}which gives (\ref{first po3}).
\medskip\qed

  \begin{lemma}\label{lower3}
Let $u,v\in \mathcal{CF}_n$. Then for every    $\lambda\in [0,1)$
  \begin{align}\label{whole3}
\frac{d \psi_{\lambda}(u,v)}{ d\lambda} \geq \frac{-\lambda n
}{(1+\lambda)^2}\psi_{\lambda}(u,v).
 \end{align}

\end{lemma}\noindent\textbf{Proof}\
Let $\lambda\in [0,1)$. Set
\begin{align}\widetilde{h}_\lambda(x,y):=\lambda(|x|^2+|y|^2)-2\lambda\langle
x, y\rangle,\ \ \  \ \forall\  x,y\in
\mathbb{R}^n.\nonumber\end{align}
 Notice that the lower level sets of $\widetilde{h}_\lambda$ are ellipsoids of
$\mathbb{R}^{2n}$. By   (\ref{e2q}) and (\ref{e1q}),
\begin{align}\int\int\widetilde{h}_\lambda(x,y) f_{2n}(x,y;\lambda)dxdy=
2\lambda(1-\lambda)n.\nonumber
\end{align}
This and (\ref{21}) give
\begin{align}
\int\int
\widetilde{h}_\lambda(x,y)u(x)v(y)f_{2n}(x,y;\lambda)dxdy\leq
2\lambda(1-\lambda)n\int\int
 u(x)v(y)f_{2n}(x,y;\lambda)dxdy. \nonumber
\end{align}
Applying   Lemma \ref{22}, Lemma  \ref{mix} and the inequality
above, we have
   \begin{align}
\frac{d \psi_{\lambda}(u,v)}{ d\lambda} =&\int\int \frac{
-\lambda(|x|^2+|y|^2)+(1+\lambda^2)\langle x,
y\rangle+\lambda(1-\lambda^2)n}{(1-\lambda^2)^2}u(x)v(y)f_{2n}(x,y;\lambda)dxdy\nonumber\\
=&\int\int \frac{ -\widetilde{h}_\lambda(x,y)+(1-\lambda)^2\langle
x,
y\rangle+\lambda(1-\lambda^2)n}{(1-\lambda^2)^2}u(x)v(y)f_{2n}(x,y;\lambda)dxdy\nonumber\\
\geq &\int\int \frac{
 -2\lambda(1-\lambda)n+\lambda(1-\lambda^2)n }{(1-\lambda^2)^2}u(x)v(y)f_{2n}(x,y;\lambda)dxdy\nonumber\\
= &\frac{
 -\lambda n }{ (1+\lambda)^2}\int\int
 u(x)v(y)f_{2n}(x,y;\lambda)dxdy\nonumber,
 \end{align}which gives (\ref{whole3}).\qed \medskip
\begin{remark}\label{re2}
The  lower bound of the derivative  given above implies that
$\mu_n(A\cap B)\geq \exp\{-(\ln 2-2^{-1})n\}\mu_n(A)\mu_n(B)$ for
every  $A,B\in \mathcal{C}_n$. This  improves the estimate
$\mu_n(A\cap B)\geq \exp\{- \frac{n}{2} \}\mu_n(A)\mu_n(B)$ given in
\cite{SSZ98}. See also \cite{Sh03} for another form of  estimate.
\end{remark}

\begin{lemma}\label{SD}
Let $u$ and $v $ be two bounded measurable functions on
$\mathbb{R}^n$. Then for every  $\lambda\in [0,1)$
  \begin{align}
&\frac{d^2 \psi_{\lambda}(u,v)}{ d\lambda^2} \nonumber\\=&\int \Big(
h_\lambda(x,y)^2+\frac{ -(1+3\lambda^2)(|x|^2+|y|^2)+
2\lambda(3+\lambda^2)\langle x, y\rangle+n(1-
\lambda^4)}{(1-\lambda^2)^3}\Big)u(x)v(y)f_{2n}(x,y;\lambda)dxdy,\nonumber
 \end{align}
where $h_\lambda$ is defined by (\ref{hl}).
\end{lemma}\noindent\textbf{Proof}\
For every $x,y\in \mathbb{R}^n$, direct calculation shows that
\begin{align}
\frac{ \partial  }{ \partial\lambda}h_\lambda (x,y)= &\frac{
-(|x|^2+|y|^2)+2\lambda\langle x,
y\rangle+n(1-3\lambda^2)}{(1-\lambda^2)^2}\nonumber\\
+&\frac{ -4\lambda^2(|x|^2+|y|^2)+4\lambda(1+\lambda^2)\langle x,
y\rangle+4n\lambda^2(1-\lambda^2)}{(1-\lambda^2)^3}\nonumber\\
=&\frac{ -(1+3\lambda^2)(|x|^2+|y|^2)+ 2\lambda(3+\lambda^2)\langle
x, y\rangle+n(1- \lambda^4)}{(1-\lambda^2)^3}. \nonumber
\end{align}
Then  we get the conclusion   by (\ref{de'}).
\medskip\qed

\subsection{\normalsize a  correlation  inequality associated with Brownian motion}

We know    that
  $d\psi_{\lambda}(u,v)/d\lambda\geq 0$  holds for  every   $\lambda\in[0,1)$
  and every
  $u,v\in \mathcal{CF}_n$
   if and only if
     for every    smooth
  functions
  $u,v\in \mathcal{CF}_n$ with bounded supports
\begin{align}
  &\int \langle\nabla u,\nabla v\rangle  d\mu_n \geq 0. \nonumber
\end{align}
  The following lemma shows that
the inequality above holds when  the Gaussian measure is replaced by
the Lebesgue measure.

\begin{lemma}\label{zz} Let $u$ and $v $ be two smooth functions  of  $\mathcal{CF}_n $  with bounded supports.  Then
\begin{align}
  &\int \langle\nabla u,\nabla v\rangle  dx \geq 0. \nonumber
\end{align}
\end{lemma}
\noindent {Proof }\   Let $(T_t)$ be the semigroup associated with
the standard Brownian motion on $\mathbb{R}^n$, which means that for
any  bounded measurable function $f$
\begin{align}
T_tf(x)=\frac{1}{(2\pi t)^{n/2}}\int
\exp\{-\frac{|y-x|^2}{2t}\}f(y)dy,\ \ \ \forall\ x\in\mathbb{R}^n,
\forall \ t>0;\ \ \ \ \ T_0f=f.\label{BM}
\end{align}
For smooth functions  $f$ and $g$      with bounded supports, we
know that
\begin{align}
\Big(\frac{d}{dt}\int fT_tgdx\Big)_{t=0}=-\frac{1}{2}\int
\langle\nabla f,\nabla g\rangle dx. \label{BM2}
\end{align}
Let $u$ and $v $ be two smooth functions  of $\mathcal{CF}_n $  with
bounded supports.  By (\ref{BM}),  (\ref{BM2}) and approximation, to
prove the lemma it is sufficient to show that for every  $t>0$
\begin{align}  \frac{d}{dt}\int uT_tvdx  \leq 0.\label{ttt}
\end{align}

Let  $t>0$ in what below. We have by (\ref{BM})
\begin{align}  \frac{d}{dt}\int uT_tvdx
=&\frac{d}{dt}\Big(\frac{1}{(2\pi t)^{n/2}}\int u(x)dx\int
\exp\{-\frac{|y-x|^2}{2t}\}v(y)dy\Big)\nonumber\\
=&-\frac{1}{2t^2(2\pi t)^{n/2}}\int u(x)dx\int
\exp\{-\frac{|y-x|^2}{2t}\}(nt-|x-y|^2)v(y)dy.\label{t}
\end{align}
Next we apply a standard  technique of   approximation, c.f.
\cite{BA95}.  We have
\begin{align}
 &\frac{1}{ (2\pi t)^{n/2}}\int u(x)dx \int
\exp\{-\frac{|y-x|^2}{2t}\} |x-y|^2 v(y)dy\nonumber\\
=& \lim_{\delta \rightarrow \infty}\frac{1}{ (2\pi
t)^{n/2}}\int\exp\{-\frac{|x|^2}{2\delta}\}u(x)dx\int
\exp\{-\frac{|y-x|^2}{2t}\} |x-y|^2 v(y)dy. \label{del}
\end{align}
For every $\delta>0$, define Gaussian measure $\gamma_{2n}(x,y;
\delta,t)$ on $\mathbb{R}^{2n}$ by
$$ d\gamma_{2n}(x,y;t,\delta)=\frac{1}{
(2\pi)^{n}(t\delta)^{n/2}} \exp\{-\frac{|x|^2}{2\delta}
-\frac{|y-x|^2}{2t}\}dxdy.$$ Then we have by  (\ref{21})
\begin{align}
 &\frac{1}{ (2\pi
t)^{n/2}}\int\exp\{-\frac{|x|^2}{2\delta}\}u(x)dx\int
\exp\{-\frac{|y-x|^2}{2t}\} |x-y|^2 v(y)dy\nonumber\\
=& (2\pi \delta)^{n/2}\int \int|x-y|^2 u(x)
v(y)d\gamma_{2n}(x,y;\delta,t) \nonumber\\
\leq & (2\pi \delta)^{n/2}\int\int  u(x)
v(y)d\gamma_{2n}(x,y;\delta,t)  \int  \int|x-y|^2 d\gamma_{2n}(x,y;\delta,t)\nonumber\\
=&\frac{nt}{ (2\pi
t)^{n/2}}\int\exp\{-\frac{|x|^2}{2\delta}\}u(x)dx\int
\exp\{-\frac{|y-x|^2}{2t}\}   v(y)dy,\nonumber
\end{align}
which gives together with (\ref{del})
\begin{align}
 &\frac{1}{ (2\pi t)^{n/2}}\int u(x)dx\int
\exp\{-\frac{|y-x|^2}{2t}\} |x-y|^2 v(y)dy\nonumber\\
\leq& \frac{nt}{ (2\pi t)^{n/2}}\int u(x)dx\int
\exp\{-\frac{|y-x|^2}{2t}\} v(y)dy. \label{del2}
\end{align}
Combing (\ref{t}) and (\ref{del2}), we get (\ref{ttt}). \qed
\medskip

\begin{remark}\label{z}
If for any smooth functions $u$ and $v $     of $\mathcal{CF}_n $
with bounded supports
\begin{align}
  &\int_{B_{n}(r)} \langle\nabla u,\nabla v\rangle  dx \geq 0,\ \ \ \forall\    r>0, \nonumber
\end{align}
then  the first conclusion of Theorem 1.1  can be verified  by
Fubini theorem.
 From the proof in
\cite{Pitt77},   the inequality above may  hold  when  $n=2$.
   \end{remark}

\section{From symmetric convex sets to certain      log-concave functions  }

\subsection{\normalsize  symmetric convex sets containing large ball  }

Recall that  $(P_t)_{t\geq0}$ is   the Ornstein-Uhlenbeck semigroup
on $\mathbb{R}^n$. We know that $(P_t)_{t\geq0}$  is a symmetric
semigroup on $L^2(\mathbb{R}^n, \mu_n)$ which gives  that
\begin{align}
  \phi_t(u,v) = &\phi_t( v,u)
 ,\ \ \ \ \ \forall\  t\geq0, \label{t,t}\\     \phi_{t+s}(u,v)
  =& \int (P_tu)
 P_sv d \mu_n,\ \ \ \ \forall\  t,s\geq 0.
 \label{t/2}\end{align}
where  $u,v\in L^2(\mathbb{R}^n, \mu_n)$.

 The following Lemma  is frequently used  in
the study of the conjecture. It is a direct consequence of  Theorem
7 in \cite{Pre73} and the fact that the density function of Gaussian
measure is log-concave.

 \begin{lemma}   \label{Pre}
Suppose that  $f\in  \mathcal{CF}_n$. Then $P_tf\in  \mathcal{CF}_n$
for every $t\geq 0$.
\end{lemma}

 \begin{lemma} \label{db0}Let $r_0,t_0\in [0,1]$ and $A \in \mathcal{C}_n$.    Then
for every  $x\in B_n(r_0 \sqrt{n})$
\begin{align}
  \exp\{  - t_0^{-1}r_0^2n \}P_{t_0} I_A(0)\leq   P_{t_0} I_A(x)\leq    P_{t_0} I_A(0).\label{db}
\end{align}
 \end{lemma}
 \noindent{\bf Proof}\ Let $r_0, t_0\in [0,1]$ and
   $x\in B_n(r_0 \sqrt{n})$. Then we have by (\ref{P1})
 \begin{align}
 P_{t_0} I_A(x)=&\frac{1}{(2\pi(1-e^{-{t_0}}))^{n/2}}\int_A
 \exp\{-\frac{|y-e^{-{t_0}/2}x|^2}{2(1-e^{-{t_0}})}\} dy\nonumber\\
\geq &\frac{1}{(2\pi(1-e^{-{t_0}}))^{n/2}}\int_A
 \exp\{-\frac{|y |^2 +r_0^2n-2e^{-{t_0}/2} \langle y ,
 x\rangle}{2(1-e^{-{t_0}})}\}dy.
  \nonumber
\end{align}
From the symmetric assumption of $A$ and the convexity of exponent
function, we further get
 \begin{align}
 P_{t_0} I_A(x)\geq &
 \frac{ 1}{2(2\pi(1-e^{-{t_0}}))^{n/2}}\exp\{  -  \frac{r_0^2n}{2(1-e^{-t_0})} \}\int_A\Big(
 \exp\{-\frac{|y |^2 -2e^{-{t_0}/2} \langle y ,
 x\rangle}{2(1-e^{-{t_0}})}\}\nonumber\\
 &\ \ \ \ \ \ \ \ \  \ \ \  \ \ \ \ \ \  \ \ \  \ \ \ \ \ \  \ \ \  \ \ \ \ \ \
  \ \ +\exp\{-\frac{|y |^2 +2e^{-{t_0}/2} \langle y ,
 x\rangle}{2(1-e^{-{t_0}})}\}\Big)dy
  \nonumber\\
  \geq &
 \frac{1}{(2\pi(1-e^{-{t_0}}))^{n/2}} \exp\{  -  \frac{r_0^2n}{2(1-e^{-t_0})} \}\int_A
 \exp\{-\frac{|y |^2 }{2(1-e^{-{t_0}})}\} dy
  \nonumber\\
  = & \exp\{  -  \frac{r_0^2n}{2(1-e^{-t_0})} \}P_{t_0} I_A(0).\nonumber
\end{align}
Noticing that $1-e^{-{t}}>{t}/2$ when  ${t}\in [0,1]$, we get the
first inequality of (\ref{db}) from the estimate above. By  Lemma
\ref{Pre}
 and $A\in \mathcal{C}_n$, we have
$P_{t_0}I_A\in \mathcal{CF}_n$. Therefore, the function $P_{t_0}I_A$
takes its maximal at the origin which gives to the second inequality
of (\ref{db}). \qed
\medskip

Set for every  $ A\in \mathcal{C}_n  $ and every $a,r\geq 0$
\begin{align} A_{a,r}=&\{x: P_aI_A(x)\geq r P_aI_A(0)\}.\label{AR}
 \end{align}

 \begin{lemma} \label{ke}Let  $ A, B \in \mathcal{C}_n$ and  $0\leq t_1\leq
 t_2$.
   Suppose that  for some constants $a>0$ and $c_0\in (0,1]$
 \begin{align}
   \phi_{t_1}(A_{a,r}, B_{a,s})\geq  c_0
   \phi_{(t_2-2a)\wedge t_1}(A_{a,r}, B_{a,s}),\ \ \ \forall \ r,s\in [0,e^{-an} ].\label{db1}
\end{align}
Then
\begin{align}
   \phi_{t_1}(A, B)\geq  c_0 e^{-3an}
 \phi_{t_2}(A, B).\label{db2}
\end{align}
 \end{lemma}
 \noindent{\bf Proof}\
Let  $t>0$ and  set    $ \lambda=e^{-t/2}$. By Lemma \ref{lower3}
\begin{align}
\frac{d\phi_{t}(A,B) }{dt}
=&-\frac{1}{2}e^{-t/2}\frac{d\psi_{\lambda}(A,B) }{d\lambda}  \leq
  \frac{1}{2}e^{-t/2}
\frac{n\lambda}{(1+\lambda)^2}\psi_{\lambda}(A,B) \leq
\frac{n}{2}\phi_{t}(A,B) ,\nonumber
\end{align}
which gives
\begin{align}
\phi_{t'}(A,B) \leq e^{(t'-t)n/2} \phi_{t}(A,B) ,\ \ \ if \ 0\leq
t\leq t'.\nonumber
\end{align}
Therefore,
 \begin{align}
 \phi_{t_1}(A,  B)\geq e^{-an}\phi_{t_1+2a}(A,B).\label{in}
 \end{align}
The estimate above shows that   the lemma  holds when  $t_2\leq
t_1+2a$. Next we
 assume that $t_2>t_1+2a$ and set $t_2'=t_2-2a$.
By  Fubini theorem and (\ref{t/2}),
 \begin{align}
  \phi_{t_1+2a}(A,B)=&\phi_{t_1}(P_aI_A,P_aI_B)\nonumber\\
  =&\int\int P_aI_A(x)P_aI_B(y)f_{2n}(x,y;e^{-t_1/2})dxdy
 \nonumber\\
      =& \int_0^{\infty}dr\int_0^{\infty}ds\int\int I_{P_a I_A(x) \geq
 r}I_{P_a I_A(y) \geq
 s} f_{2n}(x,y;e^{-t_1/2})dxdy.\nonumber
 \end{align}
 Then, by the second inequality of (\ref{db}) and definition (\ref{AR})
 \begin{align}
  \phi_{t_1+2a}(A,B)=&\int_0^{ P_aI_A(0)} dr\int_0^{ P_aI_B(0)}ds
  \int   \int   I_{P_aI_A(x) \geq
 r} I_{P_aI_B(y) \geq
 s}f_{2n}(x,y;e^{-t_1/2})dxdy
 \nonumber\\=& \int_0^{ P_aI_A(0)} dr\int_0^{ P_aI_B(0)}
\phi_{t_1}(A_{a,r/P_aI_A(0)}, B_{a,s/P_aI_B(0)}   )ds
 \nonumber\\
= &P_aI_A(0)P_aI_B(0)\int_0^{1} dr\int_0^{1}
 \phi_{t_1}(A_{a,r}, B_{a,s}) ds\nonumber\\
\geq  &P_aI_A(0)P_aI_B(0)\int_0^{a_n} dr\int_0^{a_n}
  \phi_{t_1}(A_{a,r}, B_{a,s})ds,\nonumber
 \end{align}
where  $a_n=e^{-an}$. Applying  assumption (\ref{db1}), we further
get
 \begin{align}
  \phi_{t_1+ 2a}(A,B)
 \geq &c_0P_aI_A(0)P_aI_B(0)\int_0^{a_n } dr\int_0^{a_n}
  \phi_{t_2'}(A_{a,r}, B_{a,s})ds .\label{db3}
 \end{align}
Noticing that  $ \phi_{t_2'}(A_{a,r_1}, B_{a,s_1})\geq
\phi_{t_2'}(A_{a,r_2}, B_{a,s_2})$ if $0\leq r_1\leq r_2$ and $0\leq
s_1\leq s_2$,
 we have    \begin{align}
&P_aI_A(0)P_aI_B(0)\int_0^{a_n} dr \int_0^{a_n}
  \phi_{t_2'}(A_{a,r}, B_{a,s})ds\nonumber\\
\geq & a_n^2P_aI_A(0)P_aI_B(0)\int_0^{1} dr\int_0^{1}
  \phi_{t_2'}(A_{a,r}, B_{a,s})ds\nonumber\\
  =&a_n^2\int_0^{ P_aI_A(0)} dr\int_0^{ P_aI_B(0)}
\phi_{t_2'}(A_{a,r/P_aI_A(0)}, B_{a,s/P_aI_B(0)}   )ds
 \nonumber\\
 =& e^{ -2an} \phi_{t_2}(A,B). \label{db4}
 \end{align}
 Combing  (\ref{in})-(\ref{db4}), we get
 (\ref{db2}) when  $t_2>t_1+2a$.\qed\medskip

 \begin{lemma} \label{comb2}
 Let  $0\leq t\leq s$ and $\delta_0, c_0\in (0,1]$.
  Suppose that there exists  some integer
   $n_0$ such that   for every  $A,B\in \mathcal{C}_n$ with $B_n(\delta_0\sqrt{n})\subseteq
 A \cap B $ and every  $n\geq n_0$ \begin{align}
   \phi_{t}(A,B)\geq  c_0
\phi_{(s-2\delta_0)\wedge t}(A,B). \label{comb1}
\end{align}
   Then  for every   $ A, B \in \mathcal{C}_n$  and every   $n\geq n_0$
\begin{align}
 \phi_{t}(A,B)
\geq  c_0e^{-3\delta_0n}\phi_{s}(A,B).\nonumber
\end{align}

 \end{lemma}

 \noindent{\bf Proof}\ Let
 $ A, B \in \mathcal{C}_n$ and $n\geq n_0$. Setting  $ r_0=t_0=\delta_0$, we get
 by definition (\ref{AR}) and Lemma \ref{db0}
  \begin{align}
 B_{n}(\delta_0 \sqrt{n})
 \subseteq A_{\delta_0, r} \cap  B_{\delta_0, r},\ \
 \ \ \forall\  r\in [0,e^{-\delta_0n}],
 \nonumber
\end{align}
where $ A_{\delta_0, r}$ and $ B_{\delta_0, r}$ are  defined by
(\ref{AR}). By assumption (\ref{comb1}) and the estimate above
 \begin{align}
   \phi_t( A_{\delta_0,r}, B_{\delta_0,r'})\geq  c_0
   \phi_{(s-2\delta_0)\wedge t}(A_{\delta_0,r},
   B_{\delta_0,r'}),\ \ \ \ \ \forall\  r,r'\in [0,e^{-\delta_0n} ].\label{db3*}
\end{align}
 By taking $a=\delta_0$ in (\ref{db1}), we get
 the conclusion  by (\ref{db3*}) and  Lemma \ref{ke}.\qed\medskip

 Applying (\ref{3''}) and  Lemma \ref{comb2}, we get the following
 result.
\begin{corollary} \label{comb}
 Let $\delta_0, c_0\in (0,1]$.
  Suppose that there exists  some integer
   $n_0$ such that   for every  $A,B\in \mathcal{C}_n$ with $B_n(\delta_0\sqrt{n})\subseteq
 A\cap B $ and every   $n\geq n_0$
 \begin{align}
   \mu_n(A\cap B)\geq  c_0
 \mu_n(A) \mu_n(B).\nonumber
\end{align}
Then  for every   $ A, B \in \mathcal{C}_n$ and every  $n\geq n_0$
\begin{align}
   \mu_n(A\cap B)\geq  c_0\exp\{  -3\delta_0 n\}
 \mu_n(A) \mu_n(B).\nonumber
\end{align}

 \end{corollary}
Next   we prepare some basic formulas and estimates. We have
 \begin{align} \int
|x_i|^2d\mu_n=&1,\ \ \ for\ i=1,\cdots,n;\ \ \  \ \  \ \ \ \ \int |x|^4d\mu_n=n^2+2n;\label{4n}\\
 \mu_1([s,\infty))\leq&(2\pi)^{-1/2}s^{-1} e^{-s^2/2},\ \ \ \ \ \forall\ s>0.\label{1n}
\end{align}
By the first equality of (\ref{4n}) and Chebyshev inequality
\begin{align}  \mu_n(B_n(\sqrt{n}) ) > 1/2.
\label{1/2}\end{align} Notice  that in fact   $\lim_{n\rightarrow
\infty} \mu_n(B_n(\sqrt{n}) ) =1/2$ (c.f. \cite{SSZ98}). We have for
every $r\geq 0$ (c.f. \cite{GR07})
\begin{align} m_{n-1}(S_{n-1}(r))=\frac{2\pi^{n/2}}{\Gamma(n/2)}r^{n-1}.\label{sp}
\end{align}
By    Stirling formula (c.f. \cite{GR07}),
\begin{align} \Gamma(n/2)\sim \sqrt{\pi/n}e^{-n/2}n^{n/2}2^{-(n-2)/2}. \label{st}\end{align}

The constants  $N_l$ for $l\geq 1$ defined in what below will  be
used throughout the paper.
 \begin{lemma} \label{ball,es}
 There exists   some universal integer    $N_1$ such that
  for every  $r_0\in (0,1)$ and  every  $n\geq N_1$  \begin{align}
 \mu_n(B_n(r_0^{-1}\sqrt{n})^c)\leq &(1-r_0^2)^{-1}n^{-1/2}r_0^{-n+2}\exp\{-(r_0^{-2}-1)n/2\}. \label{>}
\end{align}

 \end{lemma}
 \noindent{\bf Proof}\
Let $0<r_0<1$.   For every $r>r_0^{-1}\sqrt{n }$,  by $r_0\in (0,1)$
we have $r^2\leq  (1-r_0^2)^{-1}(r^2-(n-2))$. By (\ref{sp}), we have
\begin{align}
  \mu_n(B_n(r_0^{-1}\sqrt{n})^c)=&
   \frac{1}{2^{(n-2)/2}\Gamma(n/2)} \int_{r_0^{-1}\sqrt{n }}^\infty
 r^{n-1}\exp\{-r^2/2\}dr\nonumber\\
\leq&
   -\frac{1}{(1-r_0^2)2^{(n-2)/2}\Gamma(n/2)} \int_{r_0^{-1}\sqrt{n }}^\infty
 d(r^{n-2}\exp\{-r^2/2\})\nonumber\\
 =&
   \frac{1}{(1-r_0^2)2^{(n-2)/2}\Gamma(n/2)}
   r_0^{-n+2}n^{-1+n/2}\exp\{-r_0^{-2}n/2\}.\nonumber
\end{align}
  Then the conclusion follows by      (\ref{st}) and the estimate above.
\qed\medskip

 \begin{lemma} \label{ball,es,2}
 There exists   some universal integer    $N_2$ such that for every
$C\geq2$
 and every   $n\geq N_2$
 \begin{align} \int_{|x|>C\sqrt{n}} |x|^3u(x)d\mu_n \leq
\exp\{-C^2n/6\}. \nonumber
\end{align}
 \end{lemma}
 \noindent{\bf Proof}\ Let $C\geq 2$.  For $n$ big enough, we have
 by (\ref{sp}) and (\ref{st})
\begin{align}
\int_{|x|>C\sqrt{n}} |x|^3d\mu_n
 = &\frac{1}{2^{(n-2)/2}\Gamma(n/2)}  \int_{C\sqrt{n}}^\infty  r^{n+2}
  \exp\{-  r^2/2\}dr  \nonumber\\
 \leq &\frac{1}{n 2^{(n-2)/2}\Gamma(n/2)} \int_{C\sqrt{n}}^\infty -d( r^{n+3}
\exp\{-  r^2/2\})dr \nonumber\\
 = &\frac{C^{n+3}n^{(n+1)/2}\exp\{-C^2n/2\}}{ 2^{(n-2)/2}\Gamma(n/2)} \nonumber\\
   \leq  &n C^{n+3}\exp\{-(C^2-1)n/2\}    .\nonumber
 \end{align}
Noticing that   $C^2-1-\frac{2\ln n}{n}-\frac{2(n+3)}{n}\ln C\geq
C^2/3$ for every  $C\geq 2$ when   $n$ is big enough,  we get the
conclusion from the estimate above. \qed\medskip
 \begin{lemma} \label{ball,es}
 There exists   some universal integer    $N_3$ such that  for every
  $u\in \mathcal{CF}_n$ and every   $n\geq
N_3$
\begin{align}
\int_{|x|>2\sqrt{n}} |x|^3u(x)d\mu_n \leq
e^{-n/2}\int_{|x|<\sqrt{n}} u(x)d\mu_n . \label{>>}
\end{align}

 \end{lemma}
 \noindent{\bf Proof}\ By assumption  $u\in \mathcal{CF}_n$,  we have  $u( r_1 x)\geq u( r_2x)$
for  every $x\in \mathbb{R}^n$ if   $r_2>r_1\geq 0$.  Then we have
for every $n\geq 2$
\begin{align} \int_{|x|>2\sqrt{n}} |x|^3u(x)d\mu_n
 = &\frac{1}{(2\pi)^{n/2} } \int_{S_{n-1}} dm_{n-1}(\widetilde{x})\int_{2\sqrt{n}}^\infty  r^{n+2}
 u(r\widetilde{x})\exp\{-  r^2/2\}dr  \nonumber\\
 \leq &\frac{1}{ n(2\pi)^{n/2}}\int_{S_{n-1}}  u( \sqrt{n}\widetilde{x})dm_{n-1}
 (\widetilde{x})\int_{2\sqrt{n}}^\infty -d( r^{n+3}
\exp\{-  r^2/2\})dr \nonumber\\
 = &\frac{2^{n+3}n^{(n+1)/2}e^{-2n}}{ (2\pi)^{n/2}}\int_{S_{n-1}}
 u( \sqrt{n}\widetilde{x})dm_{n-1}(\widetilde{x}) .\nonumber
 \end{align}
We also have
\begin{align}
\int_{|x|<\sqrt{n}} u(x)d\mu_n
 = &\frac{1}{(2\pi)^{n/2}} \int_{S_{n-1}} dm_{n-1}(\widetilde{x})\int_0^{\sqrt{n}} r^{n-1}
 u(r\widetilde{x})\exp\{-  r^2/2\}dr  \nonumber\\
 \geq &\frac{1}{ n (2\pi)^{n/2}}\int_{S_{n-1}}  u( \sqrt{n}\widetilde{x})dm_{n-1}
 (\widetilde{x})\int_0^{\sqrt{n}}  d( r^{n}
\exp\{-  r^2/2\})dr \nonumber\\
 \geq &\frac{n^{n/2}e^{-n/2}}{n (2\pi)^{n/2}}\int_{S_{n-1}}
 u( \sqrt{n}\widetilde{x})dm_{n-1}(\widetilde{x}) .\nonumber
 \end{align}
Combing the two estimates above, we get (\ref{>>}) for $n$ big
enough.\qed

\subsection{\normalsize  Hessen estimates
  }

The main subject of this subsection is to give a quantity version of
   Lemma  \ref{Pre}. First we prepare some formulas for the later use. The Ornstein-Uhlenbeck semigroup $(P_t)_{t\geq0}$  can be written
also as
\begin{align}
P_tu(x)=&\int u(e^{- {t}/2}x+(1-e^{-t})^{1/2}y)d\mu_n(y), \ \ \ \ \
\forall  x\in \mathbb{R}^n, \forall t\geq0,\label{OU}
\end{align}
where $u$ is a  bounded measurable function for instance.
    For    random variables $X$ and
$Y$ on some probability space $(\Omega,\mathcal{F},P)$ with finite
second moments, we have
\begin{align}\int_{\Omega}X (Y-\int_{\Omega}YdP)dP =
&\int_{\Omega}(X- \int_{\Omega}XdP)(Y-\int_{\Omega}YdP)dP,\label{momen'}\\
\int_{\Omega}(X- \int_{\Omega}XdP)(Y-\int_{\Omega}YdP)dP
=&\int_{\Omega}XYdP- \int_{\Omega}X dP\int_{\Omega} YdP,
\label{momen} \end{align} and
\begin{align}\label{moment}
\int_{\Omega}(X-\int_{\Omega} XdP)^2dP\leq \int_{\Omega}(X-b)^2dP,\
\ \ \ \  \  \ \ \forall\  b\in \mathbb{R}.\end{align}

 For  functions $f$ and $g$ on $\mathbb{R}^n$,   denote by $f\ast g$ the
convolution of $f$ and $g$ when it exists.
 Let  $u,v $ be positive measurable function on $\mathbb{R}^n$  and suppose that $u=e^{-U}$ is  smooth.
Define function $\widetilde{U}$   by $u\ast v=e^{-\widetilde{U}}$.
Next we derive a formula for  the  partial derivatives of
$\widetilde{U}$. We assume that   all the integrals involved below
are well defined. Let  $x\in \mathbb{R}^n$.  Define probability
measure  $\gamma_x$ on $ \mathbb{R}^n$  by
\begin{align}
 d\gamma_x(y) =\frac{1}{(u\ast v)(x)} u(x-y)v(y)dy.\label{ab}
\end{align}  Let $m\geq 1$. Set
$\mathcal{A}_{n,m}=\{1,\cdots,n\}^m$.   For every  $\mathbf{a}\in
\mathcal{A}_{n,m}$,   write $\mathbf{a}=(a_1,\cdots,a_m)$ and $
\partial_{\mathbf{a}}=
\partial_{a_1}\cdots\partial_{a_m}$. Define for every  $y\in \mathbb{R}^n$, $D\subseteq \{1,\cdots,m\}$ and every $\mathbf{a}\in
\mathcal{A}_{n,m}$
\begin{align}
  & \xi_{n,\mathbf{a}, D}(x,y)  = \big(\prod_{i\in D} \partial_{a_i}\big)U(x-y)-
  \int \big(\prod_{i\in D} \partial_{a_i}\big) U(x-y)
d\gamma_x(y)  . \nonumber
\end{align}
 For every $\mathbf{a}\in \mathcal{A}_{n,m}$, $k\geq 1$
and every $\Theta =(\Theta_1,\cdots,\Theta_k)$ such that
$\Theta_l\subseteq \{1,\cdots,m\}$ for $1\leq l\leq k$,  denote
\begin{align}
 \xi_{n,\mathbf{a} }^{ \Theta}(x)  =  \int \prod_{l=1}^k\xi_{n,\mathbf{a}, \Theta_l}(x,y)
 d\gamma_x(y).
 \label{de}
\end{align}   Denote by $\mathcal{P}_m$
the set of all      partitions of $\{1,\cdots,m\}$. For $1\leq k\leq
m$, denote by $\mathcal{P}_{m,k}$ the subset of $\mathcal{P}_m$ such
that a partition belongs to $\mathcal{P}_{m,k}$ if and only if it
  contains  exactly $k$ sets. For $\Theta \in \mathcal{P}_{m,k}$, denote
$\Theta=(\Theta_1,\cdots,\Theta_k)$, where
$(\Theta_1,\cdots,\Theta_k)$ is the partition corresponding to
$\Theta$.
 Denote for every $\mathbf{a}\in \mathcal{A}_{n,m}$
\begin{align}
  K_{n,\mathbf{a}}^{(1)}(x)  = &\sum_{k=2}^m
  \sum_{\Theta \in \mathcal{P}_{m,k}} (-1)^{k+1}
  \xi_{n,\mathbf{a}}^{ \Theta}(x),
  \label{kn}
\end{align}
where  the right hand side above is assumed to be  zero when $m=1$.
Noticing that $\xi_{n,\mathbf{a}}^{ \Theta}(x)=0$ when $\Theta \in
\mathcal{P}_{m,1}$, the summation in the  right hand side of
(\ref{kn}) can start from $k=1$.
  For  $D\subseteq \{1,\cdots,m\}$  and $\mathbf{a}\in
  \mathcal{A}_{n,m}$, set
\begin{align}\mathbf{a}(D)=(a_{l_i})_{i=1}^{|D|}\in \mathcal{A}_{n,|D|},\label{a(.)}\end{align} where
  $|D|$ is the cardinal number of $D$ and $(l_i)$
  is the unique increasing map from $\{1,\cdots,|D|\}$ to $D$.
  Set for every $\mathbf{a}\in \mathcal{A}_{n,m}$ and  every  $k\in \{1,\cdots,m\} $
 \begin{align}
  K_{n,\mathbf{a}}^{(k)}(x)  = &   \sum_{\Theta \in \mathcal{P}_{m,k}}\prod_{l=1}^k
  K_{n,\mathbf{a}(\Theta_l)}^{(1)}(x)   .   \label{ln}
  \end{align}
  Notice that the  definition above is consistent with (\ref{kn}) when $k=1$.
In what below we make convention that $k!=1$ when $k=-1,0$.

\begin{lemma}\label{J1,J2J3}
 Let  $u,v $ be positive measurable function on $\mathbb{R}^n$  and suppose that $u=e^{-U}$ is  smooth.
 Assume that  $u\ast v$ is a  well defined function  and  set   function
$\widetilde{U}$   by $u\ast v=e^{-\widetilde{U}}$. Let $\gamma_x$ be
the probability measure specified by  (\ref{ab}).
  Let be $m_0$ be a positive integer,
$\mathbf{a}\in \mathcal{A}_{n,m}$ and assume that all the integrands
   below  with respect to $dx$ are controlled by some     integrable function  for all
$x\in \mathbb{R}^n$. Then for every $m=1,\cdots,m_0$ and every $x\in
\mathbb{R}^n$
\begin{align}
    \partial_{\mathbf{a}}  \widetilde{U} (x)= \int
\big(\prod_{i=1}^m \partial_{a_i}\big) U(x-y) d\gamma_x(y) +
\sum_{k=1}^m
(k-1)!K_{n,\mathbf{a}}^{(k)}(x),\label{J1J2J3}\end{align}
  where   $K_{n,\mathbf{a}}^{(k)}$ is defined by (\ref{ln}).

\end{lemma}
  \noindent{\bf Proof}\ \ Direct calculation shows that  (\ref{J1J2J3}) holds  for
  $m=1$. Applying  the  method of finite induction, next we assume that (\ref{J1J2J3}) holds  for
 some  $1\leq m \leq m_0-1$ and prove it holds also for $m+1$.
  Let $\overline{\mathbf{a}}\in \mathcal{A}_{n,m+1}$ and
  write $\overline{\mathbf{a}}=(  \mathbf{a},a_{m+1})$ with
  $\mathbf{a}\in \mathcal{A}_{n,m}$.
For every $k\geq 1$,   $\Theta=(\Theta_1,\cdots,\Theta_m )\in
\mathcal{P}_{m,k}$ and every $i \in \{1,\cdots,m\} $, denote by
$\Theta{(i;m+1)} $ the partition  in  $ \mathcal{P}_{m,k}$ which is
equal to $\Theta$ with $\Theta_i $ replaced by
$\Theta_i\cup\{m+1\}$; denote by $\Theta(\{ m+1\}) $ the partition
in $ \mathcal{P}_{m,k+1}$ which is equal to $\Theta$ with
$\Theta_{m+1}:=\{m+1\}$ added ;  denote by $\Theta(i;-) $ the
partition in  $ \mathcal{P}_{m,k-1}$ which is equal to $\Theta$ with
$\Theta_i$ removed. Here we assume that $ \mathcal{P}_{m,0}$ is an
empty set.
  For every  $ x\in \mathbb{R}^n$, we have by induction assumption
  \begin{align}
   \partial_{\overline{\mathbf{a}}}  \widetilde{U} (x) =&\sum_{k=0}^m    (k-1)!\partial_{a_{m+1}}
K_{n,\mathbf{a}}^{(k)}(x),
    \ \ \ \ \label{123}\end{align}
    where
      \begin{align}
K_{n,\mathbf{a}}^{(0)}(x) :=
   \int
\big(\prod_{i=1}^m \partial_{a_i}\big) U(x-y)
d\gamma_x(y).\nonumber\end{align} When  $1\leq i\leq n$, by
(\ref{momen'}) we have for every  $x\in \mathbb{R}^n$
\begin{align}\partial_{x_i}\gamma_x(dy) =
-\big(\partial_iU(x-y) -\int
\partial_iU(x-y)d\gamma_x(y)\big)d\gamma_x(y),\ \ \ \  \forall y\in \mathbb{R}^n.\label{f1}
  \end{align}
For every $k\geq 1$, $\mathbf{b}\in \mathcal{A}_{n,k}$,
$D\subseteq\{1,\cdots,k\}$ and every $x,y\in \mathbb{R}^n$, we have
by  (\ref{f1})
    \begin{align}&\partial_{x_i}  \xi_{n,\mathbf{b},D }(x,y) \nonumber\\=&
     \xi_{n,(\mathbf{b},i),D\cup\{k+1\} }(x,y)
+
\int   \xi_{n,\mathbf{b} ,D}(x,y) \big(  \partial_iU(x-y)- \int
 \partial_iU(x-y)\gamma_x(y) \big)d\gamma_x(y)  .\label{f2}
  \end{align}
By  (\ref{momen'}) , definition (\ref{de}) and  (\ref{f1}), we have
for every $x\in \mathbb{R}^n$
\begin{align}
  \partial_{a_{m+1}}
  K_{n,\mathbf{a}}^{(0)}(x)=&K_{n,\overline{\mathbf{a}}}^{(0)}(x)-
 \xi_{n,
\overline{\mathbf{a}}}^{\{S_m,\{m+1\}\}}(x)
:=K_{n,\overline{\mathbf{a}}}^{(0)}(x)-J_{0}(\overline{\mathbf{a}},x)
,\nonumber\end{align} where $S_m=\{1,\cdots,m \}$ and
$\{S_m,\{m+1\}\}$ is the partition in $\mathcal{P}_{m+1,2}$
containing
  two sets $S_m$ and $\{m+1\}$.
  The equality above and  (\ref{123}) show that
\begin{align}
   \partial_{\overline{\mathbf{a}}}  \widetilde{U} (x)-K_{n,\overline{\mathbf{a}}}^{(0)}(x)
    =&\sum_{k=1}^m   (k-1)! \partial_{a_{m+1}}
K_{n,\mathbf{a}}^{(k)}(x)-J_0(\overline{\mathbf{a}},x)  .
    \ \ \ \ \label{1,2,3}\end{align}

Let $x\in \mathbb{R}^n$. We have by (\ref{f1}) and (\ref{f2})
    \begin{align}
\partial_{a_{m+1}}
K_{n,\mathbf{a}}^{(1)}(x)
    =&
    \partial_{a_{m+1}}\sum_{k=2}^m (-1)^{k+1}\sum_{\Theta \in \mathcal{P}_{m,k}}
 \int     \prod_{l=1}^k\xi_{n,\mathbf{a},\Theta_l}(x,y)
  \gamma_x(dy)\nonumber\\        = &
   \sum_{k=2}^m (-1)^{k+1}\sum_{\Theta \in \mathcal{P}_{m,k}}
     \int
     \prod_{l=1}^k \xi_{n,\mathbf{a},\Theta_l}(x,y)
    \partial_{x_{a_{m+1}}} \gamma_x(dy) \nonumber\\+&
   \sum_{k=2}^m (-1)^{k+1}\sum_{\Theta \in \mathcal{P}_{m,k}}
    \sum_{i=1}^k  \int   \big( \partial_{x_{a_{m+1}}}
     \xi_{n,\mathbf{a},\Theta_i}(x,y)\big)
     \big( \prod_{l=1,l\neq
     i}^k\ \xi_{n,\mathbf{a},\Theta_l}(x,y)\big)d
     \gamma_x(y)\nonumber\\
=& \sum_{k=2}^m(-1)^{k+2}\sum_{\Theta \in \mathcal{P}_{m,k}}
     \xi_{n,\overline{\mathbf{a}}}^{\Theta(\{m+1\})}(x)+ \sum_{k=2}^m (-1)^{k+1}\sum_{\Theta \in \mathcal{P}_{m,k}}
     \sum_{i=1}^k
     \xi_{n,\overline{\mathbf{a}}}^{\Theta(i;m+1)}(x)
    \nonumber\\
     +&\sum_{k=2}^m (-1)^{k+1}\sum_{\Theta \in \mathcal{P}_{m,k}}
    \sum_{i=1}^k  \xi_{n,\overline{\mathbf{a}}}^{ \Theta(i;-)}(x)
    \xi_{n,\overline{\mathbf{a}}}^{ \{\Theta_i,\{m+1\}\}}(x)\nonumber\\
 :=&J_{1,1}(\overline{\mathbf{a}},x)+J_{1,2}(\overline{\mathbf{a}},x)+J_{1}(\overline{\mathbf{a}},x) .\label{222}
\end{align}
Noticing that  for every $k\geq 2$\begin{align}
\mathcal{P}_{m+1,k}=&\{\Theta(i;m+1):1\leq i\leq k;\ \Theta\in
\mathcal{P}_{m,k}\}\cup\{\Theta(\{m+1\}): \ \Theta\in
\mathcal{P}_{m,k-1}\}\nonumber,
\end{align}
we have by
  (\ref{222})
  \begin{align}
 \partial_{a_{m+1}}K_{n,
{\mathbf{a}}}^{(1)}(x)=K_{n,  \overline{\mathbf{a}}}^{(1)}(x)+
J_0(\overline{\mathbf{a}},x)+J_1(\overline{\mathbf{a}},x).\nonumber
\end{align}
This and (\ref{1,2,3}) imply
\begin{align}
   \partial_{\overline{\mathbf{a}}}  \widetilde{U} (x)-\sum_{k=0,1}
   K_{n,\overline{\mathbf{a}}}^{(k)}(x) =&\sum_{k=2}^m  (k-1)!  \partial_{a_{m+1}}
K_{n,\mathbf{a}}^{(k)}(x)+J_1(\overline{\mathbf{a}},x).\nonumber\end{align}
By  the method of finite induction,  we can verify for $1\leq j\leq
m-1$
\begin{align}
   \partial_{\overline{\mathbf{a}}}  \widetilde{U} (x)-\sum_{k=0}^j
   (k-1)!
   K_{n,\overline{\mathbf{a}}}^{(k)}(x) =&\sum_{k=j+1}^m    (k-1)!\partial_{a_{m+1}}
K_{n,\mathbf{a}}^{(k)}(x)+
j!J_j(\overline{\mathbf{a}},x),\label{ba}\end{align} where
\begin{align}
J_j(\overline{\mathbf{a}},x)=\sum_{\Theta \in \mathcal{P}_{m,j}}
\sum_{l'=1}^j \prod_{l=1,l\neq l'}^j
  K_{n,\mathbf{a}(\Theta_l)}^{(1)}(x) J_1
  (\overline{\mathbf{a}}(\Theta_{l'}\cup\{m+1\}),x).\nonumber
\end{align}
where $\overline{\mathbf{a}}(\cdot)$ is defined by (\ref{a(.)}). The
appearance of $j!$ is due to that there are $j$ ways  to get a
partition in $\mathcal{A}_{n,j}$ from  a partition $\Theta$ in
$\mathcal{A}_{n,j+1}$ by   combining  a given  element of $\Theta$
with one of the others. Noticing that
$K_{n,\overline{\mathbf{a}}}^{(m+1)}(x)=
K_{n,\overline{\mathbf{a}}}^{(m)}(x)=K_{n, {\mathbf{a}}}^{(m)}(x)
=0$, we get the
 conclusion by (\ref{ba}).
 \qed\medskip

Applying Lemma \ref{J1,J2J3}, we get the following result which can
also be  checked directly.
\begin{corollary} \label{derivatives,a11}
 Let  $u,v $ be positive measurable function on $\mathbb{R}^n$  and suppose that $u=e^{-U}$ is  smooth.
Assume that  $u\ast v$ is a well defined function and set function
$\widetilde{U}$  by $u\ast v=e^{-\widetilde{U}}$. Let $1\leq i,j,k
\leq n$ and let $\gamma_x$ be the probability measure given by
(\ref{ab}). Then we have for every $x\in \mathbb{R}^n$
\begin{align}
    \partial_i\partial_j  \widetilde{U} (x)=J_1(x)-J_2(x),\label{J1J2}\end{align}
  where  \begin{align} J_1(x) =&
  \int
\partial_i\partial_j
U(x-y) d\gamma_x(y),\nonumber\\
J_2(x)=&   \int \big(\partial_i  U(x-y)-\int
\partial_i  U(x-y)d\gamma_x(y)\big)\big(
\partial_jU(x-y)
  -\int  \partial_j U(x-y)d\gamma_x(y)\big)d\gamma_x(y)
  \nonumber,
\end{align}
Moreover, for every $x\in \mathbb{R}^n$\begin{align}
    \partial_i\partial_j\partial_k  \widetilde{U} (x)=M_1(x)-M_{2,1}(x)-M_{2,2}(x)-M_{2,3}(x)+M_3(x),\end{align}
  where  \begin{align}
M_1(x) =&
  \int
\partial_i\partial_j\partial_k
U(x-y) d\gamma_x(y),\nonumber\\
M_{2,1}(x)=&   \int \big(\partial_i\partial_j U(x-y)-\int
\partial_i\partial_j U(x-y)d\gamma_x(y)\big)\big(
\partial_kU(x-y)
  -\int  \partial_k U(x-y)d\gamma_x(y)\big)d\gamma_x(y),
  \nonumber\\M_{2,2}(x)=&
 \int \big(\partial_i
\partial_kU(x-y)- \int\partial_i
\partial_kU(x-y)d\gamma_x(y)\big)\big(
\partial_jU(x-y) - \int\partial_j
U(x-y)d\gamma_x(y)\big) d\gamma_x(y),\nonumber\\ M_{2,3}(x)=& \int
\big(\partial_j
\partial_k
U(x-y)- \int\partial_j
\partial_k
U(x-y)d\gamma_x(y)\big)\big(
\partial_i
U(x-y) - \int\partial_i U(x-y)d\gamma_x(y)\big)
d\gamma_x(y),\nonumber\\
M_3(x)=&   \int\Big[ \big(\partial_iU(x-y)- \int\partial_i U(x-y)
d\gamma_x(y)\big) \big(\partial_j U(x-y)- \int\partial_j U(x-y)
d\gamma_x(y)\big)
\nonumber\\
& \ \ \ \ \ \ \ \ \ \ \ \  \ \ \ \ \cdot\big(\partial_k U(x-y)-
\int\partial_k U(x-y)d\gamma_x(y)\big) \Big] d\gamma_x(y). \nonumber
\end{align}
Here  we assume  that   all the integrands above  with respect to
$dx$ are controlled by some integrable function  for all  $x\in
\mathbb{R}^n$.
\end{corollary}

 For  every  $u=e^{-U}
\in\mathcal{CF}_n$ and every $t\geq 0$,   define   functions $u_t$
and $U_{ t }$ by
\begin{align} u_t=P_tu;\ \ \ \ \ \ U_t=-\ln u_t.\label{HA}
\end{align}
 By Lemma \ref{Pre},   $U_t$ is a  convex function.

 \begin{lemma} \label{derivatives,2}
 Let  $u=e^{-U}\in \mathcal{CF}_n$.   Then for every  $ x\in \mathbb{R}^n$ every  $t>0$
\begin{align}
      0\leq  \nabla^2  U_t(x) \leq &
   2(1\wedge t)^{-1}e^{-t} I_n.\label{de,3}\end{align}
  If   further assuming  that $U$ is twice differentiable and $\nabla^2 U \leq CI_n$ on $\mathbb{R}^n$ for some
constant   $C>0$,
  then  for  every $t\geq 0$ and  every  $ x\in \mathbb{R}^n$ \begin{align}
      0\leq  \nabla^2  U_t(x) \leq &
   e^{-t}C I_n. \label{de,131}\end{align}  \end{lemma}
  \noindent{\bf Proof}\ \ The first inequality  of (\ref{de,3}) and
  the  first inequality  of
  (\ref{de,131})
 follow by Lemma \ref{Pre}. Let $t> 0, 1\leq i\leq n$ and $x\in \mathbb{R}^n $.
 Let    $\sigma_x$ be the  probability
measure on $\mathbb{R}^n$  defined  by
\begin{align}
d\sigma_x(y)= \frac{u(y)}{(2\pi (1-e^{-t}))^{n/2} P_tu(x)}
\exp\{-\frac{|y-e^{-t/2}x|^2}{2(1-e^{-t})}\}dy.\nonumber
\end{align}
By the definition of $U_t$,  (\ref{P1}) and   applying
(\ref{J1J2}), we have
\begin{align}\partial_i^2U_t(x)
= & \frac{ e^{-t}}{ 1-e^{-t} }  - \frac{ e^{-t}}{ (1-e^{-t})^2 }
\int \Big(y_i-e^{-t/2}x_i-\int (y_i-e^{-t/2}x_i) d\sigma_x(y)\Big)^2
 d\sigma_x(y ) \nonumber\\
\leq & 2 (1\wedge t)^{-1}e^{-t} . \nonumber
\end{align}
 Since the estimate above holds under any  coordinate system
 $(Q(\mathbf{e}_i))_{1\leq i \leq n}$ when $Q$ is an  orthogonal
 transformation
   of
$\mathbb{R}^n$, the estimate above implies the  second inequality of
(\ref{de,3}).

Let $\widetilde{\sigma}_x$ be the  probability measure on
$\mathbb{R}^n$  defined  by
\begin{align}
d\widetilde{\sigma}_x(y)=
\frac{u(e^{-t/2}x+(1-e^{-t})^{1/2}y)}{(2\pi)^{n/2} P_tu(x)}
\exp\{-\frac{|y|^2}{2}\}dy.\nonumber
\end{align} Applying
(\ref{OU}) and (\ref{J1J2}), we have
\begin{align}\partial_i^2U_t(x)
= &  e^{-t}\int
\partial_i^2U(e^{-t/2}x+(1-e^{-t})^{1/2}y)d\widetilde{\sigma}_x(y )\nonumber\\
-&e^{-t} \int    \Big(\partial_iU(e^{-t/2}x+(1-e^{-t})^{1/2}y)-\int
\partial_i^2U(e^{-t/2}x+(1-e^{-t})^{1/2}y)d\widetilde{\sigma}_x(y)\Big)^2
 d\widetilde{\sigma}_x(y ) \nonumber\\
\leq & Ce^{-t}, \nonumber
\end{align}
which gives the second inequality of (\ref{de,131}).
 \qed\medskip

Let    $u=e^{-U} \in\mathcal{CF}_n$,   $t,s\geq 0$ and $x\in
\mathbb{R}^n $. Define  Borel measure $\nu_{u,t,s,x}$ on
$\mathbb{R}^n$ by
\begin{align}
 d\nu_{u,t,s,x}(y)=&u_{t+s}(x)^{-1}
 u_t(e^{-s/2}x+(1-e^{-s})^{1/2}y)d\mu_n(y)
 , \label{vts}\end{align}
where $u_t$ is defined by (\ref{HA}).  In what below
$\nu_{u,t,s,x}$ is also written in short as $\nu_x$   when   it
makes  no confusion. By definition (\ref{HA}) and the semigroup
property of $(P_t)$,\begin{align}
 \nu_x(\mathbb{R}^n)=&u_{t+s}(x)^{-1}\int
 u_t(e^{-s/2}x+(1-e^{-s})^{1/2}y)d\mu_n(y)\nonumber\\
 =&u_{t+s}(x)^{-1}P_s
u_t(x)\nonumber\\
 =&1,\nonumber\end{align}
 which shows that  $\nu_x$ is a probability measure. Define
   function $U_{t,s,x}$
 on $\mathbb{R}^n$ by
\begin{align}
 d\nu_x(y)=&
 \exp\{-U_{t,s,x}(y)\}dy.\nonumber\end{align}
  From (\ref{HA}) and (\ref{vts}), we have
\begin{align}
    U_{t,s,x}(y)=-
   U_{t+s}(x)+
  U_t(e^{-s/2}x+(1-e^{-s})^{1/2}y)+\frac{|y|^2}{2},\ \ \
 \forall\  y\in \mathbb{R}^n.\label{h1}\end{align}
For  every  $t,s\geq 0$ and every  $x\in \mathbb{R}^n$,   by
(\ref{h1}) and the convexity of $U_t$,
  there exists an  unique element  $x^*=x^*(t,s,x)\in \mathbb{R}^n$
 such that
 \begin{align}
 U_{t,s,x}(x^*)=\inf_{y\in \mathbb{R}^n} U_{t,s,x}(y).\label{hess}
 \end{align}

 \begin{lemma} \label{derivatives,1}
 Let  $u=e^{-U} \in\mathcal{CF}_n$ be a smooth function
 such that  $\nabla^2  U  \leq C I_n$ on   $  \mathbb{R}^n$ for some constant   $C>0$.  Let  $t,s\geq 0$, $x\in
\mathbb{R}^n $
  and define     $x^*\in \mathbb{R}^n$
      by (\ref{hess}).
  Then   for every   $n\geq N_2$
   \begin{align}
 \int |y-x^*|^k   d\nu_x(y) \leq &((3+C)n)^{k/2} ,\ \ \ for\ k=2,3\label{<>4}.
  \end{align}
When  further assuming that $Cs n \leq 1$, we have
   \begin{align}
 \int |y_i-x^*_i|^2 d\nu_x(y) \leq & 2 ,\ \ \ \ \ for\ i=1,\cdots,n, \label{<>21}\\
 and\ \ \ \ \ \ \int |y-x^*|^4 d\nu_x(y) \leq & 6n^2 .\label{<>20}
  \end{align}

 \end{lemma}\noindent{\bf Proof}\ \ We have by (\ref{h1}), (\ref{hess})
 and the
 assumption  $\nabla^2  U \leq C I_n$ on   $\mathbb{R}^n$
\begin{align}
 U_{t,s,x}(x^*)+ \frac{1}{2}
  |y-x^*|^2 \leq  U_{t,s,x}(y)
 \leq   U_{t,s,x}(x^*)
 +\frac{1+(1-e^{-s})C}{2}  |y-x^*|^2,\ \ \ \forall\  y\in \mathbb{R}^n . \label{<>6} \end{align}
 Applying  the second inequality of (\ref{<>6}) and $
 \int \exp\{-U_{t,s,x}(y)\}dy=
 \nu_x(\mathbb{R}^n) =1$, we have
   \begin{align}
 e^{-U_{t,s,x}(x^*)}\leq & \Big(\int \exp\{-
 \frac{1+(1-e^{-s})C}{2}|y-x^*|^2\}dy\Big)^{-1}\nonumber\\
 =&(2\pi)^{-n/2}
 (1+(1-e^{-s})C )^{n/2}
.\label{ea,2}
 \end{align}
By the first inequality of (\ref{<>6}) and  (\ref{ea,2}),
\begin{align}
  &\int_{|y-x^*|> \sqrt{(2+C)n}} |y-x^*|^2 d\nu_x(y) \nonumber\\\leq&
   e^{-U_{t,s,x}(x^*)}  \int_{|y-x^*|> \sqrt{(2+C)n}} |y-x^*|^2\exp\{-\frac{|y-x^*|^2}{2}\}
   dy\nonumber\\\leq &
 (1+C )^{n/2}
    \int_{|y|>\sqrt{(2+C)n}} |y|^2 d\mu_n(y)  .\nonumber
 \end{align}
  By Lemma \ref{ball,es,2} and the
estimate above, we have for every $n\geq N_2$
\begin{align}
  \int_{|y-x^*|> \sqrt{(2+C)n}} |y-x^*|^2 d\nu_x(y) \leq
 (1+C )^{n/2}e^{-(2+C)n/6} \leq 1.\nonumber
 \end{align}
Therefore, for every $n\geq N_2$    \begin{align} \int  |y-x^*|^2
d\nu_x(y) \leq \int_{|y-x^*|\leq  \sqrt{(2+C)n}} |y-x^*|^2 d\nu_x(y)
+1\leq (3+C)n ,\nonumber
 \end{align}
   which gives (\ref{<>4}) for $k=2$.  With the same  calculation as
   above,
  (\ref{<>4}) holds also for $k=3$.

  Next we assume that $ Cs n\leq 1$.
Let $1\leq i \leq n$. We have
\begin{align}
(1+r)^{n/2}\leq e^{nr/2}\leq 1+nr,\ \ if\  r\in
(0,\frac{1}{n}].\label{ccc}\end{align} By  (\ref{<>6})-(\ref{ccc})
 and the assumption
 $Cs n\leq 1$,
       \begin{align}
  \int |y_i-x_i^*|^2 d\nu_x(y) \leq&e^{- U_{t,s,x}(x^*)}
 \int  |y_i-x^*_i|^2\exp\{-\frac{1}{2}|y-x^*|^2\} dy\nonumber\\
 \leq &(1+Cs)^{n/2} \int |y_i|^2d\mu_n(y)\nonumber\\
 \leq & 1+ Cns\nonumber\\
 \leq & 2,\nonumber
 \end{align}
which gives    (\ref{<>21}).
 Similarly, applying  the second equality of (\ref{4n}) we have
  \begin{align} \int |y-x^*|^4 d\nu_x(y) \leq &
   ( 1+Cs )^{n/2} \int |y|^4d\mu_n(y) \leq
 2(n^2+2n),\nonumber
 \end{align}
which gives   (\ref{<>20}).
 \qed
 \medskip

Let  $u=e^{-U}\in \mathcal{CF}_n$. In what below, we say that $u$
satisfies condition $\mathcal{L}(C_1,C_2)$ for some constants
$0\leq C_1<C_2$ if $U$ is a smooth function  and
\begin{align}
      C_1 I_n\leq  \nabla^2  U(x) \leq  &C_2 I_n,\ \ \ \ \forall\  x\in \mathbb{R}^n. \label{de,3;1} \end{align}
For   $0\leq C_1<C_2$ and  $C_3>0$, we say
 that $u$
satisfies condition $\mathcal{L}(C_1,C_2,C_3)$ if $u$ satisfies
condition $\mathcal{L}(C_1,C_2)$  and
\begin{align}
      |\partial_{i}\partial_{j}\partial_{k} U(x)| \leq & C_3
      ,\ \ \ \ \forall\  x\in \mathbb{R}^n,\ \forall i,j,k\in \{1,\cdots,n\}. \label{de,3;2}\end{align}

 \begin{lemma} \label{derivatives,ab11}
Let  $u=e^{-U} \in\mathcal{CF}_n$ satisfying condition
$\mathcal{L}(C_1,C_2,C_3)$ for some constants    $0\leq C_1<C_2$ and
$C_3>0$.
  Let    $1\leq i,j,k \leq n,x\in \mathbb{R}^n$ and $t\geq 0$.
 Then for every $n\geq N_2$
\begin{align}
  | \partial_i\partial_j\partial_k U_t(x)|\leq&C_4e^{-3t/2}n^{3/2} \label{j3},\end{align}
  where $U_t$ is defined by  (\ref{HA}) and $C_4=C_3+ C_2^2(6+8C_2)(3+C_2)^{3/2} $.\end{lemma}

 \noindent{\bf Proof}\   Let    $1\leq i,j,k \leq n,x\in \mathbb{R}^n$ and $t\geq
 0$. Noticing that  (\ref{j3}) holds for  $t=0$, we assume   that
 $t>0$ in what below.   Let   $\lambda_x$  be the measure    defined by  $\nu_{u,0,t,x}$ in
(\ref{vts}). More explicitly,
\begin{align}
 d\lambda_x(y)=&u_{t}(x)^{-1}
 u(e^{-t/2}x+(1-e^{-t})^{1/2}y)d\mu_n(y).\nonumber\end{align}    Define $x^*\in\mathbb{R}^n$ by (\ref{hess})
 corresponding to $\nu_x=\nu_{u,0,t,x}$.
 Write for every  $y\in
\mathbb{R}^n$\begin{align} \overline{y}= e^{-
t/2}x+(1-e^{-t})^{1/2}y.\label{-y}
\end{align}
Applying Corollary  \ref{derivatives,a11}, we have
 \begin{align}
 &\partial_i\partial_j\partial_k
U_t(x) =e^{-3t/2}
(M_1(x)-M_{2,1}(x)-M_{2,2}(x)-M_{2,3}(x)+M_3(x)),\label{I123}
\end{align}
where  \begin{align}
 M_1(x)
=&  \int
\partial_i\partial_j\partial_kU
( \overline{y})d\lambda_x(y),\nonumber\\
M_{2,1}(x)=&   \int \big(\partial_i\partial_jU( \overline{y})-\int
\partial_i\partial_jU(
\overline{y})d\lambda_x(y)\big)\big(
\partial_kU( \overline{y})
  -\int  \partial_kU(
\overline{y})d\lambda_x(y)\big)d\lambda_x(y) ,\nonumber\\
M_{2,2}(x)=&
  \int \big(\partial_i
\partial_k
U( \overline{y})- \int\partial_i
\partial_k
U( \overline{y})d\lambda_x(y)\big)\big(
\partial_j
U( \overline{y})-  \int\partial_j U( \overline{y})d\lambda_x(y)\big)
d\lambda_x(y),\nonumber\\  M_{2,3}(x)=& \int \big(\partial_j
\partial_k
U( \overline{y})- \int\partial_j
\partial_k
U( \overline{y})d\lambda_x(y)\big)\big(
\partial_i
U( \overline{y}) - \int\partial_i U( \overline{y})d\lambda_x(y)\big)
d\lambda_x(y),\nonumber\\
 M_{3}(x)=&   \int   \prod_{l=i,j,k}
\big(\partial_lU( \overline{y})- \int\partial_lU( \overline{y})
d\lambda_x(y)\big)   d\lambda_x(y)\nonumber .
\end{align}

We always assume that $n\geq N_2$ in what below. By (\ref{de,3;2}),
\begin{align} |M_1(x)| =|\int
\partial_i\partial_j\partial_k U
( \overline{y})d\lambda_x(y)|\leq C_3.\label{I1}
\end{align}

  By the second inequality of  (\ref{de,3;1}), we
have $\sum_{j=1}^n(\partial_j\partial_m
 U)^2\leq C_2^2$ on $ \mathbb{R}^n$. Then, applying (\ref{<>4}), (\ref{-y}),
  mean value theorem and Cauchy-Schwartz
inequality, we get for every $m\in \{1,\cdots,n\}$
\begin{align} &
\big(\int \big(\partial_m U( \overline{y})-
\partial_m U(  e^{- t/2}x+(1-e^{-t})^{1/2}x^*)
\big)^2  d\lambda_x(y)\big)^{1/2} \nonumber\\
\leq &(1-e^{-t})^{1/2}\big(\int| \sum_{j=1}^n\partial_j\partial_m
U(\xi)(y_j-x^*_j)|^2
d\lambda_x(y)\big)^{1/2}\nonumber\\
\leq &  \big(\int |y-x^*|^2\sum_{j=1}^n\partial_j\partial_m
U(\xi)^2d\lambda_x(y)\big)^{1/2}\nonumber\\
\leq & C_2 \big(\int   |y-x^*|^2d\lambda_x(y)\big)^{1/2}\nonumber\\
\leq &C_2(3+C_2)^{1/2}\sqrt{n}  ,\label{point2}
\end{align}
where $\xi=e^{- t/2}x+(1-e^{-t})^{1/2} (x^*+t'(y-x^*))$ for some
$t'\in [0,1]$ depending on $x$ and $y$.   By   (\ref{moment}) and
(\ref{point2}),
\begin{align}
& \max_{1\leq m\leq n}\int \big|\partial_m U( \overline{y})-
\int\partial_m U( \overline{y})d\lambda_x(y) \big|
d\lambda_x(y)\nonumber\\\leq &  \max_{1\leq m\leq n}\big(\int
\big(\partial_m U( \overline{y})- \int\partial_m U(
\overline{y})d\lambda_x(y)\big)^2
d\lambda_x(y)\big)^{1/2} \nonumber\\
\leq & \max_{1\leq m\leq n} \big(\int \big(\partial_m U(
\overline{y})-
\partial_m U(
  e^{- t/2}x+(1-e^{-t})^{1/2}x^*)
\big)^2  d\lambda_x(y)\big)^{1/2} \nonumber\\\leq
&C_2(3+C_2)^{1/2}\sqrt{n} .\nonumber
\end{align}
For $l=1,2,3$, we have by the second inequality of (\ref{de,3;1})
and the estimate above
\begin{align} |M_{2,l}(x)|\leq 2C_2\max_{1\leq m\leq n}\int \big|\partial_m
U( \overline{y})- \int\partial_m U( \overline{y})\lambda_x(y) \big|
d\lambda_x(y)
 \leq &2C_2^2(3+C_2)^{1/2}\sqrt{n}   .\label{10n}
\end{align}

By the second inequality of  (\ref{de,3;1}), the  mean value theorem
 and (\ref{<>4}), we have
\begin{align} &\max_{1\leq m\leq n}\int \big|\partial_m U(
\overline{y})-  \partial_m U(
  e^{- t/2}x+(1-e^{-t})^{1/2}x^*)
\big|^3  d\lambda_x(y) \nonumber\\
\leq &  C_2^3  \int  |y-x^*|^3d\lambda_x(y)
\nonumber\\
\leq &  C_2^3(3+C_2)^{3/2} n^{3/2}  .\label{point3}
\end{align}
By (\ref{point2}) and   Cauchy-Schwartz inequality,
\begin{align} & \max_{1\leq m\leq n}\big|\partial_m U(
  e^{- t/2}x+(1-e^{-t})^{1/2}x^*)-
\int\partial_m U( \overline{y})d\lambda_x(y) \big|\leq
C_2(3+C_2)^{1/2}\sqrt{n} . \label{point4}
\end{align}
Applying   (\ref{point3}), (\ref{point4}), H\"{o}lder inequality
and inequality $(a+b)^3\leq 4(a^3+b^3)$ for $a,b\geq 0$, we have
\begin{align} |M_3(x)|\leq&
\max_{1\leq m\leq n} \int \big|\partial_m U( \overline{y})-
\int\partial_m U(
\overline{y})\lambda_x(y) \big|^3 d\lambda_x(y)  \nonumber\\
\leq &4\max_{1\leq m\leq n} \int \big|\partial_m U( \overline{y})-
\partial_m U(
  e^{- t/2}x+(1-e^{-t})^{1/2}x^*)\big|^3 d\lambda_x(y)
   \nonumber\\+ & 4\max_{1\leq m\leq n}  \big|\partial_m U(
  e^{- t/2}x+(1-e^{-t})^{1/2}x^*)-
\int\partial_m U( \overline{y})\lambda_x(y) \big|^3 \nonumber\\\leq
& 8C_2^3(3+C_2)^{3/2} n^{3/2} \label{I3}.
\end{align}
Applying (\ref{I123}), (\ref{I1}), (\ref{10n}), (\ref{I3}) and
$C_2\geq 1$, we get
\begin{align}
 |\partial_i\partial_j\partial_k
U_t(x)|  \leq e^{-3t/2}(C_3+C_2^2(6+8C_2)(3+C_2)^{3/2} n^{3/2}
),\nonumber
\end{align}
which implies      the conclusion. \qed\medskip

 \begin{lemma} \label{derivatives,111}
Let  $u=e^{-U} \in\mathcal{CF}_n$ satisfying condition
$\mathcal{L}(C_1,C_2,C_3)$ for some constants    $0\leq C_1<C_2$ and
$C_3>0$.
  Let    $1\leq i,j \leq n,x\in \mathbb{R}^n$ and $t,s\geq 0$.
   Then we have for every $n\geq N_2$  \begin{align}
 &|\partial_i\partial_jU_{t+s}(x)-
 \partial_i\partial_j U_t(e^{-
s/2}x+(1-e^{-s})^{1/2}x^*)|\nonumber\\ \leq & C_5
s^{1/2}e^{-t}n^{5/2},\label{i,j}
\end{align}
where    $U_t$ is defined by  (\ref{HA}), $x^* $ is defined  by
(\ref{hess}) and $C_5=C_2+C_2^2(3+C_2)+(3+C_2)^{1/2}C_4$ with $C_4$
specified in Lemma \ref{derivatives,ab11}.\end{lemma}
 \noindent{\bf Proof}\  Let $1\leq i,j \leq n,x\in \mathbb{R}^n$ and $t,s\geq 0$.
 Noticing that the left hand side of (\ref{i,j}) is zero when $s=0$,
 we   assume that $s>0$ in what below.
Define  $\nu_x=\nu_{u,t,s,x}$   by  (\ref{vts}). Applying
(\ref{J1J2}), we have
\begin{align}
&\partial_i\partial_jU_{t+s}(x) = e^{-s} J_1(x)-
e^{-s}J_2(x),\label{C0,13*}
\end{align}
where
\begin{align}
J_1(x) =& \int \partial_i\partial_j U_t(e^{-
s/2}x+(1-e^{-s})^{1/2}y) d\nu_x(y),\nonumber\\
 J_2(x)=& \int \big(\partial_i
U_t(e^{- s/2}x+(1-e^{-s})^{1/2}y) -J_3(x)\big)\big(\partial_j
U_t(e^{-
s/2}x+(1-e^{-s})^{1/2}y) -J_4(x)\big) d\nu_x(y),\nonumber\\
 J_3(x)=& \int \partial_i
U_t(e^{- s/2}x+(1-e^{-s})^{1/2}y)  d\nu_x(y),\ \ \ \ \ J_4(x)= \int
\partial_j
U_t(e^{- s/2}x+(1-e^{-s})^{1/2}y)  d\nu_x(y).\nonumber
\end{align}

We always assume that $n\geq N_2$ in what below.  Applying    mean
value theorem,
\begin{align}
 &|J_1(x) -\partial_i\partial_j U_t(e^{-
s/2}x+(1-e^{-s})^{1/2}x^*)|\nonumber\\
=& \big|\int \big(\partial_i\partial_j U_t(e^{-
s/2}x+(1-e^{-s})^{1/2}y) -
\partial_i\partial_j U_t(e^{-
s/2}x+(1-e^{-s})^{1/2}x^*)\big)d\nu_x(y)\big|\nonumber\\\leq &
(1-e^{-s})^{1/2} \int \big|
\sum_{k=1}^n\partial_k\partial_i\partial_j U_t(\xi)
(y_k-x_k^*)\big|d\nu_x(y) ,\nonumber
\end{align}
where $\xi= e^{- s/2}x+(1-e^{-s})^{1/2}( x^*+t'(y-x^*))$ for some
$t'\in [0,1]$ depending on $x$ and $y$. By (\ref{<>4}), (\ref{j3}),
Cauchy-Schwartz inequality together with  the estimate above,
\begin{align}
  |J_1(x) -\partial_i\partial_j U_t(e^{- s/2}x+(1-e^{-s})^{1/2}x^*)|
 \leq &C_4s^{1/2} e^{-3t/2}n^{3/2} \sum_{k=1}^n\int
|y_k-x^*_k|d\nu_x(y)  \nonumber\\
\leq &C_4s^{1/2}e^{-3t/2} n^2 \int |y-x^*|d\nu_x(y) \nonumber\\\leq
&(3+C_2)^{1/2}C_4s^{1/2}e^{-3t/2}n^{5/2} .\label{J1}
\end{align}

By mean value theorem, Cauchy-Schwartz inequality, (\ref{de,131})
and (\ref{<>4}), we also have
\begin{align}
J_5(x):= & \max_{1\leq l\leq n} \int \big(\partial_l U_t(e^{-
s/2}x+(1-e^{-s})^{1/2}y) -\partial_l U_t(e^{-
s/2}x+(1-e^{-s})^{1/2}x^*
 )\big)^2d\nu_x(y) \nonumber\\ \leq  &
(1-e^{-s})   \int|\sum_{k=1}^n \partial_k\partial_l
U_t(\xi') (y_k-x^*_k)|^2d\nu_x(y)  \nonumber\\
\leq  & C_2^2se^{-2t}   \int
 |y-x^*|^2d\nu_x(y) \nonumber\\\leq  & (3+C_2)C_2^2se^{-2t}n
\nonumber,
\end{align}where $\xi'= e^{- s/2}x+(1-e^{-s})^{1/2} (x^*+t^{''}(y-x^*))$ for some
$t^{''}\in [0,1]$ depending on $x$ and $y$.  Applying
Cauchy-Schwartz inequality,  the two estimates above and
(\ref{moment}), we get
\begin{align}
|J_2(x)|\leq & \Big(\int \big(\partial_i U_t(e^{-
s/2}x+(1-e^{-s})^{1/2}y)
-J_3(x)\big)^2d\nu_x(y)\Big)^{1/2}\nonumber\\
&\cdot\Big( \int \big(\partial_j U_t(e^{- s/2}x+(1-e^{-s})^{1/2}y)
-J_4(x)\big)^2
d\nu_x(y)\Big)^{1/2}\nonumber\\
\leq &J_5(x) \nonumber\\
\leq &(3+C_2)C_2^2se^{-2t}n .\label{L}
\end{align}
Combing (\ref{C0,13*}),  (\ref{J1}) and (\ref{L}),  we obtain
\begin{align}
 &| \partial_i\partial_jU_{t+s}(x)-e^{-s}
 \partial_i\partial_j U_t(e^{-
s/2}x+(1-e^{-s})^{1/2}x^*)|\nonumber\\ \leq& e^{-s}|J_1(x)
-\partial_i\partial_j U_t(e^{-
s/2}x+(1-e^{-s})^{1/2}x^*)|+e^{-s}|J_2(x)|\nonumber\\
\leq &(3+C_2)^{1/2}C_4s^{1/2}e^{-3t/2}n^{5/2}
+(3+C_2)C_2^2se^{-2t}n,\nonumber
\end{align}
which further gives together with  the second inequality of
(\ref{de,131})
\begin{align}
 &|\partial_i\partial_jU_{t+s}(x)-
 \partial_i\partial_j U_t(e^{-
s/2}x+(1-e^{-s})^{1/2}x^*)|\nonumber\\
\leq &(3+C_2)^{1/2}C_4s^{1/2}e^{-3t/2}n^{5/2} +
C_2^2(3+C_2)se^{-2t}n+(1-e^{-s})C_2e^{-t}\nonumber\\
\leq & \big(C_2+(3+C_2)C_2^2+(3+C_2)^{1/2}C_4
\big)s^{1/2}e^{-t}n^{5/2},\nonumber
\end{align}
which gives the conclusion.\qed\medskip

 \begin{lemma} \label{con}
Let  $u=e^{-U} \in\mathcal{CF}_n$ satisfying condition
$\mathcal{L}(C_1,C_2,C_3)$ for some constants    $0\leq C_1<C_2$ and
$C_3>0$.
 Set   for every $t\geq 0$   \begin{align}
  \Lambda_{t} =e^{t}\cdot \inf\{ \inf_{\mathbf{e}\in S_{n-1}} \langle \mathbf{e}\cdot \nabla^2 U_t(y),\mathbf{e}\rangle:y\in  \mathbb{R}^n
 \} \label{C0,1} ,\end{align}
   Suppose that $\Lambda_{0}\geq c>\Lambda_{t_1}$
  for some constants $c,t_1>0$. Set
  $t_0=\sup\{0\leq t\leq t_1: \Lambda_t\geq  c\}$. Then
  for every $n\geq N_2$ we have  $\Lambda_{t_0}=c$.  \end{lemma}
 \noindent{\bf Proof}\ Let $n\geq N_2$.
 Applying Lemma \ref{derivatives,111}, we   have for every   $t\geq 0$
 and every
 $s_0>0$
\begin{align}
 \inf_{s\in(0,s_0)}(\Lambda_{t+s}-\Lambda_t)\geq -C_5
s_0^{1/2}n^{5/2}.\label{uni}
\end{align}
By assumption $\Lambda_{0}\geq c>\Lambda_{t_1}$, the definition of
$t_0$ and (\ref{uni}), we have $\Lambda_{t_0}\geq c$. Suppose that
$\Lambda_{t_0}> c$, applying (\ref{uni}), we have
$\Lambda_{t_0+\varepsilon}> c$ when $\varepsilon$ is small enough,
which contradicts   the definition of $t_0$. Therefore,
$\Lambda_{t_0}= c$ holds. \qed\medskip
 \begin{remark}
The function  $\Lambda_t$ above is in fact  continuous. To this end,
by (\ref{uni}),  it is sufficient to verify that $\lim_{t\downarrow
t_0}\Lambda_{t}\leq \Lambda_{t_0}$ for every $t_0\geq 0$.  This can
be done by applying (\ref{C0,13*}) for $i=j$ and Lemma
\ref{derivatives,ab11}.
 \end{remark}

 \begin{proposition} \label{C1,2,3}
Let  $u=e^{-U} \in\mathcal{CF}_n$ satisfying condition
$\mathcal{L}(C_1,C_2,C_3)$ for some constants    $0<C_1<C_2$ and
$C_3>0$. Then for every  $t\geq 0$ and every  $x\in \mathbb{R}^n$
\begin{align}
   \nabla^2 U_t(x)\geq  C_6 e^{-t} I_n, \label{C0,0}\end{align} where
$C_6=\min(e^{-3}C_1,2^{-6}e^{-3})$ and
      $U_t$ is defined by  (\ref{HA}).
    \end{proposition}
 \noindent{\bf Proof}\ \
      We   assume that     $n\geq N_2$ in
what below. Otherwise,  we can   consider  the function
$u_k(\mathbf{x}):=e^{-U_k( \mathbf{x})}$  on $\mathbb{R}^{kn}$ for
some $k$ with   $kn\geq N_2$, where $ U_k(
\mathbf{x})=\sum_{i=1}^kU(x^{(i)})$ with
$\mathbf{x}=(x^{(1)},\cdots,x^{(k)})$,  $x^{(l)}\in \mathbb{R}^n$
for $1\leq l\leq k$.   Without loss of generality we also assume
that $C_2\geq 1$. This implies $C_4,C_5\geq1$ by their  definitions.

 Define $\Lambda_t$   for every  $t\geq 0$  by
(\ref{C0,1}).
 To prove the lemma, we claim that it is sufficient to  verify the following
 conclusion:
 for every  $t\geq 0$ and every $ s\in (0,10^{-3}(C_4^2+C_5^2)^{-2}n^{-14}C_6^4]$
 \begin{align}
\Lambda_{t+s} \geq  \Lambda_{t}(1-se^{-t}),\ \ \ \ if \ C_6\leq
\Lambda_t\leq e^3C_6.\label{C0,2}
 \end{align}
Assume that $\Lambda_{t_1}<
   e^3C_6$  for some $t_1\geq 0$.   To verify the  claim  above, in what below   we only need
to  show that $\Lambda_{t_1}\geq  C_6$
   under the assumption (\ref{C0,2}).

   From the   definition of $C_6$ and the assumption of $U$, we have $\Lambda_{0}  \geq
e^3C_6$. Define $t_0=\sup\{0\leq t \leq t_1: \Lambda_{t}\geq
e^3C_6\}$.
  Then, by $\Lambda_{0} \geq e^3C_6$,
$\Lambda_{t_1}<e^3C_6$ and Lemma \ref{con}, we have
\begin{align} 0\leq t_0<t_1,\ \ \ \ \Lambda_{t_0}=e^3C_6.
\label{t0}\end{align}
  Choose integer
$k'\geq 1$ and $s'\in (0,10^{-3}(C_4^2+C_5^2)^{-2}n^{-14}C_6^4]$
such that $t_0+k's'= t_1$.  Notice  that we have $ s'\leq 10^{-3}$
  by $C_4\geq 1$ and $C_6\leq 1$. Then, applying the method of  finite
induction,  we have by (\ref{C0,2}) and (\ref{t0})
 \begin{align}
\Lambda_{t_1}=\Lambda_{t_0+k's'} \geq&
\Lambda_{t_0}\prod_{j=1}^{k'}(1-s'e^{-(t_0+(j-1)s')})\nonumber\\
 \geq&
e^3C_6\exp\{-2s'\sum_{j=1}^{k'}e^{-(t_0+(j-1)s')}\}
\nonumber\\
 \geq&
e^3C_6\exp\{-2e^{s'} \int_0^\infty e^{-(t_0+r)}dr\}
\nonumber\\
\geq & C_6,\nonumber
\end{align}
where we use $\ln(1-c)\geq -2c$ for every $c\in (0,1/2)$ in the
second inequality above.    Therefore,   the lemma holds  if we can
verify  (\ref{C0,2}).

To prove   (\ref{C0,2}),  in what below we assume that $t\geq 0$ and
\begin{align}
 s\in& (0,10^{-3}(C_4^2+C_5^2)^{-2}n^{-14}C_6^4];\ \label{t,s,C6}\\\ \ \ \
  C_6&\leq \Lambda_{t}\leq e^3C_6.\label{t,s,C6,2}
\end{align}
Recall that $\nu_x$ is defined by $\nu_{u,t,s,x}$ in (\ref{vts}).
Fix an arbitrary element
  $x\in \mathbb{R}^n$ and define  $x^*\in \mathbb{R}^n$ by
(\ref{hess}). Choose a coordinate system such that $\nabla^2
U_{t+s}(x)$ is a diagonal matrix. Let $1\leq i \leq n$.
  To prove (\ref{C0,2}), it is sufficient to  show that
\begin{align}
\partial_i^2U_{t+s}(x) \geq e^{-(t+s)}\Lambda_{t}(1-e^{-t}s).\label{t,s}
\end{align}  By  definition (\ref{C0,1})
and  the first inequality     of   (\ref{t,s,C6,2}),
\begin{align}
\partial_i^2U_{t}(e^{-
s/2}x+(1-e^{-s})^{1/2}x^*) \geq  C_6e^{-t}.\label{t,s,i}
\end{align}
From the assumption that   $\nabla^2 U_{t+s}(x)$ is  diagonal, we
have
\begin{align}
   \partial_j\partial_k U_{t+s}(x)=0
   ,\ \ \ \ \ for\ 1\leq j<k\leq n\nonumber .\end{align}
If $1\leq j<  k \leq n$,  by (\ref{i,j}) and the equality  above, we
obtain
\begin{align} &|
 \partial_j\partial_k U_{t}(e^{-
s/2}x+(1-e^{-s})^{1/2}x^*)|\leq
C_5s^{1/2}e^{-t}n^{5/2}\label{j3,}.\end{align}

Applying (\ref{J1J2}), we have
\begin{align}
\partial_i^2U_{t+s}(x)  =& e^{-s} L_1(x)-
e^{-s}L_2(x),\label{C0,3}
\end{align}
where
\begin{align}
L_1 (x)=& \int \partial_i^2 U_{t}(e^{-
s/2}x+(1-e^{-s})^{1/2}y) d\nu_x(y),\nonumber\\
 L_2(x)=& \int \big(\partial_i
U_{t}(e^{- s/2}x+(1-e^{-s})^{1/2}y) -L_3(x)\big)^2
d\nu_x(y),\nonumber
\\L_3(x)=&\int\partial_i U_{t}(e^{-
s/2}x+(1-e^{-s})^{1/2}y) d\nu_x(y).\nonumber
\end{align}

Since $\nu_x$ is a probability measure, with  assumption
(\ref{t,s,C6}) we have by    (\ref{<>21}), (\ref{j3}),
   mean value theorem and Cauchy-Schwartz inequality
\begin{align}
& \big|L_1(x)-\partial_i^2 U_t( e^{-
s/2}x+(1-e^{-s})^{1/2}x^*) \big|\nonumber\\
 =&  \big|\int \big(\partial_i^2U_t(e^{-
s/2}x+(1-e^{-s})^{1/2}y) -\partial_i^2 U_t(e^{-
s/2}x+(1-e^{-s})^{1/2}x^*)\big)d\nu_x(y) \big|\nonumber\\\leq &
C_4(1-e^{-s})^{1/2}e^{-3t/2} n^{3/2}
 \big(\int   \sum_{k=1}^n |y_k-x^*_k| d\nu_x(y)\big)
 \nonumber\\ \leq &
C_4s^{1/2}e^{-3t/2} n^{2}
 \int    |y-x^*|  d\nu_x(y)
 \nonumber\\ \leq &
 \sqrt{2}s^{1/2}C_4e^{-t} n^{5/2}\nonumber\\
 \leq & 2^{-1}C_6e^{-t}.   \nonumber
\end{align}
The estimate above and (\ref{t,s,i}) give
\begin{align}
L_1(x)&\geq \frac{1}{2}\partial_i^2 U_t( e^{-
s/2}x+(1-e^{-s})^{1/2}x^*). \label{C0,11}
\end{align}
By (\ref{de,131}) and the definitions of  $\Lambda_t$ and $L_1(x)$,
we also have
\begin{align}
   e^{-t}\Lambda_t\leq L_1 (x)\leq C_2e^{-t}. \label{C0,12}
\end{align}

Applying Talor formula and (\ref{j3}), we have for every $y\in
\mathbb{R}^n$

\begin{align}
 & \big|\partial_i
U_t(e^{- s/2}x+(1-e^{-s})^{1/2}y) -\partial_iU_t(e^{-
s/2}x+(1-e^{-s})^{1/2}x^* ) \nonumber\\- &  (1-e^{-s})^{1/2}
\sum_{j=1}^n\partial_i\partial_jU_t(e^{- s/2}x+(1-e^{-s})^{1/2} x^*
)(y_j-x_j^*)\big|\nonumber\\
=&\frac{1}{2}(1-e^{-s})
|\sum_{j=1}^n\sum_{k=1}^n\partial_i\partial_j\partial_kU_t(\xi
)(y_j-x_j^*)(y_k-x_k^*)|\nonumber\\
 \leq& C_4se^{-3t/2}n^{5/2}|y -x^*|^2. \label{C0,8}
\end{align}
where $\xi=e^{- s/2}x+(1-e^{-s})^{1/2} (x^*+t'(y-x^*))$ for some
$t'\in [0,1]$ depending on $x$ and $y$.  Applying (\ref{<>20}),
(\ref{C0,8}) together with the inequality $(a+b)^2\leq 2(a^2+b^2)$
for every  $a,b\in \mathbb{R}$
\begin{align}
 L_4(x):= &\int \big(\partial_iU_t(e^{-
s/2}x+(1-e^{-s})^{1/2}y) -\partial_iU_t(e^{-
s/2}x+(1-e^{-s})^{1/2}x^* ) \big)^2 d\nu_x(y)\nonumber\\
\leq  & 2\int \big(  (1-e^{-s})^{1/2}
\sum_{j=1}^n\partial_i\partial_jU_t(e^{- s/2}x+(1-e^{-s})^{1/2} x^*
)(y_j-x_j^*) \big)^2 d\nu_x(y)\nonumber\\+&2\int \Big(\partial_i
U_t(e^{- s/2}x+(1-e^{-s})^{1/2}y) -\partial_iU_t(e^{-
s/2}x+(1-e^{-s})^{1/2}x^*) \nonumber\\&\ \ \ \ \ \ \
-(1-e^{-s})^{1/2} \sum_{j=1}^n\partial_i\partial_jU_t(e^{-
s/2}x+(1-e^{-s})^{1/2} x^* )(y_j-x_j^*) \Big)^2
d\nu_x(y)\nonumber\\
\leq &2(1-e^{-s})\sum_{j=1}^n\sum_{k=1}^na_{i,j}a_{i,k}  \int
(y_j-x_j^*)(y_k-x_k^*) d\nu_x(y)
+2C_4^2s^2e^{-3t} n^5\int |y -x^*|^4 d\nu_x(y)\nonumber\\
\leq &2s\sum_{j=1}^n\sum_{k=1}^n|a_{i,j}a_{i,k}| \big( \int
(y_j-x_j^*)^2d\nu_x(y)\big)^{1/2}
 \big( \int
(y_k-x_k^*)^2d\nu_x(y)\big)^{1/2} +12C_4^2s^2e^{-3t}n^7 , \nonumber
\end{align}
where \begin{align} a_{i,j}=\partial_i\partial_j U_{a,t}(e^{-
s/2}x+(1-e^{-s})^{1/2}x^* ),\ \ \ for\ 1\leq j\leq n.\nonumber
\end{align}
 Applying   (\ref{de,131}), (\ref{j3,}) and $C_2\leq C_5$ to the estimate above, we further get  \begin{align}
 L_4(x)\leq& 4s\sum_{j=1}^n\sum_{k=1}^n|a_{i,j}a_{i,k}|
  +12C_4^2s^2e^{-3t}n^7 \nonumber\\ =& 4sa_{i,i}^2+8s\sum_{j=1,j\neq i}^n|a_{i,i}a_{i,j}|
  +4s\sum_{j=1,j\neq i}^n\sum_{k=1,k\neq i}^n|a_{i,j}a_{i,k}|
   +12C_4^2s^2e^{-3t}n^7 \nonumber\\
  \leq &4sa_{i,i}^2+8s(n-1)C_2e^{-t}\cdot C_5s^{1/2}n^{5/2}e^{-t}+
  4s(n-1)^2(C_5s^{1/2}e^{-t}n^{5/2})^2 +12C_4^2s^2e^{-3t}n^7 \nonumber\\\leq&
4sa_{i,i}^2+12(C_4^2+C_5^2)s^{3/2}  e^{-2t}n^7  . \nonumber
\end{align}
The estimate above and (\ref{moment}) give
\begin{align}
L_2(x)\leq L_4 (x)\leq 4sa_{i,i}^2+12(C_4^2+C_5^2)s^{3/2} e^{-2t}n^7
 .\label{L2}
\end{align}

Applying (\ref{C0,3}), (\ref{C0,11}) and (\ref{L2}), we have by
assumption $s\in(0,10^{-3}(C_4^2+C_5^2)^{-2}n^{-14}C_6^4]$
\begin{align}
\partial_i^2U_{t+s}(x) \geq & e^{-s} L_1(x)-
e^{-s}(4sa_{i,i}^2+12(C_4^2+C_5^2)s^{3/2}  e^{-2t}n^7
)\nonumber\\
 \geq & e^{-s}( L_1(x)-
2^4sL_1(x)^2-C_6^2 se^{-2t}) .\nonumber
\end{align}
By  the first inequality of (\ref{C0,12}) and $\Lambda_t\geq C_6$ in
(\ref{t,s,C6,2}), we have   $L_1(x)\geq C_6 e^{-t} $  . This and the
estimate above  further give
\begin{align}
\partial_i^2U_{t+s}(x) \geq & e^{-s}( L_1(x)-
2^5sL_1(x)^2) .\label{px}
\end{align}

When  $L_1(x)\geq 2e^3C_6e^{-t}$, we have by
    (\ref{t,s,C6}), the second inequality of
  (\ref{C0,12}) and  (\ref{px})
\begin{align}
\partial_i^2U_{t+s}(x) \geq   e^{-s}(2e^3C_6e^{-t}-
2^5sC_2^2e^{-2t})> e^3C_6e^{-t-s},\nonumber
\end{align}
which gives     (\ref{t,s})  together with    the assumption
$\Lambda_t\leq e^3C_6$ in (\ref{t,s,C6,2}). Next   we assume that
$L_1(x)< 2e^3C_6e^{-t}$. By (\ref{px})  and $C_6\leq 2^{-6}e^{-3}$,
we have
\begin{align}
\partial_i^2U_{t+s}(x) >   e^{-s}L_1(x)( 1-
2^6e^3C_6se^{-t})\geq  e^{-s}L_1(x)( 1- e^{-t}s),\nonumber
\end{align}
which further  gives by the first inequality of (\ref{C0,12})
\begin{align}
\partial_i^2U_{t+s}(x) \geq    e^{-t-s}\Lambda_t( 1-
e^{-t}s). \nonumber
\end{align}
Combing the two cases above we  completes the proof of (\ref{t,s}).
  \qed\medskip

\subsection{\normalsize some
log-concave functions associated with symmetric convex sets}

 For $A\subseteq
\mathbb{R}^n$, define for every $x\in \mathbb{R}^n$
\begin{align}
\rho_A(x)=&\inf\{|x-y|: y\in A\}.\label{rho}
\end{align}
Let  $\alpha >0$. Define for every $A\in \mathcal{C}_n$ and every
$x\in \mathbb{R}^n$
\begin{align}
H_{A,\alpha}(x)=&\frac{\alpha}{2}|x|^2+  n\rho_{A }(x), \label{hA}\\
h_{A,\alpha}(x)=&\exp\{-H_{A,\alpha}(x)\}.\label{HA'}
\end{align}
We see that $ h_{A,\alpha}\in \mathcal{CF}_n$.  For every  $t\geq
0$, define functions   $ h_{A,\alpha,t}$ and $H_{A,\alpha,t}$
  as follows:
\begin{align}
h_{A,\alpha,t}=\exp\{- H_{A,\alpha,t}\}= P_th_{A,\alpha}.
\label{HA''}
\end{align}
Notice that  $h_{A,\alpha,0}=h_{A,\alpha}$ and
$H_{A,\alpha,0}=H_{A,\alpha}$  from the definitions above. We may
prove the result of Proposition \ref{C1,2,3} under a more natural
assumption  that  $u$ satisfies condition $\nabla^2 U\geq CI_n$ on
$\mathbb{R}^n$ for some constant $C>0$. Next we only prove this for
some special cases which is enough for our purpose.

 \begin{lemma} \label{C1,2,3;}
Let $\alpha>0$ and $A\in \mathcal{C}_n$.  Then for every   $x\in
\mathbb{R}^n$ and every $t>0 $
\begin{align}
C(\alpha)e^{-t} I_n\leq   \nabla^2 H_{A,\alpha,t}(x)\leq 2(1\wedge
t)^{-1}e^{-t} I_n,\label{C0,0}\end{align} where
$C(\alpha)=\min(e^{-3}\alpha ,2^{-6}e^{-3})$.
 \end{lemma}
  \noindent{\bf Proof}\
  Let $\beta\in(0,1)$ and define for every $y\in \mathbb{R}^n$
 \begin{align}
H_{A,\alpha}^{(\beta)}(y)=P_\beta H_{A,\alpha}(y).\nonumber
 \end{align}
 Set  $H(y)=2^{-1}\alpha|y|^2$ for
$y\in\mathbb{R}^n$.  Let $x\in\mathbb{R}^n$. We have by definition
 \begin{align}
\nabla^2 H_{A,\alpha}^{(\beta)}(x)\geq \nabla^2 P_\beta
H(x)=e^{-\beta}\alpha I_n.\label{bet}
 \end{align}
   Let  $1\leq i,j,k\leq n$.
Notice that  $\rho_A$ is a Lipschitz function and $|
\nabla\rho_A|\leq 1$ almost everywhere. Then, we have
 \begin{align}
&\partial_i P_\beta \rho_A (x) = e^{-\beta/2} \int
 \partial_i  \rho_A (e^{-\beta/2}x+(1-e^{-\beta})^{1/2}y) d\mu_n(y),\nonumber\end{align}
which further gives
\begin{align}
\partial_i^2 P_\beta \rho_{A}(x) =& \frac{
e^{-\beta}}{(2\pi)^{n/2}(1-e^{-\beta})^{(n+2)/2}} \int
(y_i-e^{-\beta/2}x_i)
  \partial_i  \rho_A (y)
\exp\{-\frac{|y-e^{-\beta/2}x|^2}{2(1-e^{-\beta})}\}dy \nonumber\\
\leq & \frac{ e^{-\beta}}{(2\pi)^{n/2}(1-e^{-\beta})^{(n+2)/2} }\int
 |y_i-e^{-\beta/2}x_i|
\exp\{-\frac{|y-e^{-\beta/2}x|^2}{2(1-e^{-\beta})}\}dy\nonumber\\
\leq &(1-e^{-\beta})^{-1/2}.\nonumber\end{align} By (\ref{hA}), the
equality in  (\ref{bet}) and the estimate above
\begin{align}
\partial_i^2 H_{A,\alpha}^{(\beta)}(x) \leq
& n(1-e^{-\beta})^{-1/2}+\alpha:=c(n,\alpha, \beta)
.\nonumber\end{align} Since $P_\beta H_{A,\alpha}$ is a convex
function and the estimate above holds under    any  coordinate
system
 $(Q(\mathbf{e}_i))_{1\leq i \leq n}$ when $Q$ is an  orthogonal
 transformation of $\mathbb{R}^n$, we get
 \begin{align}
 |\partial_j\partial_k H_{A,\alpha}^{(\beta)}(x)|\leq
 c(n,\alpha, \beta). \label{beta}
 \end{align}
By the semigroup property of $(P_t)$, we have
 \begin{align}
& \partial_i \partial_j\partial_k
H_{A,\alpha}^{(\beta)}(x)\nonumber\\
 =&\partial_i \partial_j\partial_k P_{\beta/2}H_{A,\alpha}^{(\beta/2)}(x)\nonumber\\
 =& \frac{e^{-\beta/2}}{(2\pi(1-e^{-\beta/2}))^{n/2}}\partial_i \int \partial_j\partial_k
  H_{A,\alpha}^{(\beta/2)}(e^{-\beta/4}x+(1-e^{-\beta/2})^{1/2}y)\exp\{-\frac{|y|^2}{2}\}dy \nonumber\\
 =& \frac{e^{-3\beta/4}(1-e^{-\beta/2})^{-1}}{(2\pi(1-e^{-\beta/2}))^{n/2}} \int
  (y_i-e^{-\beta/4}x_i )\partial_j\partial_k
 H_{A,\alpha}^{(\beta/2)}( y)\exp\{-\frac{|y-e^{-\beta/4}x|^2}{2(1-e^{-\beta/2})}\}dy, \nonumber
 \end{align}
 which implies  together with  (\ref{beta})
  \begin{align}
&| \partial_i \partial_j\partial_k
H_{A,\alpha}^{(\beta)}(x)|\nonumber\\
\leq & \frac{ (1-e^{-\beta/2})^{-1}c(n,\alpha,\beta/2)}{
(2\pi(1-e^{-\beta/2}))^{n/2}} \int
  |y_i-e^{-\beta/4}x_i | \exp\{-\frac{|y-e^{-\beta/4}x|^2}{2(1-e^{-\beta/2})}\}dy\nonumber\\
\leq & \frac{c(n,\alpha,\beta/2)}{ (1-e^{-\beta/2})^{1/2} }
.\nonumber
 \end{align}
 From  (\ref{bet}), (\ref{beta}) and  the estimates obtained above, we see that
 $H_{A,\alpha}^{(\beta)}$ satisfies condition $\mathcal{L}(
 C_1,C_2,C_3)$ with $
 C_1=
 e^{-\beta}\alpha,C_2=c(n,\alpha, \beta)$ and $C_3=(1-e^{-\beta/2})^{-1/2}c(n,\alpha,\beta/2)$.

  Let  $t> 0$ and  define
 $H_{A,\alpha,t}^{(\beta)}$
 by $\exp\{-H_{A,\alpha,t}^{(\beta)}\}=P_t\exp\{-H_{A,\alpha}^{(\beta)}\}$.
 Applying
 Proposition  \ref{C1,2,3}, we get
 \begin{align}
\nabla^2 H_{A,\alpha,t}^{(\beta)}(x)\geq \min(e^{-3}e^{-\beta}\alpha
,2^{-6}e^{-3})  e^{-t}I_n .\label{limit}
 \end{align}
Applying (\ref{J1J2}), we have
\begin{align}
\lim_{\beta\rightarrow
0}\partial_i\partial_jH_{A,\alpha,t}^{(\beta)}(x)
=\partial_i\partial_jH_{A,\alpha,t} (x) .\nonumber\end{align} Then,
we get
 the first inequality of
 $(\ref{C0,0})$ by taking $\beta \rightarrow 0$ in (\ref{limit}). The second inequality of (\ref{C0,0}) follows by Lemma
 \ref{derivatives,2}.\qed\medskip

\section{Derivative estimates for time parameter }

\subsection{\normalsize second derivative  estimates    for large time}

 Denote by $Q^\tau$  the transpose  of a matrix $Q$. For every
function $u$ on $\mathbb{R}^n$ and every $t\geq 0$, denote
\begin{align}P_t\nabla u=&(P_t\partial_1 u,\cdots,P_t\partial_n u),
\ \ \ \ \ P_t\nabla^2
u=(P_t\partial_i\partial_j u )_{1\leq i,j\leq n},\nonumber\\
\nabla^\tau u=&(\nabla u )^\tau ,\ \ \ \  \ \ \ \  \ \ \ \ \ \ \ \ \
\nabla^\tau(P_t\nabla u)=(
\partial_i( P_t
\partial_j u) )_{1\leq i,j\leq n},\nonumber
\end{align}
provided that    the right hand sides above are well defined. By
 (\ref{OU}), for   smooth function  $u$  with
  gradient  controlled by some polynomial for instance,
\begin{align}\nabla P_t
u=e^{-t/2}P_t\nabla  u,\ \ \ \ \ \forall\  t\geq 0. \label{gred}
\end{align}

\begin{lemma} \label{secd,0}Let  $u$ and $v$ be smooth functions
on $\mathbb{R}^n$ with bounded second derivatives.
 Then for every  $t\geq 0$
\begin{align}
\frac{d^2}{dt^2}\phi_t(u,v)=  - \frac{1}{2}
\frac{d}{dt}\phi_t(u,v)+\frac{1}{4} \int trace\Big(  \nabla^2P_{t/2}
u \cdot \nabla^2P_{t/2}v\Big)
  d\mu_n.\label{SAB2}
\end{align}

 \end{lemma}
 \noindent{\bf Proof}\
 Applying (\ref{d}), (\ref{t/2}) and
   (\ref{gred}), we have
  \begin{align}
  \frac{d}{dt}\phi_t(u,v) =&\frac{d}{dt} \big( \int    P_{t/2} u
  P_{t/2}v
  d\mu_n\big)\nonumber\\
   =&- \frac{1}{2}\int   \langle\nabla P_{t/2} u ,
 \nabla P_{t/2}v \rangle
  d\mu_n
  \nonumber\\  =&- \frac{1}{2}e^{-t/2}\int   \langle P_{t/2} \nabla  u ,
 P_{t/2}\nabla v \rangle
  d\mu_n ,\nonumber
        \end{align}
        and hence
 \begin{align}
  \frac{d^2}{dt^2}\phi_t(u,v)
   =&-\frac{1}{2}\frac{d}{dt}\big( e^{-t/2}\int   \langle P_{t/2} \nabla  u ,
 P_{t/2}\nabla v \rangle
  d\mu_n \big)
  \nonumber\\
   =& \frac{1}{4} e^{-t/2}\int   \langle P_{t/2} \nabla  u ,
 P_{t/2}\nabla v\rangle
  d\mu_n+\frac{1}{4}e^{-t/2}\int trace\Big(  \nabla^\tau (P_{t/2}\nabla u)
\cdot  \nabla^\tau (P_{t/2}\nabla v)\Big)
  d\mu_n\nonumber\\
  =&-\frac{1}{2}
\frac{d}{dt}\phi_t(u,v)+ \frac{1}{4} \int trace\Big( \nabla^2P_{t/2}
u\cdot \nabla^2P_{t/2}v\Big)
  d\mu_n,\nonumber
        \end{align}
 which gives (\ref{SAB2}).
\qed\medskip

\begin{lemma} \label{low}  Let  $u=e^{-U}  $ be a smooth function of $\mathcal{CF}_n$
and assume that
    $u$
satisfies   condition $\mathcal{L}(C_1,C_2)$ for some constants
 $0<C_1<C_2$. Then
  \begin{align}
C_1|x|\leq |\nabla U(x)|\leq & C_2|x|
 ,\ \ \ \ \ \forall\  x\in \mathbb{R}^n,\label{fi}\\ \langle \frac{\nabla U(x)}
 {|\nabla U(x)|},
  \frac{x}{|x|}\rangle\geq&
 \frac{C_1}{C_2},\ \ \ \ \ \ \  \forall\  x\in \mathbb{R}^n\ with \ x\neq 0 .\label{trx} \end{align}

 \end{lemma}\noindent{\bf Proof}\ Let $x\in \mathbb{R}^n$.
We have $\nabla U(0)
 =0$  by the symmetric assumption of $U$. This gives
(\ref{fi}) when $x=0$. Next we assume   $x\neq 0$ and denote
$\mathbf{e}= \frac{x}{|x|}$.  By $\nabla U(0)
 =0$ and  the assumption  $\mathcal{L}(C_1,C_2)$ of $u$,    \begin{align}
   \partial_{\mathbf{e}}U(x)= \langle \nabla U(x) ,
  \mathbf{e}\rangle \geq C_1 |x|,\label{U}\end{align}
 which implies the first inequality of  (\ref{fi}). For every  $\textbf{e}'\in
 S_{n-1}$,  by $\nabla U(0)
 =0$ and  the assumption  $\mathcal{L}(C_1,C_2)$ of $u$  we also have
 \begin{align}
   \partial_{\mathbf{e}'}U(x) \leq C_2 |x|,\nonumber\end{align}
   which implies the second  inequality of  (\ref{fi}).
Applying (\ref{U}) and  the second inequality of (\ref{fi}),   we
get (\ref{trx}). \qed\medskip

\begin{lemma} \label{low}  Let
 $\alpha>0 $ and $A,B\in \mathcal{C}_n$.  Then   for every  $t>  4(2\ln 2-\ln C(\alpha))$
 \begin{align}
\frac{d^2}{dt^2}\phi_t(h_{A,\alpha},h_{B,\alpha})
>& -\frac{1}{2}
\frac{d}{dt}\phi_t(h_{A,\alpha},h_{B,\alpha})  ,\nonumber
\end{align}
where $C(\alpha)=\min(e^{-3}\alpha,2^{-6}e^{-3})$.
 \end{lemma}
 \noindent{\bf Proof}\ Let  $t>  0$. Recall that $h_{A,\alpha}$ and
$h_{A,\alpha,t}$ are defined by (\ref{HA'}) and (\ref{HA''})
respectively. We have by  (\ref{SAB2})
 \begin{align}
\frac{d^2}{dt^2}\phi_t(h_{A,\alpha},h_{B,\alpha})  =& -\frac{1}{2}
\frac{d}{dt}\phi_t(h_{A,\alpha},h_{B,\alpha})+\frac{1}{4} \int
trace\Big( \nabla^2 h_{A,\alpha,t/2} \cdot \nabla^2 h_{B,\alpha,t/2}
\Big)
  d\mu_n.\label{SAB1**}
\end{align}
Recall that  $H_{A,\alpha}$ and $H_{A,\alpha,t}$ are defined by
(\ref{hA}) and (\ref{HA''}), respectively. Direct calculation shows
that for every $x\in \mathbb{R}^n$
 \begin{align}
  &trace\Big( \nabla^2 h_{A,\alpha,t/2} \cdot \nabla^2
h_{B,\alpha,t/2}
\Big)(x)\nonumber\\
=&h_{A,\alpha,t/2}(x)h_{B,\alpha,t/2}(x)\big(
K_{1,t}(x)+K_{2,t}(x)-K_{3,t}(x)-K_{4,t}(x)\big) ,\label{SAB1}
\end{align}
where
  \begin{align}
  K_{1,t}(x)=& trace \big( \nabla^2  H_{A,\alpha,t/2}
  \cdot  \nabla^2 H_{B,\alpha,t/2}   \big)(x)  ,\nonumber\\
  K_{2,t} (x)=&trace\Big( \big(  \nabla^\tau H_{A,\alpha,t/2} \cdot \nabla
 H_{A,\alpha,t/2}
\big)\cdot \big( \nabla^\tau H_{B,\alpha,t/2}  \cdot \nabla
H_{B,\alpha,t/2} \big) \Big)(x)
 ,\nonumber\\
 K_{3,t}(x)=&trace\big(\nabla^2 H_{A,\alpha,t/2}   \cdot
   \nabla^\tau H_{B,\alpha,t/2}  \cdot \nabla H_{B,\alpha,t/2}
 \big)(x),\nonumber\\
K_{4,t}(x)=& trace\big(\nabla^2 H_{B,\alpha,t/2}  \cdot \nabla^\tau
H_{A,\alpha,t/2} \cdot \nabla H_{A,\alpha,t/2} \big)(x). \nonumber
\end{align}

  Applying  (\ref{t/2}) and   Lemma
\ref{C1,2,3;}, we have
\begin{align}
 \int K_{1,t} h_{A,\alpha,t/2}h_{B,\alpha,t/2}d\mu_n
 \geq&  C(\alpha)^2 e^{-t} n\phi_t(h_{A,\alpha},h_{B,\alpha}).\label{K4}\end{align}
Since the product of two suitable matrixes can be switched under
trace operation, we have for every $x\in
     \mathbb{R}^n$
\begin{align}
   K_{2,t} (x)=&trace\Big( \big(  \nabla^\tau H_{A,\alpha,t/2} \cdot \nabla
  H_{A,\alpha,t/2}
\big)\cdot \big( \nabla^\tau H_{B,\alpha,t/2} \cdot \nabla
H_{B,\alpha,t/2}\big)  \Big)(x)
 \nonumber\\=& trace \big(    \nabla
  H_{A,\alpha,t/2}
\cdot \nabla^\tau H_{B,\alpha,t/2} \cdot \nabla
H_{B,\alpha,t/2}\cdot\nabla^\tau H_{A,\alpha,t/2} \big)(x)
 \nonumber\\
 =&|\langle \nabla
  H_{A,\alpha,t/2}(x), \nabla
  H_{B,\alpha,t/2}(x)\rangle|^2, \nonumber
\end{align}
which gives
\begin{align}
 &\int K_{2,t}
 h_{A,\alpha,t/2}h_{B,\alpha,t/2}d\mu_n\geq0.\label{K2,}\end{align}
 By
(\ref{de,3}) and the second inequality of (\ref{fi}), we have for
every  $x\in
     \mathbb{R}^n$
 \begin{align}
 K_{3,t}(x)=&trace\big(\nabla^2 H_{A,\alpha,t/2}  \cdot
   \nabla^\tau H_{B,\alpha,t/2} \cdot \nabla H_{B,\alpha,t/2}
 \big)(x)\nonumber\\
 =& \big(\nabla H_{B,\alpha,t/2}\cdot\nabla^2 H_{A,\alpha,t/2}  \cdot
   \nabla^\tau H_{B,\alpha,t/2}
     \big)(x)\nonumber\\\leq &8(1\wedge\frac{t}{2})^{-3} e^{-3t/2} |x|^2.\label{K3}
\end{align}By (\ref{21}) and  (\ref{K3}), we have for every   $t\geq 2$
 \begin{align}
 &\int K_{3,t}
 h_{A,\alpha,t/2}h_{B,\alpha,t/2}d\mu_n\nonumber\\
 \leq&
8(1\wedge\frac{t}{2})^{-3} e^{-3t/2} \int |x|^2d\mu_n
 \int   h_{A,\alpha,t/2}h_{B,\alpha,t/2}d\mu_n
\nonumber\\ =& 8
 e^{-3t/2}n\phi_t(h_{A,\alpha},h_{B,\alpha}).\label{K,3}
\end{align}
 Similarly, we have  for every  $t\geq 2$\begin{align}
 \int K_{4,t} h_{A,\alpha,t/2}h_{B,\alpha,t/2}d\mu_n
 \leq&  8
 e^{-3t/2}n\phi_t(h_{A,\alpha},h_{B,\alpha}).\label{K4}\end{align}

Applying   (\ref{SAB1})-(\ref{K2,}), (\ref{K,3}) and (\ref{K4}), we
have for every   $t\geq 2$
 \begin{align}
  &\int trace\Big( \nabla^2 h_{A,\alpha,t/2 } \cdot \nabla^2
h_{B,\alpha,t/2} \Big) d\mu_n\geq e^{-t}\big( C(\alpha)^2
-2^4e^{-t/2} \big)n\phi_t(h_{A,\alpha},
  h_{B,\alpha})  .\nonumber
\end{align}
Since  $ 4(2\ln 2-\ln C(\alpha))>2 $, the estimate above shows that
for  every $t> 4(2\ln 2-\ln C(\alpha)) $
\begin{align}
  &\int trace\Big( \nabla^2 h_{A,\alpha,t/2} \cdot \nabla^2
h_{B,\alpha,t/2} \Big) d\mu_n > 0,\nonumber
\end{align}
which gives the conclusion with the help of
(\ref{SAB1**}).\qed\medskip

\subsection{\normalsize   derivative estimates    for moderate  time}
Let    $U$ be  a convex function   on $\mathbb{R}^n$. Define
probability measure $\nu$ on $\mathbb{R}^n$ as follows:\begin{align}
d\nu(y)=\big(\int\exp\{-U(y)\}d\mu_n\big)^{-1}\exp\{-U(y)\}d\mu_n.\label{nup}
 \end{align}
 For any   set $A\subseteq \mathbb{R}^n$ and $r>0$,
 denote
 \begin{align}
A[r]=&\{x:\rho_A(x)\leq r\},\ \ \ \forall r>
 0,\label{the}
\end{align}
where  $\rho_A$ is defined by (\ref{rho}). Define
\begin{align}\Phi(r)=\mu_1((-\infty,r]),\ \ \ \ \ \
\forall r\in \mathbb{R}.\end{align} Next we introduce  the
Poincar\'{e} inequality and the isoperimetric inequality of $\nu$.

The following    inequality of $\nu$ is a consequence of Theorem 4.1
of \cite{BL76}.

\textbf{Poincar\'{e} inequality of $\nu$}: \emph{For every
differential function $H $ on $ \mathbb{R}^n $ such that $ |\nabla
H|$ controlled by some polynomial,
\begin{align}
 \int \big( H(x)-\int H(x) d\nu \big)^2 d\nu \leq   \int |\nabla H|^2 d\nu.\label{Poinc}
 \end{align}}

 The following     inequality of $\nu$ is a consequence
 of
Corollary 2.2 in     \cite{BL96}. The formulation below   is taken
from Theorem 1.1 in  \cite{Le99}.

\textbf{Isoperimetric inequality of $\nu$}: \emph{For every
measurable set $A\subseteq \mathbb{R}^n$ and every $a\in
\mathbb{R}$, we have
\begin{align}
\nu(A[r])\geq \Phi(a+r) , \ \ \ \ \ \ \forall r\geq 0,\label{Iso}
\end{align}
provided that $\nu(A)\geq \Phi(a)$. }

 \begin{lemma} \label{concern} Let $\nu$
 be the probability measure on   $\mathbb{R}^n$ defined by (\ref{nup})
 for some symmetric convex function  $U$   on $\mathbb{R}^n$.
Let   $C>0$ and assume that  $H$ is a differentiable function on
$\mathbb{R}^n$ satisfying    $|\nabla H(x)|\leq C|x|$   for every
$x\in \mathbb{R}^n$. Then
 \begin{align}
 \int \big( H-\int H d\nu \big)^2 d\nu \leq
C^2n. \label{popp}
 \end{align}
 Moreover, for every $a_0>0$ there exists some integer
 $N_4=N_4(a_0,C)$ such that   \begin{align}
 \nu(x: H(x)-\int Hd\nu \geq an)\leq
 \exp\{- 2^{-7}C^{-2}a^2n\}+e^{-n/2},\ \ \ \ \ \forall a\geq a_0, \forall  n\geq N_4.\label{inc4}
 \end{align}
 \end{lemma} \noindent{\bf Proof}\
Since $e^{-U}\in \mathcal{CF}_n$, we have  by (\ref{21})
\begin{align}
\int |x|^2d\nu\leq  n.\label{two}
 \end{align}
Applying the estimate above and   the    Poincar\'{e} inequality in
(\ref{Poinc}), we get
\begin{align}
 \int \big( H-\int H d\nu \big)^2 d\nu \leq   \int |\nabla H|^2 d\nu
  \leq &C^2\int |x|^2 d\nu
  \leq C^2n,\nonumber
 \end{align}
which gives the first  conclusion.

Let $a\geq a_0$. To prove the second conclusion, without loss of
generality, we assume that $\int Hd\nu =0$ in what below.  Set
\begin{align}
M=&\sup \{s:\nu(x:H(x)   \leq s)\leq 1/2\},\nonumber\\
D_M=& \{x:H(x)  \leq M\} .\nonumber
\end{align} By (\ref{popp}) and
Chebyshev inequality,
\begin{align}
 \nu(x:|H(x)  |  \geq 2 C\sqrt{n})\leq 1/4 ,\nonumber
\end{align}
which gives $M\leq 2C\sqrt{n}$. Therefore, when $n>2^4a^{-2}C^2$ we
have
\begin{align}
 \nu(x:H(x)  \geq an )\leq & \nu(x:H(x)-M \geq an-2C\sqrt{n}
 )\nonumber\\
 \leq& \nu(x:H(x)-M \geq \frac{an}{2};|x|\leq 2\sqrt{n})+\nu(x: |x|> 2\sqrt{n})
 \nonumber\\
 \leq& \nu(x:H(x)-M \geq \frac{an}{2};|x|\leq 2\sqrt{n})+\mu_n(x: |x|>
 2\sqrt{n}).\label{inc3}
 \end{align}
where we use   (\ref{21}) in the last step above.

 By assumption  $|\nabla H(x)|\leq  C|x|$ on $\mathbb{R}^n$, when   $r<\sqrt{n}$ we have
  \begin{align}
 D_M[r]\cap B_n(2\sqrt{n})\subseteq
  \{x:H(x)-M \leq 3C\sqrt{n}r \}\cap B_n(2\sqrt{n}),\nonumber
 \end{align}
where $D_M[r]$ is defined by (\ref{the}).  By (\ref{two}) and
Chebyshev inequality, we have $B_n(\sqrt{n})\cap D_M \neq
\emptyset$,
 which implies that $H(x_0)\leq M$ for some $x_0\in B_n(\sqrt{n})$.
 Therefore, when   $r\geq \sqrt{n}$, by assumption $|\nabla H(x)|\leq  C|x|$ on $\mathbb{R}^n$ we have
  \begin{align}
 B_n(2\sqrt{n})\subseteq
  \{x:H(x)-M \leq 4C\sqrt{n}r \}.\nonumber
 \end{align}
The two relations above  show that  for every $r>0$
  \begin{align}
D_M[r]\cap B_n(2\sqrt{n})\subseteq
  \{x:H(x)-M \leq 4C\sqrt{n}r \}\cap B_n(2\sqrt{n}),\nonumber
 \end{align}
 which  implies that
 \begin{align}
 \{x:H(x)-M \geq \frac{an}{2};|x|\leq 2\sqrt{n}\}
 \subseteq&
D_M[\frac{a}{8C}\sqrt{n}]^c.\nonumber
 \end{align}
Applying (\ref{1n}), $\nu(D_M)=1/2$, the isoperimetric inequality
  (\ref{Iso}) and the relation above, we get
 \begin{align}
 \nu(x:H(x)-M \geq \frac{an}{2};|x|\leq 2\sqrt{n})\leq &\mu_1([
 \frac{a}{8C}\sqrt{n},\infty))\nonumber\\
 \leq & \frac{8C}{a\sqrt{2\pi n}}\exp\{-2^{-7}C^{-2}a^2n\} .\label{inc1}
 \end{align}
Applying  (\ref{>}), (\ref{inc3}) and (\ref{inc1}), we have for
every  $n>N_1\vee2^4a^{-2}C^2$\begin{align}
 \nu(x:H(x)  \geq an )\leq& \frac{8C}{a\sqrt{ 2\pi n}}\exp\{-2^{-7}C^{-2}a^2n\} +e^{-n/2}\nonumber\\
 \leq &  \exp\{- 2^{-7}C^{-2}a^2n\}+e^{-n/2} ,\nonumber
 \end{align}
which gives the conclusion. \qed\medskip

Next we study  some assistant functions with dilation parameter.

\begin{lemma} \label{upper,p}
Suppose that  $u=e^{-U}$ and $v=e^{-V}$ be   differentiable
functions of $\mathcal{CF}_n$ with $|\nabla U|+|\nabla V|$
controlled by some polynomial. Let   $k$ be a positive integer and
$r\in(0,1]$. Define constant $M_r$    and probability measure
$\nu_r$ on $\mathbb{R}^n$ by
   \begin{align} M_r= &   \int  u(rx)  v(rx) d\mu_n, \nonumber  \\
d\nu_r(x)= &     M_r^{-1}u(rx)  v(rx) d\mu_n(x) .
\label{nu}\end{align} Then  \begin{align} &\frac{d}{dr}\int
H_{1,r}(x)u(rx)  v(rx) d\mu_n
\nonumber\\=&M_r^{-1}\big(\frac{dM_r}{dr}\big) \int H_{1,r}(x)u(rx)
v(rx) d\mu_n+ r^{-1}M_r\int
(H_{1,r}-a_{1,r})(H_{2,r}-a_{2,r})d\nu_r\nonumber\\
-&kr^{-(k+1)} M_r\int \frac{ \langle \nabla U(rx), \nabla
 V(rx)\rangle^2}{1+r^{-k}\langle \nabla U(rx), \nabla
 V(rx)\rangle^2} d\nu_r,\nonumber
\end{align} where
\begin{align}
 H_{1,r}(x)=&\ln(1+r^{-k}\langle \nabla
U(rx), \nabla
 V(rx)\rangle^2),\ \ \ \ \ \forall\  x\in \mathbb{R}^n,\nonumber\\
H_{2,r}(x)=& |x|^2-n,\ \ \ \ \ \ \ \ \ \ \ \ \ \ \ \ \ \ \  \ \ \ \
\ \ \ \ \ \ \ \ \ \ \ \forall\ x\in \mathbb{R}^n,\nonumber
\end{align}
and
\begin{align}
 a_{1,r} =& \int H_{1,r}d\nu_r ,\ \ \ \ \ \ \ a_{2,r}= \int  H_{2,r}d\nu_r
 .\nonumber
\end{align}

 \end{lemma}
 \noindent{\bf Proof}\ \ Let $r\in (0,1]$.
We have
\begin{align}
&\frac{d}{dr}\int\ln(1+r^{-k}\langle \nabla U(rx), \nabla
 V(rx)\rangle^2)u(rx)  v(rx)
d\mu_n\nonumber\\=&\frac{d}{dr}\Big(\frac{r^{-n}}{(2\pi)^{n/2}}\int
\ln(1+r^{-k}\langle \nabla U(x), \nabla
 V(x)\rangle^2)u(x)v(x) \exp\{-\frac{|x|^2}{2r^2}\} dx
\Big)\nonumber\\=& -k\frac{r^{-n}}{(2\pi)^{n/2}}\int
\frac{r^{-(k+1)}\langle \nabla U(x), \nabla
 V(x)\rangle^2}{1+r^{-k}\langle \nabla U(x), \nabla
 V(x)\rangle^2}u(x)v(x) \exp\{-\frac{|x|^2}{2r^2}\} dx
\nonumber\\
+& \frac{r^{-(n+1)}}{(2\pi)^{n/2}}\int(r^{-2}|x|^2-n)
\ln(1+r^{-k}\langle \nabla U(x), \nabla
 V(x)\rangle^2)\exp\{-\frac{|x|^2}{2r^2}\}u(x)v(x)
dx\nonumber\\=& -kr^{-(k+1)} \int \frac{ \langle \nabla U(rx),
\nabla
 V(rx)\rangle^2}{1+r^{-k}\langle \nabla U(rx), \nabla
 V(rx)\rangle^2}u(rx)v(rx) d\mu_n
\nonumber\\
+&  r^{-1}\int(|x|^2-n) \ln(1+r^{-k}\langle \nabla U(rx), \nabla
 V(rx)\rangle^2)u(rx)v(rx)
  d\mu_n \nonumber\\
=& r^{-1}M_r\int(|x|^2-n) \ln(1+r^{-k}\langle \nabla U(rx), \nabla
 V(rx)\rangle^2)
  d\nu_r\nonumber\\
  -&kr^{-(k+1)} M_r\int \frac{ \langle \nabla U(rx), \nabla
 V(rx)\rangle^2}{1+r^{-k}\langle \nabla U(rx), \nabla
 V(rx)\rangle^2} d\nu_r.\label{s1}
\end{align}
Similarly,
\begin{align}
&\frac{d}{dr}\int u(rx)  v(rx) d\mu_n = r^{-1}M_r\int(|x|^2-n)
  d\nu_r =r^{-1}M_ra_{2,r}.\label{s3}
\end{align}
 By (\ref{momen}) and (\ref{s1}),
\begin{align} & r^{-1}M_r\int(|x|^2-n) \ln(1+r^{-k}\langle \nabla U(rx), \nabla
 V(rx)\rangle^2)
  d\nu_r \nonumber\\=&r^{-1}M_r  a_{1,r}a_{2,r}+\ r^{-1}M_r\int
(H_{1,r}-a_{1,r})(H_{2,r}-a_{2,r})d\nu_r ,\nonumber
\end{align}
which gives the conclusion with the help of (\ref{s3}). \qed\medskip

\begin{lemma} \label{upper,pp}
Suppose that  $u=e^{-U}$ and $v=e^{-V}$ be twice differentiable
functions of $\mathcal{CF}_n$.  Let $C>0$ and assume that for every
$ x\in \mathbb{R}^n$\begin{align} \nabla^2 U(x) \leq  C I_n,\ \ \ \
\ \nabla^2 V(x) \leq  C I_n.\label{he}
\end{align}  Let $\varepsilon_0,r_0\in (0,1]$ and
 $\varepsilon\in[\varepsilon_0,\infty ),r\in[r_0,1]$. Assume that
\begin{align}
 \int \langle \nabla U(rx), \nabla
V(rx)\rangle u(rx)v(rx) d\mu_n=-\varepsilon  n \int  u(rx)v(rx)
d\mu_n\label{conct}.
\end{align}
 Then for every integer $k\geq 2^8\varepsilon_0^{-1}C^2$,
 there exists  some integer $N_5=N_5(\varepsilon_0,r_0,k,C)$ such that  for every
$n\geq N_5$
 \begin{align} &\frac{d}{dr}\int\ln(1+r^{-k}\langle \nabla
U(rx), \nabla
 V(rx)\rangle^2)u(rx)v(rx)
d\mu_n\nonumber\\\leq &M_r^{-1}\big(\frac{dM_r}{dr}\big) \int\ln
(1+r^{-k}\langle \nabla U(rx), \nabla
 V(rx)\rangle^2)u(rx)  v(rx)
d\mu_n\nonumber.
\end{align}

 \end{lemma}
 \noindent{\bf Proof}\ \ Let $\varepsilon\in[\varepsilon_0,\infty )$, $r\in [r_0,1]$
  and $k\geq 2^8\varepsilon_0^{-1}C^2$. Define constants
 $a_{1,r},a_{2,r}$, functions $H_{1,r},H_{2,r}$ and probability
 measure $\nu_r$ the same as  those     in  Lemma \ref{upper,p}.
   By  the second inequality of (\ref{fi}) and
assumption (\ref{he}), we have  for every $ x\in \mathbb{R}^n$
\begin{align}
  |\nabla U(rx)|\leq Cr|x|,\ \ \  |\nabla V(rx)|\leq Cr|x|.\nonumber
\end{align}
From  the estimate above and     assumption (\ref{he}), we have  for
  every  $x\in \mathbb{R}^n$
\begin{align}
 &|\nabla H_{1,r}(x)|\nonumber\\\leq &2r^{-k}|1+r^{-k}\langle \nabla U(rx), \nabla
V(rx)\rangle^2|^{-1}|\langle \nabla U(rx), \nabla V(rx)\rangle||
r\nabla
V(rx)\cdot  \nabla^2 U(rx) +r  \nabla U(rx)\cdot \nabla^2 V(rx)| \nonumber\\
\leq & 4rCW(x)(|\nabla V(rx)|+
  |\nabla U(rx)|) \nonumber\\ \leq & 8r^2C^2W(x)  |x|, \label{H1}\
\end{align}
where
\begin{align}
W(x)=& r^{-k},\ \ \ \ \ \ \ \ \ \ \ \ \ \ \ \ \ \ \ \ \ \ \ \ \ \ \
\ \ \ \ \ \ \ \ \ \ if\ |\langle \nabla U(rx), \nabla
V(rx)\rangle|\leq 1,\nonumber\\W(x)=&(1+|\langle \nabla U(rx),
\nabla V(rx)\rangle|)^{-1},\ \  \ \ \ \ otherwise.\nonumber
\end{align}

 Similarly, we have  for   every  $x\in \mathbb{R}^n$
\begin{align}
|\nabla  H_{3,r}(x)|\leq 2r^2C^2|x| ,\label{H3}
\end{align}
where
\begin{align}
   H_{3,r}(x)=    \langle\nabla U(rx), \nabla
V(rx)\rangle,\ \ \ \ \ \forall\ x\in \mathbb{R}^n.\nonumber
\end{align}
 We also have for every  $x\in \mathbb{R}^n$
\begin{align}
 |\nabla H_{2,r}(x)|=&2 |x| .\label{H2}
\end{align}

When $n\geq 2\varepsilon^{-1}$, applying    Poincar\'{e} inequality
(\ref{Poinc}),  (\ref{H1}) and (\ref{H2}), we obtain
\begin{align} \int
(H_{1,r}-a_{1,r})^2d\nu_r  \leq&4(2rC)^4 \int W(x)^2|x|^2
d\nu_r\nonumber\\\leq&4(2rC)^4 \int_{\langle \nabla U(rx), \nabla
V(rx)\rangle\leq -\frac{\varepsilon n}{2}}  \frac{|x|^2
}{(1+|\langle \nabla U(rx), \nabla
V(rx)\rangle|)^2}d\nu_r\nonumber\\
+ &4(2rC)^4r^{-2k} \int_{\langle \nabla U(rx), \nabla
V(rx)\rangle>-\frac{\varepsilon n}{2}}  |x|^2  d\nu_r\nonumber\\
\leq&2^8\varepsilon^{-2}r^4C^4n^{-2}\int  |x|^2
d\nu_r+4(2rC)^4r^{-2k} \int_{\langle \nabla U(rx), \nabla
V(rx)\rangle>-\frac{\varepsilon n}{2}}  |x|^2  d\nu_r
 .\label{ii}
\end{align}
By Lemma \ref{Ho19} and (\ref{4n}),\begin{align}   \int |x|^4d\nu_r
\leq
 n^2+2n.\label{pf}
\end{align}
By    (\ref{inc4}),  (\ref{conct}),  (\ref{H3}) and $r\in (0,1]$, we
have for every $n>  N_4(2^{-1}\varepsilon_0, 2 C^2)$
\begin{align} & \nu_r(x:\langle \nabla U(rx), \nabla
V(rx)\rangle>-\frac{\varepsilon n}{2})  \nonumber\\
\leq &\nu_r(x:\langle \nabla U(rx), \nabla V(rx)\rangle- \int
\langle \nabla U(rx), \nabla V(rx)\rangle d\nu_r>\frac{\varepsilon
n}{2}) \nonumber\\ \leq&  \exp\{- 2^{-11}
C^{-4}\varepsilon^2n\}+e^{-n/2}.\label{ii''j}
\end{align}

By Cauchy-Schwartz inequality,   (\ref{pf}) and (\ref{ii''j}), we
have for every   $n> N_4( 2^{-1}\varepsilon_0, 2 C^2)$
\begin{align}  \int_{\langle \nabla U(rx), \nabla
V(rx)\rangle>-\frac{\varepsilon n}{2}} |x|^2d\nu_r
 \leq& \big(
\int |x|^4d\nu_r\big)^{1/2}\nu_r(x:\langle \nabla U(rx), \nabla
V(rx)\rangle>-\frac{\varepsilon n}{2})^{1/2} \nonumber\\\leq&
\sqrt{n^2+2n} (\exp\{- 2^{-12}
C^{-4}\varepsilon^2n\}+e^{-n/4}).\label{ii''}
\end{align}
By (\ref{ii}), the first  inequality of (\ref{pf})  and
(\ref{ii''}),  there exists  some integer
$N_5'=N_5'(\varepsilon_0,r_0,k,C)$ such that for every  $n\geq N_5'$
\begin{align}
\int (H_{1,r}-a_{1,r})^2d\nu_r\leq
2^9\varepsilon^{-2}r^4C^4n^{-1}.\nonumber
\end{align}
By    (\ref{popp}), we also have
 \begin{align} &\int (H_{2,r}-a_{2,r})^2d\nu_r \leq 4n
.\nonumber
\end{align}
Applying Cauchy-Schwartz inequality and     the  two estimates
above,  we have  for every $n\geq N_5'$ \begin{align} \int
(H_{1,r}-a_{1,r})(H_{2,r}-a_{2,r})d\nu_r \leq & \Big(\int
(H_{1,r}-a_{1,r})^2d\nu_r\Big)^{1/2}\int\Big(H_{2,r}-a_{2,r})^2d\nu_r\Big)^{1/2}\nonumber\\
\leq &2^6\varepsilon^{-1}r^2C^2\nonumber\\
\leq &2^6\varepsilon_0^{-1}C^2. \nonumber
\end{align}
By (\ref{ii''j}), there exists some  integer
$N_5''=N_5''(\varepsilon_0,C)\geq 2\varepsilon_0^{-1}$ such that
for every $n\geq N_5''$
 \begin{align} & r^{-k}  \int \frac{ \langle \nabla U(rx), \nabla
 V(rx)\rangle^2}{1+r^{-k}\langle \nabla U(rx), \nabla
 V(rx)\rangle^2} d\nu_r\nonumber\\
\geq & r^{-k}  \int_{ \langle \nabla U(rx), \nabla V(rx)\rangle\leq
-\frac{\varepsilon n}{2}} \frac{ \langle \nabla U(rx), \nabla
 V(rx)\rangle^2}{1+r^{-k}\langle \nabla U(rx), \nabla
 V(rx)\rangle^2} d\nu_r\nonumber\\
\geq & \frac{1}{2}  \int_{ \langle \nabla U(rx), \nabla
V(rx)\rangle\leq -\frac{\varepsilon n}{2}}  d\nu_r\nonumber\\
\geq & 1/4.\nonumber
\end{align}  Then we get the conclusion by the two estimates above and   Lemma
\ref{upper,p}.\qed\medskip

\begin{corollary} \label{upper,ppp}
Suppose that  $u=e^{-U}$ and $v=e^{-V}$ be twice differentiable
functions of $\mathcal{CF}_n$ satisfying condition (\ref{he}) for
some constant $C>0$. Let $\varepsilon_0,r_0\in (0,1]$,
$\varepsilon\in [ \varepsilon_0,\infty),r\in [r_0,1]$ and assume
that
\begin{align}
 \int \langle \nabla U(rx), \nabla
V(rx)\rangle u(rx)v(rx) d\mu_n=-\varepsilon  n \int  u(rx)v(rx)
d\mu_n\label{conct,n}.
\end{align}
 Let $k$ be a positive integer,  and let    $H_{1,r} $, $\nu_r$ be  the function
 and the  measure defined  in  Lemma \ref{upper,p} respectively. Then there
 exists some  integer $N_6=N_6(\varepsilon_0,r_0,k,C)$ such that  for every
$n\geq N_6$
 \begin{align}
  \ln(1+r^{-k}(3^{-1}\varepsilon n)^2)<  \int H_{1,r} d\nu_r <
\ln(1+r^{-k}( 2\varepsilon n)^2).\label{ulb}
\end{align}

 \end{corollary}
 \noindent{\bf Proof}\ \ Let $r\in [r_0,1]$ and  $\varepsilon\in [ \varepsilon_0,\infty)$.
 Applying (\ref{ii''j}) and (\ref{conct,n}), we have  for every
$n\geq N_4(2^{-1}\varepsilon_0 , 2 C^2)$
 \begin{align}
\int  H_{1,r} d\nu_r \geq&\int_{\langle \nabla U(rx), \nabla
V(rx)\rangle\leq -\frac{\varepsilon n}{2}} \ln (1+r^{-k}\langle
\nabla U(rx), \nabla
 V(rx)\rangle^2)d\nu_r\nonumber\\
\geq &(1-\exp\{- 2^{-11}
C^{-4}\varepsilon^2n\}-e^{-n/2})\ln(1+r^{-k}(2^{-1}\varepsilon
n)^2),\nonumber
 \end{align}
which implies the first inequality of (\ref{ulb}). By (\ref{H3}) and
$\nabla U(0)=0$,
\begin{align}
|\langle\nabla U(rx), \nabla V(rx)\rangle|\leq (rC|x|)^2 ,\ \ \ \
\forall  x\in \mathbb{R}^n.\label{H.3}
\end{align}
Applying  (\ref{H3}),  (\ref{pf}), (\ref{ii''j}),  (\ref{conct,n})
and (\ref{H.3}), we get for every $ n\geq N_4(2^{-1}\varepsilon_0 ,
2 C^2)$
 \begin{align}
&\int  H_{1,r}  d\nu_r\nonumber\\ =&\int_{|\langle \nabla U(rx),
\nabla V(rx)\rangle+\frac{\varepsilon n}{2}|\leq \frac{\varepsilon
n}{2}} H_{1,r} (x)d\nu_r+ \int_{|\langle \nabla U(rx), \nabla
V(rx)\rangle+\frac{\varepsilon n}{2}|> \frac{\varepsilon n}{2}}
H_{1,r} (x)d\nu_r\nonumber\\\leq  &\ln(1+r^{-k}(
\frac{3}{2}\varepsilon n)^2)+2\int_{|\langle \nabla U(rx), \nabla
V(rx)\rangle+\frac{\varepsilon n}{2}|> \frac{\varepsilon
n}{2}}r^{-k/2}(rC|x|)^2 d\nu_r \nonumber\\\leq &\ln(1+r^{-k}(
\frac{3}{2}\varepsilon n)^2)+2r^{-k/2}(rC)^2\nu_r(x:{|\langle \nabla
U(rx), \nabla V(rx)\rangle+\frac{\varepsilon n}{2}|>
\frac{\varepsilon n}{2}})^{1/2}(\int
|x|^4d\nu_r)^{1/2}\nonumber\\\leq &\ln(1+r^{-k}(
\frac{3}{2}\varepsilon n)^2)+2^2r^{-k/2}(rC)^2(n^2+2n)^{1/2}\exp\{-
2^{-12} C^{-4}\varepsilon^2n\},\nonumber
 \end{align}
which implies  the second inequality of (\ref{ulb}). \qed\medskip

\begin{lemma} \label{assis}
Let
 $\beta \in \mathbb{R}$ and $r_0\in (0,1)$. Let   $f$  and $g$ be    differentiable
functions on $[r_0,1]$ with $f(1)=\beta g(1) $. Suppose that $g(r)
>0$ for every  $r\in [r_0,1]$. Suppose also  that   for every   $
r\in [r_0,1]$
 \begin{align}
 \frac{d}{dr} f(r) \leq& f(r) g(r)^{-1} \frac{d}{dr}g(r)  .\label{fd}
\end{align}
Then
 \begin{align}
  f(r) \geq&   \beta g(r) ,\ \ \ \forall\ r\in[r_0,1].\nonumber
\end{align}

 \end{lemma}
 \noindent{\bf Proof}\
Set  $h(r)=f(r)/g(r)$ for every $r\in [r_0,1]$. From assumption
$f(1)=\beta g(1)  $,
 we have $h(1)=\beta $.
 Applying (\ref{fd}) and the assumption that  $g>0$ on $[r_0,1]$, we have for every  $r\in [r_0,1]$\begin{align}
 \frac{d}{dr} h(r) =&   g(r)^{-1} \frac{d}{dr}f(r)-f(r)g(r)^{-2}\frac{d}{dr} g(r)
  \leq 0. \nonumber
\end{align}
The estimate above,  the assumption that $g>0$ on $[r_0,1]$ and
$h(1)=\beta $ show that $h(r) \geq  \beta$ for every  $r\in
[r_0,1]$. With  the assumption that $g>0$ on $[r_0,1]$, this implies
the conclusion.\qed\medskip

 \begin{lemma} \label{dir}
Let  $\alpha\in (0,1)$ and let $A\in \mathcal{C}_n$ such that
$B_n(\delta \sqrt{n})\subseteq A$ for some  $\delta\in (0,1)$.  Then
there exists some universal integer  $N_7 $ such that for every
$t\in (0,2^{-4}\delta^2)$, $ x\in B_n(\frac{1}{2}\delta \sqrt{n})$
and every $n\geq N_7$
\begin{align}
  |\nabla H_{A,\alpha,t}(x)-\frac{\alpha e^{-t}x}{1+\alpha(1-e^{-t})}|\leq
   e^{-n/6}
  . \end{align}
 \end{lemma}
  \noindent{\bf Proof}\ Recall that function
$H_{A,\alpha,t}$ is defined by (\ref{HA''}).  Set
  $h
(z)=\exp\{- \alpha |z|^2/2\}$ for $z\in \mathbb{R}^n$.   Let $t\in
(0,2^{-4}\delta^2)$ and $x\in \mathbb{R}^n$ such that $
|x|<\frac{1}{2}\delta \sqrt{n}$.  By  $ |x|<\frac{1}{2}\delta
\sqrt{n}$ and $0\leq t\leq 2^{-4}$, we have
$e^{-t/2}x+(1-e^{-t})^{1/2}y\in B_n(\delta\sqrt{n})$ for every $y\in
B_n(2\sqrt{n})$.  This and the assumption $B_n(\delta
\sqrt{n})\subseteq A$ show that for every $y\in B_n( 2 \sqrt{n})$
\begin{align}h_{A,\alpha}
(e^{-t/2}x+(1-e^{-t})^{1/2}y)=h(e^{-t/2}x+(1-e^{-t})^{1/2}y)
.\label{2y}\end{align} With  definition (\ref{hA}), we have
\begin{align}|\rho_A(y)-\rho_A(z)|\leq |y-z|,\ \  for\ y,z\in \mathbb{R}^n;\ \ \ |\nabla h_{A,\alpha}(z)|\leq
(\alpha|z|+ n)h_{A,\alpha}(z),\  \ for\ z\in \mathbb{R}^n\
a.e..\label{alm}\end{align}  Then, by (\ref{OU}) and (\ref{2y}),
  \begin{align}
 \nabla H_{A,\alpha,t}   (x)  =&-
 \frac{e^{-t/2}}{h_{A,\alpha,t} (x)} \int
 \nabla
 h_{A,\alpha} (e^{-t/2}x+(1-e^{-t})^{1/2}y) d\mu_n(y)\nonumber\\
 =&-
 \frac{e^{-t/2}}{h_{A,\alpha,t} (x)}  \int_{|y|>
 2\sqrt{n}}
 \nabla
 h_{A,\alpha} (e^{-t/2}x+(1-e^{-t})^{1/2}y) d\mu_n(y)\nonumber\\
 +&\frac{e^{-t/2}}{h_{A,\alpha,t} (x)}
   \int_{|y|> 2\sqrt{n}}
 \nabla
 h (e^{-t/2}x+(1-e^{-t})^{1/2}y) d\mu_n(y)\nonumber\\-&
 \frac{e^{-t/2}}{h_{A,\alpha,t} (x)}
\int
 \nabla
 h (e^{-t/2}x+(1-e^{-t})^{1/2}y) d\mu_n(y)\nonumber\\
: =&-I_1(x)+I_2(x)-I_3(x).\label{I12}\end{align}

Noticing  that $ |x|\leq \frac{1}{2} \sqrt{n}$ and $0<h_{A,\alpha}
(z),h(z)\leq e^{-\alpha|z|^2/2}$ for all  $z\in \mathbb{R}^n$,  we
have by  (\ref{alm}) and  (\ref{I12})
\begin{align} &|I_1(x)-I_2(x)|\nonumber\\  \leq &
 \frac{2 e^{-t/2}}{h_{A,\alpha,t} (x)} \int_{|y|>
 2\sqrt{n}}
 (\alpha|e^{-t/2}x+(1-e^{-t})^{1/2}y|+ n)\exp\{-\frac{\alpha|e^{-t/2}x+(1-e^{-t})^{1/2}y|^2}{2}\}d\mu_n(y)
\nonumber\\ \leq &
 \frac{4n}{h_{A,\alpha,t} (x)}    \int_{|y|>
 2\sqrt{n}}
  d\mu_n(y) , \label{I2}\end{align}
  where we use assumption  $\alpha\in (0,1)$ and the inequality $c\cdot e^{-c^2/8}\leq 2$ for $c>0$ in
  the last step above.

  Direct calculation shows that
 \begin{align}
P_th(x)= & (1+\alpha(1-e^{-t}))^{-n/2}
  \exp\{-\frac{\alpha e^{-t}  |x|^2}{2(1+\alpha(1-e^{-t}))}\}
   \nonumber,\end{align}
which gives
  \begin{align}
I_3(x) =&
   \frac{ 1}{h_{A,\alpha,t} (x)}  \nabla \int  h(e^{-t/2}x+(1-e^{-t})^{1/2}y)   d\mu_n(y)\nonumber\\
 =&-
 \frac{P_th(x)}{h_{A,\alpha,t} (x)}\frac{\alpha e^{-t}x}{1+\alpha(1-e^{-t})}\nonumber\\
=&-\frac{\alpha e^{-t}x}{1+\alpha(1-e^{-t})}-
 \frac{P_th(x)-h_{A,\alpha,t} (x)}{h_{A,\alpha,t} (x)}
\frac{\alpha e^{-t}x}{1+\alpha(1-e^{-t})} . \label{I1,3}\end{align}
By $ |x|\leq \frac{1}{2} \sqrt{n}$ and $\alpha\in (0,1)$, we have
 \begin{align}
|\frac{\alpha e^{-t}x}{1+\alpha(1-e^{-t})}|\leq \frac{1}{2}
\sqrt{n}.
 \end{align}
By (\ref{2y}),
  \begin{align}
  &  \frac{|P_th(x)-h_{A,\alpha,t} (x)|}{h_{A,\alpha,t} (x)}
 \nonumber\\
\leq & \frac{ 1}{h_{A,\alpha,t} (x)}  \int_{|y|>
 2\sqrt{n}}
\big(h_{A,\alpha} (e^{-t/2}x+(1-e^{-t})^{1/2}y)+h
(e^{-t/2}x+(1-e^{-t})^{1/2}y)\big)
   d\mu_n(y)\nonumber\\
\leq &   \frac{ 2}{h_{A,\alpha,t} (x)}    \int_{|y|>
 2\sqrt{n}}
   d\mu_n(y).
\label{has}\end{align}

Combing  (\ref{I12}) and  (\ref{I1,3})-(\ref{has}), we have
\begin{align}
  |\nabla H_{A,\alpha,t}(x)-\frac{\alpha e^{-t}x}{1+\alpha(1-e^{-t})}|
  \leq \frac{\sqrt{n}}{h_{A,\alpha,t} (x)}    \int_{|y|>
 2\sqrt{n}}
   d\mu_n(y)+|I_1(x)-I_2(x)|. \label{ch}\end{align}
By (\ref{2y}),   $ |x|\leq \frac{1}{2} \sqrt{n}$ and applying
(\ref{1/2}), (\ref{>}) and  (\ref{2y}), we get for $n\geq N_1$
\begin{align}
 &\frac{1}{h_{A,\alpha,t} (x)}   \int_{|y|>
 2\sqrt{n}} d\mu_n(y)\nonumber\\
 \leq  & 2^ne^{-3n/2} \big(\int h_{A,\alpha}(e^{-t/2}x+(1-e^{-t})^{1/2}y) d\mu_n(y)\big)^{-1}\nonumber\\
 \leq  & 2^ne^{-3n/2}\big(\int_{|y|<
 \sqrt{n}}
 \exp\{-\frac{\alpha}{2}|e^{-t/2}x+(1-e^{-t})^{1/2}y|^2\} d\mu_n(y)\big)^{-1}\nonumber\\
  \leq  & 2^ne^{-3n/2}\big(\int_{|y|<
 \sqrt{n}}
 \exp\{-\frac{\alpha}{2}( |x|^2+ |y|^2)\} d\mu_n(y)\big)^{-1}\nonumber\\
 \leq  &2^ne^{-3n/2} e^{5  n/8}   \big(\int_{|y|<
 \sqrt{n}} d\mu_n(y)\big)^{-1}\nonumber\\
 \leq  &2^{n-1} e^{-7  n/8}     .\nonumber
  \end{align}
Applying    (\ref{I2}),  (\ref{ch})   and  the estimate above, we
get the conclusion.
  \qed\medskip

  \begin{lemma} \label{pos}
Let  $u=e^{-U} \in\mathcal{CF}_n$ satisfying condition
$\mathcal{L}(C_1,C_2)$ for some constants    $0<C_1<C_2$.  Let
$\delta,\alpha\in (0,1)$ and $A\in \mathcal{C}_n$ such that
$B_n(\delta \sqrt{n})\subseteq A$.  Then there exist some  constant
$C_7=C_7(\alpha, C_1,C_2)>0$ and integer $N_8=N_8(\alpha, C_1,C_2)$
such that for every $n\geq N_8$
\begin{align} \int \langle \nabla h_{A,\alpha,t}(rx), \nabla
u(rx)\rangle d\mu_n
 >C_7r^2   n \int   h_{A,\alpha,t}(rx)
u(rx)  d\mu_n,\ \ \ \  \end{align}
   provided that  $
 e^{-n/6}<r< \delta/4 $ and $t\in
(0,2^{-4}\delta^2)$.
 \end{lemma}
  \noindent{\bf Proof}\ Let $n^{-1/6}<r<   \delta/4
$ and $t\in (0,2^{-4}\delta^2)$.  We have
\begin{align}
\int  \langle \nabla h_{A,\alpha,t}(rx), \nabla u(rx)\rangle \
d\mu_n  =&\Big(\int_{|x|\leq 2\sqrt{n}} +\int_{|x|>
2\sqrt{n}}\Big)\langle \nabla h_{A,\alpha,t}(rx), \nabla
u(rx)\rangle  d\mu_n  \nonumber\\:=&J_1+J_2 .\label{J12}
\end{align}

By the first inequality of (\ref{fi}) and (\ref{trx})
\begin{align}
&\int_{|x|\leq  2\sqrt{n}} \langle x, \nabla U(rx)\rangle
h_{A,\alpha,t}(rx) u(rx)  d\mu_n \nonumber\\ \geq
&\frac{C_1}{C_2}\int_{|x|\leq 2\sqrt{n}}  |\nabla U(rx)||x|
h_{A,\alpha,t}(rx) u(rx)  d\mu_n \nonumber\\\geq
&\frac{rC_1^2}{C_2}\int_{|x|\leq  2\sqrt{n}} |x|^2
h_{A,\alpha,t}(rx) u(rx)  d\mu_n .\label{UV1}
\end{align}

For every  $\widetilde{x}\in S_{n-1}$ and every  $s\in
(0,2\sqrt{n})$, set
$H(s;\widetilde{x})=H_{A,\alpha,t}(sr\widetilde{x})+U(sr\widetilde{x})+\frac{s^2}{2}-(n-1)\ln
s$. Denote  $s_0=(2+C_2)^{-1}\sqrt{n}$. By Lemma \ref{dir} and the
assumptions of $r,t,\alpha$ and $\delta$,
 we have for every  $n\geq 2\vee N_7$  and every  $s\in (0,s_0)$
\begin{align}
 \frac{\partial}{ \partial s}H(s;\widetilde{x})\leq&
\frac{\alpha e^{-t}sr^2}{1+\alpha(1-e^{-t})}  +e^{-n/6}r+C_2sr^2+s-\frac{n-1}{s}\nonumber\\
\leq & sr +e^{-n/6}r+C_2sr+s-  (1+\frac{C_2}{2})\sqrt{n} \nonumber\\
\leq &0.\nonumber
\end{align}
 Then   for every  $n\geq 2\vee N_7$,
\begin{align}
 & \int_{|x|\leq  2\sqrt{n}} |x|^2 h_{A,\alpha,t}(rx) u(rx)
 d\mu_n  \nonumber\\ \geq & (\frac{s_0 }{2})^2
\int_{S_{n-1}}dm_{n-1}(\widetilde{x})\int_{\frac{s_0}{2}}^{2\sqrt{n}}
s^{n-1}
\exp\{-\frac{s^2}{2}-H_{A,\alpha,t}(sr\widetilde{x})-U(sr\widetilde{x})
\} ds\nonumber\\
 \geq &
\frac{s_0^2}{2^3}\int_{S_{n-1}}dm_{n-1}(\widetilde{x})\int_{0}^{2\sqrt{n}}
s^{n-1}
\exp\{-\frac{s^2}{2}-H_{A,\alpha,t}(sr\widetilde{x})-U(sr\widetilde{x})
\} ds\nonumber\\
=&\frac{n}{2^3(2+C_2)^2}\int_{|x|\leq  2\sqrt{n}} h_{A,\alpha,t}(rx)
u(rx)  d\mu_n .\label{UV2}
\end{align}

Set $a_t= \alpha e^{-t}(1+\alpha(1-e^{-t}))^{-1} $. Since   $
      0<r< \delta/4$ by assumption,  we have $|rx|\leq \frac{\delta\sqrt{n}}{2} $ if $|x|\leq 2\sqrt{n}$.
      Then, by  the second inequality of (\ref{fi}),  Lemma \ref{dir}, (\ref{J12})   and (\ref{UV1}) we have for every  $n\geq   N_7$
\begin{align}
J_1 =&  -\int_{|x|\leq  2\sqrt{n}} \langle ra_t x, \nabla
u(rx)\rangle h_{A,\alpha,t}(rx)  d\mu_n\nonumber\\
-&\int_{|x|\leq 2\sqrt{n}} \langle\nabla H_{A,\alpha,t}(rx)- ra_t x,
\nabla u(rx)\rangle
h_{A,\alpha,t}(rx)  d\mu_n \nonumber\\
 \geq &ra_t\int_{|x|\leq  2\sqrt{n}} \langle x, \nabla
U(rx)\rangle h_{A,\alpha,t}(rx) u(rx)  d\mu_n-e^{-n/6}\int_{|x|\leq
2\sqrt{n}}|\nabla U(rx)|
  h_{A,\alpha,t}(rx) u(rx)  d\mu_n \nonumber\\
 \geq &\frac{r^2a_tC_1^2}{C_2}\int_{|x|\leq  2\sqrt{n}}|x|^2  h_{A,\alpha,t}(rx) u(rx)  \
d\mu_n-2re^{-n/6}C_2\sqrt{n}\int_{|x|\leq 2\sqrt{n}}
  h_{A,\alpha,t}(rx) u(rx) d\mu_n .\nonumber
\end{align}
Noticing  that $a_t\geq \alpha/2$ by assumption of $t$ and $\alpha$,
   for every $n\geq 2\vee N_7$ we get by   applying (\ref{UV1}),
(\ref{UV2}) and  the estimate above
\begin{align}
J_1 \geq &\big(\frac{\alpha r^2  C_1^2n}{2^4C_2(2+C_2)^2}
-2re^{-n/6}C_2\sqrt{n}\big)\int_{|x|\leq 2\sqrt{n}}
  h_{A,\alpha,t}(rx) u(rx) d\mu_n .\label{UV4}
\end{align}

Let $x\in \mathbb{R}^n$. Applying
  (\ref{alm}), we have
 \begin{align}
 |\nabla H_{A,\alpha,t}   (x)|  =&\frac{e^{-t/2}}{h_{A,\alpha,t} (x)}\big|
  \int
 \nabla
 h_{A,\alpha} (e^{-t/2}x+(1-e^{-t})^{1/2}y) d\mu_n(y)\big|\nonumber\\
 \leq &\frac{  e^{-t/2}}{h_{A,\alpha,t} (x)}
  \int
(|e^{-t/2}x+(1-e^{-t})^{1/2}y|+n)
 h_{A,\alpha} (e^{-t/2}x+(1-e^{-t})^{1/2}y) d\mu_n(y) \nonumber\\
 \leq &n+
 \frac{1}{h_{A,\alpha,t} (x)(2\pi(1-e^{-t}))^{n/2}}\int
 | y|
 h_{A,\alpha} ( y)  \exp\{-\frac{|y-e^{-t/2}x|^2}{2(1-e^{-t})}\} dy.
 \label{Hd1}\end{align}
Noticing that  $ h_{A,\alpha}( r\widetilde{x})$ is a  decreasing
function of $r\in [0,\infty)$, we have by $t<2^{-4}$
 \begin{align}L_1(x):=& \int_{|y|>4(|x|+\sqrt{n})}
 | y|
 h_{A,\alpha} ( y)  \exp\{-\frac{|y-e^{-t/2}x|^2}{2(1-e^{-t})}\} dy
 \nonumber \\
 =&\int_{S_{n-1}} dm_{n-1}(\widetilde{y})\int_{4(|x|+\sqrt{n})}^\infty
r^n
 h_{A,\alpha} ( r \widetilde{y})  \exp\{-\frac{|r\widetilde{y}-e^{-t/2}x|^2}{2(1-e^{-t})}\}dr  \nonumber \\
\leq &  \int_{S_{n-1}}
dm_{n-1}(\widetilde{y})\int_{4(|x|+\sqrt{n})}^\infty h_{A,\alpha} (
(|x|+\sqrt{n}) \widetilde{y}) r^n
\exp\{-\frac{r^2}{4(1-e^{-t})}\}dr \nonumber \\
\leq &  2(4(|x|+\sqrt{n}))^{n-1}
\exp\{-\frac{4(|x|+\sqrt{n})^2}{1-e^{-t}}\}\int_{S_{n-1}}
h_{A,\alpha} ( (|x|+\sqrt{n}) \widetilde{y})
dm_{n-1}(\widetilde{y}). \nonumber
\end{align}
We also have for $n\geq 4$
\begin{align}L_2(x):=& \int_{|y|<|x|+\sqrt{n}}
 | y|
 h_{A,\alpha} ( y)  \exp\{-\frac{|y-e^{-t/2}x|^2}{2(1-e^{-t})}\} dy
 \nonumber \\
 \geq &\int_{S_{n-1}}h_{A,\alpha} ( (|x|+\sqrt{n}) \widetilde{y}) dm_{n-1}(\widetilde{y})\int^{|x|+\sqrt{n}}_0
r^n \exp\{-\frac{|r\widetilde{y}-e^{-t/2}x|^2}{2(1-e^{-t})}\}dr
\nonumber \\
 \geq &\int_{S_{n-1}} h_{A,\alpha} ( (|x|+\sqrt{n}) \widetilde{y}) dm_{n-1}(\widetilde{y})\int^{|x|+\sqrt{n}}_0
r^n\exp\{-\frac{r^2 +|x|^2}{2(1-e^{-t})}\}dr
 \nonumber \\
\geq & (|x|+\sqrt{n}-1)^{n}
\exp\{-\frac{(|x|+\sqrt{n})^2+|x|^2}{2(1-e^{-t})}\}\int_{S_{n-1}}
h_{A,\alpha} ( (|x|+\sqrt{n}) \widetilde{y})
dm_{n-1}(\widetilde{y})     \nonumber\\
\geq &  2^{-n}(|x|+\sqrt{n})^{n} \exp\{-\frac{
(|x|+\sqrt{n})^2}{1-e^{-t}}\}\int_{S_{n-1}} h_{A,\alpha} (
(|x|+\sqrt{n}) \widetilde{y}) dm_{n-1}(\widetilde{y})  . \nonumber
\end{align}

 Applying the  two estimates above and $t<2^{-4}$ we
have for $n\geq 4$
 \begin{align}L_1(x)/L_2(x)
\leq   2^{3n} \exp\{-\frac{3(|x|+\sqrt{n})^2}{1-e^{-t}}\}   \leq
1.\nonumber
\end{align}
This gives  for $n\geq 4$
 \begin{align}
 \int
 | y|
 h_{A,\alpha} ( y)  \exp\{-\frac{|y-e^{-t/2}x|^2}{2(1-e^{-t})}\} dy
  \leq &2\int_{|y|\leq 4(|x|+\sqrt{n})}
 | y|
 h_{A,\alpha} ( y)  \exp\{-\frac{|y-e^{-t/2}x|^2}{2(1-e^{-t})}\} dy \nonumber\\
 \leq &8(|x|+\sqrt{n})\int
  h_{A,\alpha} ( y)  \exp\{-\frac{|y-e^{-t/2}x|^2}{2(1-e^{-t})}\} dy.\nonumber\end{align}
Applying (\ref{Hd1}) and the estimate above we have for $n\geq 4$
\begin{align}
 |\nabla H_{A,\alpha,t}   (x)|
 \leq &n+8(|x|+\sqrt{n}).\nonumber\end{align}
 Applying
   Lemma \ref{ball,es}, the second inequality of (\ref{fi}) and the estimate above, we have
 for every  $n\geq N_3\vee 4$
\begin{align}
|J_2|\leq &  \int_{|x|> 2\sqrt{n}} |\nabla
H_{A,\alpha,t}(rx)|| \nabla U(rx)|h_{A,\alpha,t}(rx)u(rx) d\mu_n \nonumber\\
 \leq &rC_2 \int_{|x|> 2\sqrt{n}} (n+8(|x|+\sqrt{n})|x|)h_{A,\alpha,t}(rx) u(rx)  d\mu_n\nonumber\\
 \leq &10rC_2   \int_{|x|> 2\sqrt{n}}|x|^3  h_{A,\alpha,t}
 (rx) u(rx)  d\mu_n
 \nonumber\\
 \leq &10 re^{-n/2}C_2
\int_{|x|<\sqrt{n}}h_{A,\alpha,t}(rx)u(rx) d\mu_n.\label{UV3}
\end{align}
By Lemma \ref{Ho19} and (\ref{1/2}),
\begin{align}
 \int_{|x|\leq 2\sqrt{n}}h_{A,\alpha,t}(rx)u(rx) d\mu_n
 \geq  \frac{1}{2}\int h_{A,\alpha,t}(rx)u(rx)
 d\mu_n.\label{UV5}
\end{align}
 Applying (\ref{J12}), (\ref{UV4}) and
(\ref{UV3}), we get    for every  $n\geq N_3\vee 4$
\begin{align}
&\int  \langle \nabla h_{A,\alpha,t}(rx), \nabla u(rx)\rangle d\mu_n
\nonumber\\ \geq &r\big(\frac{\alpha r C_1^2n}{2^4C_2(2+C_2)^2}
-2e^{-n/6}C_2\sqrt{n}-10 e^{-n/2}C_2 \big)\int_{|x|\leq 2\sqrt{n}}
  h_{A,\alpha,t}(rx) u(rx) d\mu_n  . \nonumber
\end{align}
Since  $r>e^{-n/6}$, we get the conclusion by (\ref{UV5}) and the
estimate above.\qed\medskip

\begin{lemma} \label{upper}
 Let $\delta,\alpha\in (0,1)$,  $\varepsilon \in (0,\delta^2/8)$ and  $A,B\in \mathcal{C}_n$
 with
    $B_n(\delta
\sqrt{n})\subseteq A$. Assume that for some constant $t_0>0$
 \begin{align}
\frac{d}{dt}\phi_t(h_{A,\alpha},h_{B,\alpha}) <0,\ \ \ \ \forall
t\in (t_0,\infty).\label{dt}
\end{align}Then  there
exists some integer $N_9=N_9(\alpha,\delta,\varepsilon)$ such that
for every $n\geq N_9$
 \begin{align}
\frac{d}{dt}\phi_t(h_{A,\alpha},h_{B,\alpha}) < & \varepsilon n
 \phi_t(h_{A,\alpha},h_{B,\alpha}),\ \ \ \ \ \forall t\geq \varepsilon .\label{gap}
\end{align}

 \end{lemma}
 \noindent{\bf Proof}\ Let  $C(\alpha)=\min(e^{-3}\alpha,2^{-6}e^{-3})$.   By
 Lemma \ref{low} and assumption (\ref{dt}),  we have
  \begin{align}
\frac{d}{dt}\phi_t(h_{A,\alpha},h_{B,\alpha}) <0,\ \ \ \ \forall
t\in (  4(2\ln 2-\ln C(\alpha)),\infty).\label{dtt}
\end{align}Suppose that (\ref{gap}) does not hold. Then, by
(\ref{dtt}) and  the intermediate value theorem,   there exists some
$t'\in (\varepsilon, 4(2\ln 2-\ln C(\alpha)))$
 such that
 \begin{align}
\big(\frac{d}{dt}\phi_t(h_{A,\alpha},h_{B,\alpha})\big)_{t=t'}
=\varepsilon n
 \phi_{t'}(h_{A,\alpha},h_{B,\alpha}).\nonumber
\end{align}
Applying (\ref{d}) and  (\ref{t/2}),    the equality above shows
that
\begin{align} \int \langle \nabla
h_{A,\alpha,\frac{\varepsilon}{2}}, \nabla
h_{B,\alpha,t'-\frac{\varepsilon}{2}} \rangle  d\mu_n =-2\varepsilon
n
 \phi_{t'}(h_{A,\alpha},h_{B,\alpha}).\label{0}
\end{align}
To finish the proof of the lemma, in what below we assume that
(\ref{0}) holds and
 show  a contradiction when $n$ is big enough depending on $\alpha,\delta$
 and $\varepsilon$.

Let   $r_0=\delta/8$.  First we show that for $n$ big enough
depending on $\varepsilon$ and $\delta$
  \begin{align}
\int \langle \nabla U(rx), \nabla V(rx)\rangle u(rx)v(rx) d\mu_n
\leq -  5^{-1}\varepsilon n \int  u(rx)v(rx) d\mu_n,\ \ \ \forall
r\in [r_0,1].\label{contr}
\end{align}
where $U,u$ and $V,v$ are defined by  \begin{align}
u(x)=e^{-U(x)}=h_{A,\alpha,\frac{\varepsilon}{2}}(x),\ \ \ \ \ \
v(x)=e^{-V(x)}=h_{B,\alpha,t'-\frac{\varepsilon}{2}}(x),\ \ \ \
\forall\ x\in \mathbb{R}^n.\nonumber
\end{align}

For  every $r\in (0,1]$, define probability measure $\nu_r$ by
(\ref{nu}). Set for every $k\geq 1$ and every $r\in (0,1]$
\begin{align} f_k(r)=\int\ln(1+r^{-k}\langle \nabla U(rx), \nabla
 V(rx)\rangle^2)u(rx)v(rx)
d\mu_n ,\ \ \ \
 g(r)= &\int
u(rx)  v(rx)  d\mu_n,\nonumber
\end{align}
By $\varepsilon< t'< 4(2\ln 2-\ln C(\alpha))$,  Lemma
\ref{derivatives,2} and Lemma \ref{C1,2,3;}, we have for every $x\in
\mathbb{R}^n$
\begin{align} C(\alpha)   e^{-T(\alpha)} I_n\leq \nabla^2U(x)\leq
4\varepsilon^{-1}I_n,\ \ \ \  C(\alpha)  e^{-T(\alpha)} I_n\leq
\nabla^2V(x)\leq 4\varepsilon^{-1}I_n  , \label{ul}
\end{align}
where $T(\alpha)=4(2\ln 2-\ln C(\alpha))$.  Let  $k_0=\lfloor5\cdot
2^{12}\varepsilon^{-3}\rfloor+1$. By the upper bounds in (\ref{ul})
and Lemma \ref{upper,pp}, there exists some integer $
N^{(1)}=N^{(1)}(\varepsilon,\delta)$ such that for every
$r\in[r_0,1]$ and every $n\geq N^{(1)}$
 \begin{align}
\frac{d}{dr} f_{k_0}(r) \leq& f_{k_0}(r) g(r)^{-1} \frac{d}{dr} g(r)
 \label{f'}
 \end{align}
 provided  that   \begin{align}
 \int \langle \nabla U(rx), \nabla V(rx)\rangle d\nu_r \leq - 5^{-1}\varepsilon n.
 \label{ccontr}
\end{align}

By the first inequality of (\ref{ulb}), (\ref{0})   and the upper
bound in  (\ref{ul}), there exists
  some integer $N^{(2)}=N^{(2)}(\varepsilon,\delta)$ such that for every $n\geq N^{(2)}$
 \begin{align} f_{k_0}(1)= &\int \ln(1+ \langle \nabla U(x), \nabla
 V(x)\rangle^2)d\nu_1 \int u(x)v(x)
d\mu_n  \nonumber\\
 \geq &  (\ln(1+ (2\cdot 3^{-1}\varepsilon n)^{2})) g(1). \label{f(1)}
\end{align}
where  $\nu_1$ is defined by (\ref{nu}).   By assumption (\ref{0}),
there exists some $r_1\in [r_0,1)$ such that (\ref{ccontr}) holds
for every $r\in [r_1,1]$. Then, applying Lemma \ref{assis},
(\ref{f'}) and (\ref{f(1)}), we have for every $ r\in[r_1,1]$ and
every $n\geq N^{(1)}\vee N^{(2)}$
 \begin{align}
 f_{k_0}(r) \geq  (\ln(1+r^{-k_0}(2\cdot 3^{-1}\varepsilon n)^{2}))
 g(r).\label{e2}
 \end{align}
This implies that
  for every $n\geq N^{(1)}\vee N^{(2)}\vee N^{(3)}$ with $N^{(3)}=N_6(3^{-1}\varepsilon,r_0,k_0,4\varepsilon^{-1})$
   \begin{align}
 \int \langle \nabla U(rx), \nabla V(rx)\rangle  d\nu_r < -
3^{-1}\varepsilon n ,\ \ \ \forall r\in [r_1,1].\label{e1}
\end{align}
In fact, if (\ref{e1}) does not hold, then by assumption (\ref{0})
and  the intermediate value theorem,  there exists some some $n\geq
N^{(1)}\vee N^{(2)}\vee N^{(3)}$ and some $r'\in [r_1,1]$ such that
(\ref{e1}) is an equality for $r=r'$. This and  the second
inequality of (\ref{ulb}) give
 \begin{align}
 f_{k_0}(r') <  (\ln(1+r'^{-k_0}(2\cdot 3^{-1}\varepsilon n)^{2}))
 g(r').\label{e3}
 \end{align}
which contracts (\ref{e2}).

Let  $n\geq N^{(1)}\vee N^{(2)}\vee N^{(3)}$.  From the conclusion
(\ref{e1}), we see that if (\ref{ccontr}) holds for every $r\in
[r_1,1]$ with  some $r_1\in [r_0,1)$,  then there exists some
$r_2<r_1$ such that (\ref{ccontr}) holds for every $r\in [r_2,1]$.
By standard continuity arguments, this implies that
   \begin{align}
\int \langle \nabla U(rx), \nabla V(rx)\rangle u(rx)v(rx) d\mu_n
\leq -  5^{-1}\varepsilon n \int  u(rx)v(rx) d\mu_n,\ \ \ \forall
r\in [r_0,1], \nonumber
\end{align}
which gives  conclusion of (\ref{contr}).

 By
(\ref{ul}) and Lemma \ref{pos}, there  exists some integer
$N^{(4)}=N^{(4)}(\varepsilon, \alpha)$ such that for every  $n\geq
N^{(4)}$
  \begin{align}
\int \langle \nabla U(rx), \nabla V(rx)\rangle u(rx)v(rx) d\mu_n
> 0,\ \ \ \ \ \ for\ e^{-n/6}<r< \frac{1}{4}\delta.\nonumber
\end{align}
Since $r_0=\delta/8$, this    contradicts   (\ref{contr}) when
$n\geq N^{(1)}\vee N^{(2)}\vee N^{(3)}\vee N^{(4)}\vee (-6\ln
(2^{-3}\delta)) $.
 \qed\medskip

 \section{Proof of   Theorem 1.1 }

\subsection{\normalsize   further derivative estimates  }

To characterize the equality in (\ref{1}), we prepare the following
lemmas.
 \begin{lemma}\label{Deco} Suppose  that $A$  is a
 closed,
  non-degenerate    element of $\mathcal{C}_n$  and   it is  not equal to  $\mathbb{R}^n$.
  Then, either  $A$ is  bounded or
   $A$ is of the form $\widetilde{A}\times\mathbb{R}^{n-k}$    after some   orthogonal
  transformation, where $k\in\{1,\cdots, n-1\}$ and
  $\widetilde{A}\in \mathcal{C}_k$ is bounded.
\end{lemma}
\noindent{\bf Proof}\ To prove the lemma we can assume  that $A$ is
 unbounded in what below.  Since $A$ is a unbounded, symmetric   and      convex, there
exists     $( \overrightarrow{\theta}_l)_{l\geq 1}\in S_{n-1}$ such
that
  $\{y:y=t \overrightarrow{\theta}_l,\ t\in
[0,l]\}\subseteq A$ for every  $l\geq 1$. Therefore,  from the
assumption that $A$ is closed and symmetric we have $\{y:y=t
\overrightarrow{\theta},\ t\in \mathbb{R}\}\subseteq A$ for  some
element $ \overrightarrow{\theta}\in S_{n-1}$.

Let  $k$ be  the smallest integer of $\{1,\cdots,n-1\}$ such that
 $  K\subseteq A$ for some  $(n-k)$-dimensional subspace $K$ of $\mathbb{R}^n$.
  By orthogonal transformation, we   assume in what below that
$K=\mathbf{0}_k\times \mathbb{R}^{n-k}  $, where $\mathbf{0}_k$ is
the zero point of $\mathbb{R}^k$. To prove the lemma,   it is
sufficient to verify that
\begin{align}( \tilde{x}_k,z_{k+1},\cdots,z_n)\in  A,\ \ \ \ \
\forall\
(z_{k+1},\cdots,z_{n})\in\mathbb{R}^{n-k}\label{ve}\end{align}
provided that $  (\tilde{x}_k,x_{k+1},\cdots,x_n)\in A$ for some
$\tilde{x}_k\in \mathbb{R}^k$ and some $(x_{k+1},\cdots,x_n)\in
\mathbb{R}^{n-k}$. Let
 $x=(
\tilde{x}_k,x_{k+1},\cdots,x_n)\in A$. For every $r\in (0,1)$ and
every  $(z_{k+1},\cdots,z_{n})\in\mathbb{R}^{n-k}$, we have
 by the assumption  of $A$ and $K\subseteq A$
 \begin{align}
&( r\widetilde{x}_k,z_{k+1},\cdots,z_{n})\nonumber\\= & (1-r) \big(
\mathbf{0}_k,(1-r)^{-1}
(z_{k+1},\cdots,z_{n})-(1-r)^{-1}r(x_{k+1},\cdots,x_{n})\big)+ r
x\in A.\nonumber
 \end{align}
Then we get  (\ref{ve}) by letting  $r\rightarrow 1$ in the above
formula.
  \qed\medskip

 \begin{lemma}\label{po}
 Suppose that $A=\widetilde{A}\times\mathbb{R}^{n-k}$ for some
 bounded set
 $\widetilde{A}\in \mathcal{A}_k$ with  $1\leq k\leq n $.
 Let $\mathbf{a}=(a_1,\cdots,a_n)\in \mathbb{R}^n$ such that
  $|(a_1,\cdots,a_k)|>0$.
Let   $B=\{x:|\langle x,\mathbf{a}\rangle|\leq a\}$ for some
 $a>0$.
 Suppose that
$A$ is     non-degenerate.
  Then
 \begin{align}   \mu_{n}(A\cap B)>\mu_n(A)\mu_n(B).\nonumber\end{align}\end{lemma}
 \noindent{\bf Proof}\
From the assumption that $(a_1,\cdots,a_k)\neq 0$ and the previous
proof of $\mu_{n}(A\cap B)\geq \mu_n(A)\mu_n(B)$ in, e.g.,
\cite{Sid67} and \cite{JOG},  to verify the conclusion  we only need
to show that
\begin{align}\mu_{k}(\widetilde{A}+r(a_1,\cdots,a_k)), \ \ \ \ for\ r\geq
0,\nonumber
\end{align}
is a strictly decreasing function of $r\geq 0$. Since
$\widetilde{A}$ is a non-degenerate bounded set in $\mathbb{R}^k$
from the assumption, by Corollary 2 in \cite{An55}, we know that
this property holds.
\medskip\qed

 \begin{lemma}\label{De} Let   $A,B\in {\mathcal{C}}_n$.
 Suppose that $\overline{A}$ and $\overline{B}$ are  not  unlinked   and
 both of them are  non-degenerate.
  Then
 \begin{align}\label{greater} \sum_{i=1}^n
 \int_A (1-|x_i|^2)d\mu_n(x)\int_B (1-|y_i|^2) d\mu_{n}(
y)>0.\end{align}\end{lemma} \noindent{\bf Proof}\  Without loss of
generality, we assume that $A$ and $B$ are both closed sets in what
below. Since $A$ and $B$ are not unlinked, we have by definition
that neither
 $A$ nor  $B$ is equal to  $ \mathbb{R}^n$.
 By Lemma \ref{Deco}, there exist bounded   sets $\widetilde{A}\in \mathcal{C}_{k_1},\widetilde{B}\in\mathcal{C}_{k_2}$
 for some   integers $1\leq
k_1,k_2\leq n-1$ and  orthogonal transformations $\Psi_1$ and
$\Psi_2$ of $\mathbb{R}^n$ such  that
  $\Psi_1(A)=\widetilde{A}\times\mathbb{R}^{n-k_1}$
  and  $\Psi_2(B)=\widetilde{B}\times\mathbb{R}^{n-k_2}$. Let $r_0$ be a constant  such that
\begin{align}\label{bound}
|\widetilde{x}|<r_0,\ \  \forall\   \widetilde{x}\in \widetilde{A} .
\end{align}

For $1\leq k\leq n$, denote by  $\mathbf{0}_k$ the zero point of
$\mathbb{R}^k$.  For $O_1,O_2\subseteq\mathbb{R}^n$, denote
$O_1\perp O_2$ if $\langle x,y\rangle=0$ for every  $x\in O_1$ and
every $y\in O_2$; denote  $O_1 \oplus O_2=\{x+y: x\in O_1,y\in O_2
\}$ if $O_1\perp O_2$.
 Since  $\Psi_1,\Psi_2$ are orthogonal transformations,
 we have
\begin{align}
\Psi_1(\Psi_2^{-1}(\widetilde{B}\times \textbf{0}_{n-k_2}) ) \perp
\Psi_1(\Psi_2^{-1}(\textbf{0}_{k_2} \times
\mathbb{R}^{n-k_2})).\nonumber\end{align} Then    $\Psi_1(A)$ and
$\Psi_1(B)$ are of the form
\begin{align}
\widetilde{A}\times\mathbb{R}^{n-k_1}\ \ and \ \
\Psi_1(\Psi_2^{-1}(\textbf{0}_{k_2} \times \mathbb{R}^{n-k_2})
)\oplus \Psi_1(\Psi_2^{-1}(\widetilde{B}\times \textbf{0}_{n-k_2}) )
,\ \ \label{00}
\end{align} respectively.  Notice  that
$ \Psi_1( \Psi_2^{-1}(\widetilde{B}\times \textbf{0}_{n-k_2}) )
\perp\widetilde{A}\times \textbf{0}_{n-k_1}$ does not hold.
Otherwise, by (\ref{00}), $A$ and $B$ are  unlinked which
contradicts the assumption.
 Therefore,  there
exist some $1\leq j_1\leq k_1$ and $1\leq j_2\leq k_2$ such that
 \begin{align}\nonumber\langle \mathbf{e}_{j_1},\Psi_1(\Psi_2^{-1}(\mathbf{e}_{j_2}))\rangle\neq
 0,\end{align}
   which implies   that there exists some  $i_0\in\{1,\cdots,n\}$ such
that
\begin{align}\label{or}\langle \Psi_1(\mathbf{e}_{i_0}), \mathbf{e}_{j_1}\rangle\neq
0,\ \ \ \ \ \ \langle
\Psi_1(\mathbf{e}_{i_0}),\Psi_1(\Psi_2^{-1}(\mathbf{e}_{j_2})\rangle\neq
0.\end{align}

Applying orthogonal transformation $\Psi_1$, Fubini theorem, the
first equality of (\ref{4n}), the first property  of  (\ref{or}) and
Lemma \ref{po}, we have
\begin{align}\int_{A} (1-|x_{i_0}|^2)d\mu_{n}(x)= & \int_{\Psi_1(A)} (1-|
(\Psi_1^{-1}(y))_{i_0}|^2)d\mu_{n}(y)\nonumber\\
=& \int_{\widetilde{
A} \times \mathbb{R}^{n-k_2}} (1-|(\Psi_1^{-1}(y))_{i_0}|^2) d\mu_{n}(y)\nonumber\\
 =& \mu_n(A)-\int_{\widetilde{
A} \times \mathbb{R}^{n-k_2}} | \langle y,  \Psi_1 (\mathbf{e}_{i_0}))\rangle|^2d\mu_{n}(y)  \nonumber\\
 =&\mu_n(A)- \int_0^\infty dr\int_{\widetilde{
A} \times \mathbb{R}^{n-k_2}}  I_{| \langle y, \Psi_1 (\mathbf{e}_{i_0})\rangle|2>r}d\mu_{n}(y)  \nonumber\\
>&\mu_n(A)- \mu_n(A)\int_0^\infty dr\int  I_{| \langle y,
\Psi_1 (\mathbf{e}_{i_0})\rangle|2>r}d\mu_{n}(y)  \nonumber\\
=&0.\label{phi2}
\end{align}
Similarly, by the second  property  of  (\ref{or})  we also  have
\begin{align}\int_{B} (1-|x_{i_0}|^2)d\mu_{n}(x)>0.\label{pphi2}
\end{align}
By (\ref{21}),  we  also  have for $i=1,\cdots,n$
\begin{align}
 \int_A (1-|x_i|^2)d\mu_n(x)\geq 0,\ \ &\  \int_B (1-|x_i|^2) d\mu_{n}(
x)\geq 0, \label{bi}. \end{align} Applying (\ref{phi2})-(\ref{bi}),
we get
\begin{align}
 &\sum_{i=1}^n\int_A (1-|x_i|^2)d\mu_n(x) \int_B (1-|y_i|^2) d\mu_{n}(
y)\nonumber\\
\geq &\int_A (1-|x_{i_0}|^2)d\mu_n(x) \int_B (1-|y_{i_0}|^2)
d\mu_{n}( y) \nonumber\\
>&0,\nonumber
\end{align} which gives (\ref{greater}).\qed\medskip

\begin{lemma}\label{second}
Let   $A,B\in {\mathcal{C}}_n$.  Then  \begin{align} \frac{d^2
\psi_{\lambda}(A,B)}{ d\lambda^2}_{\lambda=0} \geq 0.\label{D=0}
 \end{align}
 Suppose further that $\overline{A}$ and $\overline{B}$ are  not unlinked    and neither     of
 them is
  degenerate. Then
  \begin{align}
\frac{d^2 \psi_{\lambda}(A,B)}{ d\lambda^2}_{\lambda=0}
>0.\label{D>0}
 \end{align}

\end{lemma}\noindent\textbf{Proof}\
Applying  Lemma \ref{SD} and   Fubini theorem, we have
    \begin{align}
&\frac{d^2}{d\lambda^2}\psi_{\lambda}(A,B)_{\lambda=0}\nonumber\\
=&\int_A\int_B
(n+|\langle x,y\rangle|^2-|x|^2-|y|^2) d\mu_{2n}(x,y)\nonumber\\
=&\sum_{i=1}^n \int_A\int_B  (1-|x_i|^2)(1-|y_i|^2) d\mu_{2n}(x,y)
+2\sum_{i,j=1,i\neq j}^n\int_A\int_B x_ix_jy_iy_jd\mu_{2n}(x,y)
\nonumber\\
=&\sum_{i=1}^n \int_A(1-|x_i|^2)d\mu_n(x)\int_B   (1-|y_i|^2)
d\mu_{n}( y) +2\sum_{i,j=1,i\neq j}^n\int_A
x_ix_jd\mu_n(x)\int_By_iy_jd\mu_{n}(y). \label{dec}\end{align} Set
\begin{align}
a_{i,j}=\int_A x_ix_jd\mu_n(x),\ \ \ for\
i,j\in\{1,\cdots,n\}.\nonumber
\end{align}
Since $\int_A x_id\mu_n(x)=0$ for $1\leq i \leq n$,  $(a_{i,j})$ is
the  covariance matrix of random variables $(I_Ax_i)_{1\leq i\leq
n}$ under probability space   $(\mathbb{R}^n,\mu_n)$. Then,
$Q(a_{i,j})Q^{\tau}$ is  a diagonal matrix for  some orthogonal
transformation $Q$ of $\mathbb{R}^n$. Under the coordinate   system
$(Q(\mathbf{e}_1),\cdots,Q(\mathbf{e}_n))$, we have
  by  (\ref{dec})
\begin{align}
\frac{d^2}{d\lambda^2}\psi_{\lambda}(A,B)_{\lambda=0}=&\sum_{i=1}^n
\int_A(1-|x_i|^2)d\mu_n(x)\int_B   (1-|y_i|^2) d\mu_{n}( y).
\nonumber
    \end{align}
Then we get (\ref{D=0}) and (\ref{D>0}) by Lemma \ref{Ho19} and
Lemma \ref{De}, respectively. \qed\medskip

\subsection{\normalsize proof of   Theorem 1.1}

The following result  is a generalization of   Proposition 4 in
\cite{SSZ98}.

\begin{lemma}\label{Appr}
Let $0\leq \lambda_0 <\lambda_1\leq 1 $. Suppose that for    any
$\varepsilon>0$, there exists some   integer $n_0=n_0(\varepsilon)$
such that for every $n\geq n_0$
\begin{align}\label{ar}
\psi_{\lambda_1}(A,B)\geq \exp\{-\varepsilon
n\}\psi_{\lambda_0}(A,B),\ \ \ \ \forall\  A,B\in \mathcal{C}_n.
\end{align}
Then
\begin{align}\label{re}
\psi_{\lambda_1}(A,B)\geq  \psi_{\lambda_0}(A,B),\ \ \ \ \ \forall\
A,B\in \mathcal{C}_n,\ n\geq 1.
\end{align}
In particular, if (\ref{ar}) holds for $\lambda_0=0$ and
$\lambda_1=1$, then
\begin{align}\label{re2}
\mu_n(A\cap B)\geq  \mu_n(A)\mu_n(B),\ \ \ \  \forall\  A,B\in
\mathcal{C}_n,\ n\geq 1.
\end{align}
\end{lemma}
\noindent\textbf{Proof}\   Assume in what below that $0\leq
\lambda_0\leq \lambda_1\leq 1$. Let $A,B\in \mathcal{C}_n$ and
$m\geq 1$. Set $A_m=\prod_{1\leq i\leq m}A\in \mathcal{C}_{mn},
B_m=\prod_{1\leq i\leq m}B\in \mathcal{C}_{mn}$.
  By Fubini theorem and (\ref{definition}), we have for $\lambda\in [0,1)$
\begin{align} \psi_\lambda(A_m,B_m)=&\int_{A_m}\int_{B_m}
f_{2mn}((\textbf{x}_i)_{i=1}^m,(\textbf{y}_i)_{i=1}^m;\lambda)
\prod_{i=1}^m d\textbf{x}_id\textbf{y}_i
\nonumber\\
=&\int_{A_m}\int_{B_m} \prod_{i=1}^m
f_{2n}(\textbf{x}_i,\textbf{y}_i;\lambda)d\textbf{x}_id\textbf{y}_i
\nonumber\\
=& \Big(\int_{A}\int_{B}
f_{2n}(\textbf{x}_i,\textbf{y}_i;\lambda)d\textbf{x}_id\textbf{y}_i\Big)^m\nonumber\\=&\psi_\lambda(A,B)^m,
\label{mn}
\end{align}
where $\textbf{x}_i,\textbf{y}_i\in \mathbb{R}^n$ for $1\leq i \leq
m$.  Notice that the  inequality above holds also for $\lambda=1$.

Let $\varepsilon>0$. By (\ref{ar}) and (\ref{mn}), when  $mn\geq
n_0$ we have
\begin{align}
 \psi_{\lambda_1}(A,B)=&\psi_{\lambda_1}(A_m,B_m)^{1/m}\nonumber\\\geq &\big(\exp\{-\varepsilon
mn\}\psi_{\lambda_0}(A_m,B_m)\big)^{1/m}\nonumber\\
 =&\exp\{-\varepsilon
n\}\psi_{\lambda_0}(A,B). \label{ine}
\end{align}
Noticing that $\varepsilon$ can be taken arbitrary small,    the
conclusion  (\ref{re}) follows  by  (\ref{ine}).  The last
conclusion is a direct consequence of (\ref{3'}) and
(\ref{re}).\qed\medskip

\begin{lemma}  \label{ab,1,1} Let
 $\alpha>0$ and  $\delta\in (0,1)$.
  Let $A\in \mathcal{C}_n$, $u\in \mathcal{CF}_n$ and assume that
  $B_n(\delta\sqrt{n})\subseteq
 A $. Then, there exists some integer  $N_{10}=N_{10}(\delta)$ such that  for every  $n\geq N_{10}$
\begin{align}
    \int h_{A, \alpha}u
  d\mu_n \leq 4  \int I_A u
  d\mu_n . \label{HAB1,0,1}\end{align}
For every  $n\geq 1$, we also have \begin{align}
  \int
   h_{A, \alpha} ud\mu_n\geq &
(1+2\alpha)^{-n/2}\int    I_{A }u d\mu_n.\label{HAB1,0,2}
\end{align}

 \end{lemma}
 \noindent{\bf Proof}\
Let   $\widetilde{x}\in S_{n-1}$. Set
\begin{align}
r_0(\widetilde{x})=&\inf \{r\geq 0:r\widetilde{x}\in A  \}.\nonumber
\end{align}
  From  the definition above and the
assumption $B_n(\delta\sqrt{n})\subseteq
 A$, we get
\begin{align}
r_0( \widetilde{x})\widetilde{x} \in & \partial A,\ \ \ \ \  r_0(
\widetilde{x})\geq \delta \sqrt{n}.\ \ \ \ \
 \label{t1,t2}
\end{align}
 Let $\Sigma_{\widetilde{x}}$ be a  supporting  hyperplane of $A$
  containing  $r_0( \widetilde{x})\widetilde{x}$ and let
$\mathbf{e}_{\widetilde{x}}$ be the unit normal vector of
$\Sigma_{\widetilde{x}}$ such that $\langle
\mathbf{e}_{\widetilde{x}}, \widetilde{x}\rangle>0$. Then, from
$B_n(\delta\sqrt{n})\subseteq
 A $ and $r_0( \widetilde{x})\widetilde{x}\in \partial A$, we have   $\Sigma_{\widetilde{x}}\cap
B_n(\delta \sqrt{n})=\emptyset$.  Therefore,
   \begin{align}
  \langle \mathbf{e}_{\widetilde{x}},r_0( \widetilde{x}){\widetilde{x}}\rangle \geq  \delta\sqrt{n} .\label{v,x}
 \end{align}
 Denote  $\Lambda_{\widetilde{x}}=\{y:\langle y- r_0( \widetilde{x})\widetilde{x}, \mathbf{e}_{\widetilde{x}}\rangle \leq
 0\}$.
 We have
 $A\subseteq \Lambda_{\widetilde{x}}$ by the convexity of $A$ and definition of $\mathbf{e}_{\widetilde{x}}$.
  Then, with the help of (\ref{v,x}), we have for every   $s>0$
    \begin{align}
  \rho_A( (r_0( \widetilde{x})+s)\widetilde{x} ) \geq
    \rho_{\Lambda_{\widetilde{x}}}((r_0( \widetilde{x})+s)\widetilde{x}  )
     =\rho_{\Sigma_{\widetilde{x}}}((r_0( \widetilde{x})+s)\widetilde{x})
  =s\langle \mathbf{e}_{\widetilde{x}}, \widetilde{x}\rangle  \geq
  \frac{\sqrt{n}}{r_0( \widetilde{x})}\delta s,\nonumber
 \end{align}
which gives
\begin{align}\rho_A( (r_0( \widetilde{x})+s)\widetilde{x} ) \geq
\delta     s,\ \ \ if\ r_0(\widetilde{x})\leq \sqrt{n}\ and\ s\geq
0.\label{t12}\end{align}

Noticing that      $ u( r\widetilde{x})$ is a  decreasing function
of $r\in [0,\infty)$, we have
    \begin{align}
   \int_{0}^{r_0(\widetilde{x})}
  u(r\widetilde{x})r^{n-1}
\exp\{-\frac{r^2}{2}\} dr
   \geq &   u(
r_0(\widetilde{x})\widetilde{x})\exp\{-\frac{r_0(\widetilde{x}
)^2}{2}\}\int_{0\vee(r_0(\widetilde{x})
-n^{-1/2})}^{r_0(\widetilde{x})}
  r^{n-1}
dr. \label{t0,5}
  \end{align}
By the inequality in (\ref{t1,t2}) and $0<\delta<1$, we have for
every  $n\geq 2\delta^{-1}$
  \begin{align}
(r_0(\widetilde{x})-n^{-1/2})^{n-1}\geq &
r_0(\widetilde{x})^{n-1}(1-r_0(\widetilde{x})^{-1}n^{-1/2})^{n-1}
\nonumber\\\geq&
r_0(\widetilde{x})^{n-1}(1-\delta^{-1}n^{-1})^n\nonumber\\\geq&
 e^{-2/\delta}r_0(\widetilde{x})^{n-1}.\nonumber
  \end{align}
   This   and (\ref{t0,5})  show that  for every  $n>2\delta^{-1}$
     \begin{align}
   \int_{0}^{r_0(\widetilde{x})}
  u(r\widetilde{x})r^{n-1}
\exp\{-\frac{r^2}{2}\} dr  \geq &\frac{1}{2}e^{-1/\delta} n^{-1/2}
u( r_0(\widetilde{x})\widetilde{x})
\exp\{-\frac{r_0(\widetilde{x})^2}{2}\} r_0(\widetilde{x})^{n-1},
\nonumber
  \end{align}
which gives
   \begin{align}
  &\int_{A\cap B_n(\sqrt{n})}
    u d\mu_n\nonumber\\
   \geq &\frac{1}{(2\pi)^{n/2}}
\int_{S_{n-1}} dm_{n-1}(\widetilde{x})
\int_{0}^{r_0(\widetilde{x})}u(r\widetilde{x})r^{n-1}
\exp\{-\frac{r^2}{2}\}I_{r_0(\widetilde{x})\leq \sqrt{n}}
dr  \nonumber \\
  \geq &\frac{ e^{-2/ \delta }}{n^{1/2}(2\pi)^{n/2}}
\int_{S_{n-1}} u(r_0(\widetilde{x})\widetilde{x})
\exp\{-\frac{r_0(\widetilde{x})^2}{2}\}
r_0(\widetilde{x})^{n-1}I_{r_0(\widetilde{x})\leq \sqrt{n}}
dm_{n-1}(\widetilde{x}).\label{t0,4,r}
  \end{align}

 If  $\delta^2 \sqrt{n}\geq 2$,  by
   the inequality in  (\ref{t1,t2}),  we have  for every $r\geq r_0(\widetilde{x})$
 \begin{align}
\frac{d}{dr}\big(-(n-1)\ln r+\frac{1}{2}(r-r_0(\widetilde{x}))\delta
 n\big) =& - \frac{n-1}{r}+\frac{1}{2}\delta  n\nonumber\\
\geq & - \frac{n}{r_0(\widetilde{x})}+\delta^{-1}
 \sqrt{n}\nonumber\\
 \geq &0,\nonumber
\end{align}
which gives
\begin{align}
 -(n-1)\ln r+\frac{1}{2}(r-r_0(\widetilde{x}))\delta  n \geq&
-(n-1)\ln r_0(\widetilde{x}),\ \ \ \ for\ r\geq r_0(\widetilde{x})
.\nonumber
\end{align}
If  $\delta^2 \sqrt{n}\geq 2$, applying (\ref{t12}) and the estimate
above, we obtain
\begin{align}
 & \int_{A^c\cap B_n(\sqrt{n})}
   \exp\{-   n \rho_{A }(x)\}u(x)
   d\mu_n \nonumber\\
   \leq &\frac{1}{(2\pi)^{n/2}}
\int_{S_{n-1}}
dm_{n-1}(\widetilde{x})\int_{r_0(\widetilde{x})}^{r_0(\widetilde{x})\vee
\sqrt{n}}\exp\{- \delta n(r-r_0(\widetilde{x})) \}
u(r\widetilde{x})r^{n-1} \exp\{-\frac{r^2}{2}\}
dr  \nonumber \\
  \leq &\frac{ 1}{(2\pi)^{n/2}}
\int_{S_{n-1}} u( r_0(\widetilde{x})\widetilde{x})
\exp\{-\frac{r_0(\widetilde{x})^2}{2}\}dm_{n-1}(\widetilde{x})\int_{r_0(\widetilde{x})}^{r_0(\widetilde{x})\vee
\sqrt{n}}\exp\{- \delta n(r-r_0(\widetilde{x})) \} r^{n-1}dr\nonumber \\
  \leq &\frac{ 1}{(2\pi)^{n/2}}
\int_{S_{n-1}}u (r_0(\widetilde{x})\widetilde{x})
\exp\{-\frac{(r_0(\widetilde{x}))^2}{2}\}
r_0(\widetilde{x})^{n-1}dm_{n-1}(\widetilde{x})
\int_{r_0(\widetilde{x})}^{r_0(\widetilde{x})\vee
\sqrt{n}}\exp\{-\frac{1}{2} \delta n(r-r_0(\widetilde{x})) \} dr\nonumber \\
  \leq &\frac{ 2}{ \delta n(2\pi)^{n/2}}
\int_{S_{n-1}}u (r_0(\widetilde{x})\widetilde{x})
\exp\{-\frac{(r_0(\widetilde{x}))^2}{2}\}
r_0(\widetilde{x})^{n-1}I_{r_0(\widetilde{x})\leq
\sqrt{n}}dm_{n-1}(\widetilde{x})  .\label{t0,4}
  \end{align}

Combing (\ref{t0,4,r}) and (\ref{t0,4}), we have for $n$ big enough
depending on $\delta$
\begin{align}
 & \int_{A^c\cap B_n(\sqrt{n})}
   \exp\{-  n \rho_{A }(x)\}u(x)
   d\mu_n\leq  \int_{A\cap B_n(\sqrt{n})}
   ud\mu_n.\nonumber
  \end{align}
Applying (\ref{21}), (\ref{1/2}) and the estimate above, for  $n
\geq N_{10}$ we have
\begin{align}
 \int
   h_{A,\alpha}u\mu_n
   \leq &
   \frac{1}{\mu_n(B_n(\sqrt{n}))}\int_{B_n(\sqrt{n})}
    h_{A,\alpha}u d\mu_n
  \nonumber\\ \leq  & \frac{1}{\mu_n(B_n(\sqrt{n}))}\Big(\int_{A\cap B_n(\sqrt{n})}
 u  d\mu_n
   +\int_{A^c\cap B_n(\sqrt{n})}
   \exp\{-  n \rho_{A } \}u
   d\mu_n\Big)\nonumber \\ \leq  & \frac{2}{\mu_n(B_n(\sqrt{n}))}
   \int_A
   u  d\mu_n\nonumber \\ \leq  & 4
   \int_A
   u  d\mu_n  \nonumber,
  \end{align}
which gives   (\ref{HAB1,0,1}).

 Applying Lemma \ref{B}  and Fubini theorem , we get  \begin{align}
         \int h_{A, \alpha}u
  d\mu_n\geq &\int \exp\{- \alpha|x|^2\}I_A(x) u(x)
  d\mu_n\nonumber\\
  \geq& \int \exp\{-\alpha|x|^2\}
  d\mu_n \int  I_A u
  d\mu_n\nonumber\\ =& (1+2\alpha)^{-n/2}\int  I_Au
  d\mu_n,\nonumber
\end{align}
 which gives (\ref{HAB1,0,2}).
\qed\medskip

\begin{corollary}  \label{ab,1,1,1} Let
 $\alpha>0$ and  $\delta\in
(0,1)$.
  Let $A,B\in \mathcal{C}_n$  and assume that
  $B_n(\delta\sqrt{n})\subseteq
 A\cap B $. Then  we have  for every   $t\geq 0$
 and   every $n\geq N_{10}$
\begin{align}
    2^4\int I_A P_tI_{B }
  d\mu_n \geq &
    \int h_{A, \alpha}P_th_{B, \alpha}
  d\mu_n, \label{HAB1,0,1,n1}\end{align}
Moreover,  for every  $n\geq 1$ and every $t\geq 0$ we have
\begin{align}
  \int
   h_{A, \alpha}P_th_{B, \alpha}d\mu_n\geq &
(1+2\alpha)^{-n}\int    I_{A }P_tI_{B } d\mu_n.\label{HAB1,0,2,n2}
\end{align}

 \end{corollary}
 \noindent{\bf Proof}\ Let $t\geq 0$.  By (\ref{t/2}) and  (\ref{HAB1,0,1}),
 we have for  every  $n\geq N_{10}$ \begin{align}
         \int h_{A, \alpha}P_th_{B,\alpha}
  d\mu_n
  \leq&  4\int  I_AP_th_{B, \alpha}
  d\mu_n \nonumber\\=&4\int h_{B, \alpha} P_tI_A
  d\mu_n \nonumber\\\leq& 2^4\int I_B P_tI_A
  d\mu_n= 2^4\int I_A P_tI_B
  d\mu_n.\nonumber
\end{align}which gives
(\ref{HAB1,0,1,n1}).
 Similarly,  we have  by (\ref{t/2}) and  (\ref{HAB1,0,2})  \begin{align}
         \int h_{A, \alpha}P_tH_{B,\alpha}
  d\mu_n
  \geq&  (1+2\alpha)^{-n/2}\int  I_AP_th_{B, \alpha}
  d\mu_n \nonumber\\=& (1+2\alpha)^{-n/2}\int h_{B, \alpha} P_tI_A
  d\mu_n
  \nonumber\\\geq &(1+2\alpha)^{-n}\int I_B P_tI_A
  d\mu_n  \nonumber\\= &(1+2\alpha)^{-n}\int I_AP_tI_B
  d\mu_n.\nonumber
\end{align}which gives
(\ref{HAB1,0,2,n2}). \qed\medskip

\begin{lemma}\label{unlink}
Let $\alpha>0 $ and $A,B\in \mathcal{C}_n$. Then for every $r
> \max\{ H_{A,\alpha}(0), H_{B,\alpha}(0)\}$, the set  $\{x :H_{A,\alpha}(x)\leq r\}$ and
the set  $\{ :H_{B,\alpha}(x)\leq  r\}$ are not unlinked.
\end{lemma}
\noindent\textbf{Proof}\  If two   sets are unlinked, then by
definition at least one of  them  is    unbounded.  Notice that for
every  $r
> \max\{ H_{A,\alpha}(0), H_{B,\alpha}(0)\}$, the two sets
$\{x\in \mathbb{R}^n:H_{A,\alpha}(x)\leq r\}$ and $\{x\in
\mathbb{R}^n:H_{B,\alpha}(x)\leq r\}$ are both bounded. Therefore,
$\{x\in \mathbb{R}^n:H_{A,\alpha}(x)\leq r\}$ and $\{x\in
\mathbb{R}^n:H_{B,\alpha}(x)\leq  r\}$ are not unlinked.
\qed\medskip

  \noindent\textbf{Proof  for the first conclusion  of Theorem
1.1}\    Define for every  $D\in \mathcal{C}_n$ and every $r>0$
\begin{align}
D^{[r]}=\{x :H_{D,\alpha}(x)\leq r\}.\nonumber
\end{align}
Let $A,B\in \mathcal{C}_n$ and   $\alpha>0$. We have for every $t>0$
 \begin{align}
 h_{A, \alpha,t}(x)  =&P_t h_{A,\alpha}(x)\nonumber\\
 =&\frac{1}{(1-e^{-t})^{n/2}(2\pi)^{n/2}}\int \exp\{-H_{A,\alpha}(y)\}\exp\{-\frac{
 |y-e^{-t/2}x|^2}{2(1-e^{-t})}\}dy
 \nonumber\\
 =&\frac{1}{(1-e^{-t})^{n/2}(2\pi)^{n/2}}
 \int \int_0^{\exp\{-H_{A,\alpha}(y)\}}
 dr \exp\{-\frac{
 |y-e^{-t/2}x|^2}{2(1-e^{-t})}\}dy\nonumber\\
 =&\frac{1}{(1-e^{-t})^{n/2}(2\pi)^{n/2}}
 \int_0^1dr \int I_{A^{[-\ln r]}}(y)
  \exp\{-\frac{
 |y-e^{-t/2}x|^2}{2(1-e^{-t})}\}dy\nonumber\\
  =&
 \int_0^1P_t I_{A^{[-\ln r]}}(x) dr. \nonumber
\end{align}
Notice that the equality above holds also for $t=0$. By Fubini
theorem and the equality above
 \begin{align}
 \phi_t(h_{A,\alpha},h_{B,\alpha})
=&\int h_{A, \alpha,t}
h_{B,\alpha}d\mu_n\nonumber\\
=&\int \big(\int_0^1P_t I_{A^{[-\ln r_1]}}(x) dr_1\big)
\big(\int_0^1  I_{B^{[-\ln r_2]}}(x) dr_2\big)d\mu_n(x)\nonumber\\
=& \int_0^1dr_1\int_0^1dr_2 \int P_t I_{A^{[-\ln r_1]}}(x)
 I_{B^{[-\ln r_2]}}(x) d\mu_n(x)\nonumber
 \\
=& \int_0^1dr_1\int_0^1 \phi_t( I_{A^{[-\ln r_1]}},
 I_{B^{[-\ln r_2]}})dr_2, \nonumber
\end{align}
which gives
 \begin{align}
 \psi_\lambda (h_{A,\alpha},h_{B,\alpha})
=&   \int_0^1dr_1\int_0^1  \psi_\lambda( I_{A^{[-\ln r_1]}},
 I_{B^{[-\ln r_2]}})dr_2. \nonumber
\end{align}
Then we have  for every  $\lambda\in [0,1)$
 \begin{align}
 \frac{d^2}{d\lambda^2}\psi_\lambda (h_{A,\alpha},h_{B,\alpha})
=&   \int_0^1dr_1\int_0^1\frac{d^2}{d\lambda^2} \psi_\lambda(
I_{A^{[-\ln r_1]}},
 I_{B^{[-\ln r_2]}})dr_2. \label{sed}
\end{align}

Set $s_0=\max(H_{A,\alpha}(0),H_{B,\alpha}(0))$. Notice that for
every  $s>s_0$, $A^{[s]}$ and $B^{[s]}$  are both non-degenerate
symmetric convex sets. Moreover, by   Lemma \ref{unlink}, $A^{[s]}$
and $B^{[s]}$ are  not  unlinked   for every   $s
>s_0$. Therefore,  we have by
Lemma \ref{second}
\begin{align}
\big(\frac{d^2}{d\lambda^2}\psi_\lambda( I_{A^{[ -\ln r_1]}},
 I_{B^{[-\ln  r_2]}})\big)_{\lambda=0}>0,\ \ \ \ \ \forall\  r_1,r_2\in
 (0,e^{-s_0}),\nonumber
\end{align}
which further gives together with    (\ref{sed})
 \begin{align}
 \big(\frac{d^2}{d\lambda^2}\psi_\lambda (h_{A,\alpha},h_{B,\alpha})
\big)_{\lambda=0}>0. \nonumber
\end{align}
By   Lemma \ref{mid,a} and the estimate above, there exists some
constant $\lambda_0\in(0,1)$ such that
\begin{align}
 \frac{d}{d\lambda}
 \psi_\lambda (h_{A,\alpha},h_{B,\alpha}) >0,\ \ \ \ \ \  \forall\  \lambda\in
 (0,\lambda_0).\nonumber
\end{align}
Therefore, \begin{align}
 \frac{d}{dt}
 \phi_t (h_{A,\alpha},h_{B,\alpha}) <0,\ \ \ \ \ \ \forall\  t\in (
-2\ln \lambda_0,\infty).\label{inf}
\end{align}

 Let $\varepsilon\in (0,2^{-6})$ and $T(\alpha)=4(2\ln2-\ln C(\alpha))$
 with    $C(\alpha)=\min(e^{-3}\alpha,2^{-6}e^{-3})$. Set  $\alpha=\varepsilon$.
 We have
 $C(\alpha)=e^{-3}\varepsilon$ and   $T(\varepsilon) =
4(3+2\ln 2-\ln\varepsilon)$.  By     Lemma \ref{low} and
(\ref{inf}),
 \begin{align}
\frac{d^2}{dt^2}\phi_t(h_{A,\alpha},h_{B,\alpha})
>&0,\ \ \ \ \ \forall\  t\in ( T(\varepsilon),\infty),\nonumber
\end{align}which further  gives by  Lemma \ref{low}   and (\ref{inf})
 \begin{align}
\frac{d}{dt}\phi_t(h_{A,\alpha},h_{B,\alpha})
>&0,\ \ \ \ \ \    \   \forall\  t\in  ( T(\varepsilon),\infty).\nonumber
\end{align}
This
   implies
 \begin{align}
 \phi_t(h_{A,\alpha},h_{B,\alpha}) \geq
\phi_s(h_{A,\alpha},h_{B,\alpha}),\ \ \ \ if\ T(\varepsilon) \leq
t<s.\label{lim}
\end{align}

Next we assume   $n\geq
N_9(\varepsilon,4\sqrt{\varepsilon},\varepsilon)\vee
N_{10}(4\sqrt{\varepsilon})$. When $B_n(4\sqrt{\varepsilon
n})\subseteq A\cap B$, we have  by (\ref{gap}) and (\ref{lim})
 \begin{align}
 \phi_t(h_{A,\alpha},h_{B,\alpha})
>& \exp\{-\varepsilon (s-t) n\}
 \phi_s(h_{A,\alpha},h_{B,\alpha}) ,\ \ \ \ \ if \ \varepsilon \leq t\leq
 s,\label{lim2}
\end{align}
which gives
 \begin{align}
 \phi_\varepsilon(h_{A,\alpha},h_{B,\alpha})
\geq & \exp\{-\varepsilon (T(\varepsilon)-1) n\}
 \phi_{T(\varepsilon)}(h_{A,\alpha},h_{B,\alpha})
 .\nonumber
\end{align}
When $B_n(4\sqrt{\varepsilon n})\subseteq A\cap B$, applying
(\ref{lim}) and the estimate above, we further get
 \begin{align}
 \phi_0(h_{A,\alpha},h_{B,\alpha})
\geq & \exp\{-\varepsilon T(\varepsilon)
n\}\lim_{s\rightarrow\infty}
 \phi_s(h_{A,\alpha},h_{B,\alpha})
 .\nonumber
\end{align}
When $B_n(4\sqrt{\varepsilon n})\subseteq A\cap B$, the estimate
above   and Corollary  \ref{ab,1,1,1} give  \begin{align}
    2^4\int I_A I_{B }
  d\mu_n \geq  &
    \int  h_{A,\alpha} h_{B,\alpha}
  d\mu_n\nonumber\\
  \geq&  \exp\{-\varepsilon T(\varepsilon)n\}\lim_{s\rightarrow\infty}
    \int
   h_{A,\alpha}P_sh_{B,\alpha}d\mu_n
\nonumber\\
  \geq& (1+2\alpha)^{-n}  \exp\{-\varepsilon T(\varepsilon) n\}\lim_{s\rightarrow\infty}
 \int    I_{A }P_sI_{B } d\mu_n\nonumber\\
  \geq&  \exp\{-\varepsilon  (T(\varepsilon)+2)
n\}\lim_{s\rightarrow\infty}
 \int    I_{A }P_sI_{B } d\mu_n\nonumber\\
=& \exp\{-\varepsilon  (T(\varepsilon)+2) n\}
 \int    I_A  d\mu_n \int    I_B  d\mu_n\nonumber\\
  \geq& \exp\{-4\varepsilon  (4+2\ln 2-\ln\varepsilon)
n\}
 \int    I_A  d\mu_n \int    I_B  d\mu_n,\nonumber
 \end{align} For  every  $A,B\in \mathcal{C}_n$, the  estimate above and Corollary
\ref{comb}  give
\begin{align}
     \int I_A I_{B }
  d\mu_n \geq  2^{-4} \exp\{-4\big(3\sqrt{\varepsilon}+\varepsilon  (4+2\ln 2-\ln\varepsilon)
\big)n\}
 \int    I_A  d\mu_n \int    I_B  d\mu_n,\label{T1}
 \end{align} With the estimate above,
 (\ref{1}) follows by   Lemma \ref{Appr}.
\qed
\medskip

 \noindent\textbf{Proof  for the second  conclusion  of Theorem
1.1}\
 Let $\varepsilon>0$. From the calculation to prove (\ref{T1}), for
 $n$ big enough depending on $\varepsilon$  we have
\begin{align}
 \phi_t(A,B)\geq  2^{-4}\exp\{-4\big(3\sqrt{\varepsilon}+\varepsilon  (4+2\ln 2-\ln\varepsilon)
\big)n\}\phi_s(A,B),\ \ \ if\ 0\leq t\leq s,\ A,B\in \mathcal{C}_n.
\nonumber
\end{align}   The
estimate above  and Lemma \ref{Appr} show that $\phi_t(A,B)$ is a
non-increasing function of $t\geq0$ for every $A,B\in \mathcal{C}_n$
with  $n\geq 1$. Let   $A,B\in \mathcal{C}_n$ and assume that
$\overline{A}$ and $\overline{B}$ are not unlinked and neither of
them     is
 degenerate. Then,  by Lemma \ref{mid,a} and Lemma \ref{second},
  $\phi_t(A,B)$ is strictly decreasing on $[c,\infty)$ for some  $c>0$.
 Combing these two facts above,   we get  the second conclusion of Theorem
1.1.\qed\medskip

 \begin{corollary}\label{><}
 For $u,v\in \mathcal{CF}_n$, $\psi_{\lambda}(u,v)$ is a nondecreasing
 function of $\lambda\in [0,1]$. When further assuming  that  $u$
and $v$ are both differentiable, we have
\begin{align}   \int \langle \nabla
u, \nabla v\rangle d\mu_n \geq 0 ,\nonumber
\end{align}
provided that the integral above is well defined.  \end{corollary}

 \begin{remark}\label{com}
Applying Theorem 1.1, we can prove Harg\'{e}'s correlation
inequality in \cite{Har04} when the log-concave function is
symmetric and the Gaussian measure is centered. Let $u\in
\mathcal{CF}_n$ and let  $v$ be a convex function. Following the
proof of Theorem 1.2 in \cite{Har99}, we define
\begin{align}
\xi_t=\int u(x)\exp\{-t(v(x)+v(-x))\}d\mu_n -\int u(x)d\mu_n \int
\exp\{-t(v(x)+v(-x))\}d\mu_n.\nonumber
\end{align}
By Theorem 1.1, we have $\xi_t\geq 0$ for all $t\geq0$. Since
$\xi_0=0$, we have $d\xi_t/dt\geq 0$ for $t=0$, which gives
\begin{align}
 \int  u(x) v(x)d\mu_n=&\frac{1}{2}\int  u(x)
 (v(x)+v(-x))d\mu_n\nonumber\\
 \leq & \frac{1}{2}\int u(x)d\mu_n \int
(v(x)+v(-x))d\mu_n\nonumber\\= &\int  u(x)\int  v(x)d\mu_n.\nonumber
\end{align}
Here we assume that all the integrals above are well defined.
\end{remark}

\section{Some applications}

\subsection{\normalsize the Gaussian correlation inequality on Wiener space}

 The
following Theorem  verifies a conjecture which is stated   in (2.6)
of \cite{LS01} and conjecture 6.1 of  \cite{La} for instance,  with
additional measurable assumption.
 \begin{theorem}\label{C}Let $F$ be a separable Banach space. Assume that $\mu$ is  a Wiener measure on $F$
 and $\mathcal{F}$ is the    Borel $\sigma$-algebra of
 $(F,\mu)$. Then for   any     symmetric convex sets $A$ and $B$ in $\mathcal{F}$
 \begin{align}\label{hi}
\mu(A\cap B)\geq \mu(A)\mu(B).
 \end{align}

\end{theorem}\noindent{\bf Proof}\ Since
the Wiener measure on $F$ is a Radon measure, c.f.  \cite{LT91},
 there exists a sequence of compact sets $A_n$ of $F$ such that
$A_n\subseteq A$ and $\lim_{n\rightarrow \infty}\mu(A_n)=\mu(A)$.
Similarly, there exists a sequence of compact sets $B_n$ such that
$B_n\subseteq A$ and $\lim_{n\rightarrow \infty}\mu(B_n)=\mu(B)$.
Denote the convex hull of a set $D$ by $Conv(D)$.  We see  that
$Conv((-A_n)\cup A_n  )$ is a symmetric compact set. Noticing that
$\lim_{n\rightarrow \infty}\mu(Conv((-A_n)\cup A_n ))=\mu(A)$
 and $\lim_{n\rightarrow \infty}\mu(Conv((-B_n)\cup B_n ))=\mu(B)$,
  it is sufficient to  prove (\ref{hi}) for symmetric
compact sets.

Suppose in what below that $A$ and $B$ are symmetric compact sets of
$F$. Since $F$ is separable, there exist continuous    linear
functions $(l_n)_{n\geq 1}$ and $(l_n')_{n\geq 1}$ such that
$A=\cap_{n\geq 1} \{ w\in F:|l_n(w)|\leq 1\}$ and  $B=\cap_{n\geq 1}
\{ w\in F:|l_n'(w)|\leq 1\}$. Noticing that, for any $n\geq 1$,
$(l_1,\cdots,l_n, l_1',\cdots,l_n')$ is a  $2n-$dimensional Gaussian
random vector defined on  $(F,\mu)$, we have by Theorem 1.1,
\begin{align}&\mu( w\in F:|l_{k}(w)|\leq 1,\ |l_k'(w)|\leq 1,1\leq k\leq
n\}) \nonumber\\\geq& \mu(w\in F:|l_{k}(w)|\leq 1,\
 1\leq k\leq n)\mu( w\in F:|l_{k}'(w)|\leq 1,\
 1\leq k\leq n) .\nonumber \end{align} Therefore we get  (\ref{hi})   by taking
$n\rightarrow \infty$ in the   inequality above. \qed \medskip

\subsection{\normalsize
 a  spectral gap inequality   of Dirichlet  Laplacian and  a correlation inequality for  subordinate Brownian motion }

First we show that for any open  convex set $A$ which is not equal
to $\mathbb{R}^n$,   $ \overline{A}$ is also not equal to
$\mathbb{R}^n$.  Otherwise there exists some   $x_0\in \overline{A}$
with $x_0\in A^c$. Choose $y_0\in A$ and $\varepsilon$ small enough
such that $y_0+B_{n}(\varepsilon)\subseteq A$. From the convexity
assumption of $A$, $(2x_0-y_0-B_{n}(\varepsilon))\cap A =\emptyset$
which shows that  $ \overline{A}$ is not   equal to $\mathbb{R}^n$.
For any open set $A\subseteq
 \mathbb{R}^n$, denote by  $\lambda_1(A)$ the the first nonzero eigenvalue of the
Laplacian on $A$ under  Dirichlet boundary condition when it exists.
By Lemma \ref{Deco} and the property above, we see that
$\lambda_1(A)>0$ if $ A$ is an open set belonging to $\mathcal{C}_n$
and it is not equal to $\mathbb{R}^n$.

\begin{corollary}\label{La}Let $A$ and $B$ be  two   open  sets  in $\mathcal{C}_n$
and assume that  neither  of them is  equal to $\mathbb{R}^n$. Then
  \begin{align}\label{1.1}
\lambda_1(A\cap B)\leq \lambda_1(A)+\lambda_1(B).
\end{align}
\end{corollary}
\noindent{\bf Proof}\   Denote for every $T>0$
\begin{align}C([0,T];\mathbb{R}^n)=\{f:f\ is \ a \ continuous\
function\ from \ [0,T]\ to\ \mathbb{R}^n\}.\end{align} Similarly,
denote by  $C([0,\infty);\mathbb{R}^n)$ the set of continuous pathes
in  $\mathbb{R}^n$ parameterized by $[0,\infty)$. Denote by $(B_t)$
the standard coordinate  Brownian motion on
$C([0,\infty);\mathbb{R}^n)$.  The  distribution of $(B_t)$ on
$C([0,\infty);\mathbb{R}^n)$ and $C([0,T];\mathbb{R}^n)$,  $T>0$,
are denoted by $\mathbb{P}$ and $\mathbb{P}^T$, respectively. For
every $T>0$, we take $C([0,T];\mathbb{R}^n)$   as a separable Banach
space with  $L^\infty$ norm. For an open domain $D\subseteq
\mathbb{R}^n$ and $s>0$, set $\tau_D=\inf \{t\geq 0: B_t \in D^c \}$
and $\tau_D^s=\inf \{0\leq t\leq s: B_t \in D^c \}$. Here the
infimum of an empty set is assumed to be  infinity.  We know that
(c.f. \cite{BA95})
\begin{align}\label{stop}
-\lim_{t\rightarrow \infty}t^{-1}\ln \mathbb{P} (\tau_D>t )
=\lambda_1(D).
\end{align}
provided that the spectral gap exists. Noticing that
$\{\tau_A^t=\infty \}$ and $\{\tau_B^t=\infty \}$ are both    open
and symmetric convex sets of $C([0,t];\mathbb{R}^n)$, we have by
(\ref{hi}) and  (\ref{stop})
\begin{align}
\lambda_1(A\cap B)=&-\lim_{t\rightarrow \infty}t^{-1} \ln
\mathbb{P}(\tau_{A\cap B}>t )   \nonumber\\ = & -\lim_{t\rightarrow
\infty}t^{-1} \ln \mathbb{P}^t(\tau_A^t=\infty,\tau_B^t=\infty)
 \nonumber\\ \leq& -\lim_{t\rightarrow
\infty}t^{-1} \ln \mathbb{P}^t(\tau_A^t=\infty ) -\lim_{t\rightarrow
\infty}t^{-1} \ln \mathbb{P}^t(\tau_B^t=\infty)\nonumber\\
=& -\lim_{t\rightarrow \infty}  t^{-1} \ln \mathbb{P}(\tau_{A }>t \}
-\lim_{t\rightarrow \infty}  t^{-1} \ln \mathbb{P}(\tau_{B
}>t \}\nonumber\\
 = & \lambda_1(A)+\lambda_1(B),\nonumber
\end{align}
which gives the conclusion.\qed
\medskip

In what below, a measurable function $f$ on $\mathbb{R}^n$ is called
non-decreasing if $f(x_1,\cdots,x_n)\geq f(y_1,\cdots,y_n)$ provided
that  $x_i\geq y_i $  for  every   $ i=1,\cdots,  n$. The following
inequality is a special case of   FKG inequality on product spaces,
c.f. \cite{KR80}.
\begin{lemma}\label{FKG}
Let   $d\nu= \prod_{1\leq i\leq n}d\nu_i$, where $\nu_i$ is a
probability measures on $\mathbb{R}$ for each $i\in \{1,\cdots,
n\}$. Let $f$ and $g$ be non-decreasing functions on $\mathbb{R}^n$.
Then
  \begin{align}\label{FKGH}
\int fg d\nu\geq \int f d\nu\int gd\nu,
\end{align}
provided that both sides above are well defined.
\end{lemma}

\begin{lemma}\label{FKF}
Let $m\geq1$   and $(B_{i,t})_{1\leq i \leq m}$ be $m$ independent
Brownian motions  on $\mathbb{R}^{n}$. Let   $(T_i)_{1\leq i \leq
m}$ be $m$  independent  nonnegative random variables and assume
that all of  them are  independent with $(B_{i,t})_{1\leq i \leq
m}$. Set $X=(B_{T_1},\cdots, B_{T_m})$ and denote the distribution
of $X$ by $P_X$. Then  for every  $A,B\in \mathcal{C}_{nm}$
  \begin{align}\label{1.1.1}
P_X(A\cap B)\geq P_X(A)P_X(B).
\end{align}
\end{lemma}
\noindent{\bf Proof}\   By the scaling property of $(B_t)$,   for
every  $D\in \mathcal{C}_{nm}$ and
$t_1,\cdots,t_m,s_1,\cdots,s_m>0$,
 \begin{align}\label{increa}
  \mathbb{P}\Big(  (B_{1,t_1},\cdots,B_{m,t_m})\in D
\Big) =& \mathbb{P}
\Big(  ((\frac{t_1}{s_1})^{1/2}B_{1,s_1},\cdots, (\frac{t_m}{s_m})^{1/2} B_{m,s_m})\in D  \Big)\nonumber\\
=& \mathbb{P}\Big(  (B_{1,s_1},\cdots,  B_{m,s_m})\in D'  \Big),
\end{align}
where $D'=\{\mathbf{x}:
\mathbf{x}=\big((\frac{s_1}{t_1})^{1/2}y^{(1)},\cdots,
(\frac{s_m}{t_m})^{1/2}y^{(m)}\big),\ \    (y^{(i)})_{i=1}^m\in D\
with\ y^{(i)}\in \mathbb{R}^n\ for\ 1\leq i\leq m\ \}$. When $0<
t_i< s_i$ for each $i\in (1,\cdots,m)$, we have $D\subseteq D'$ and
hence by (\ref{increa})
\begin{align}\label{Compare}
\mathbb{P}\Big( (B_{1,t_1},\cdots,B_{m,t_m})\in D \Big) \geq&
\mathbb{P}\Big( (B_{1,s_1},\cdots,  B_{m,s_m})\in D \Big).
\end{align}
Set \begin{align}f_D(t_1,\cdots,t_m)=&\mathbb{P}\big(
 (B_{1,t_1},\cdots,B_{m,t_m})\in D \big),\ \ \ if\ t_i>0\ for\ i=1,\cdots,m;\nonumber\\
 f_D(t_1,\cdots,t_m)=&0,\ \ \
 \ \ \ \ \ \ \ \ \ \ \ \ \ \ \ \ \ \ \ \
 \ \ \ \ \ \ \ \ \ \ \ \ \ \ \ otherwise.\nonumber\end{align} Applying  (\ref{Compare}),
we have that $f_D(t_1,\cdots,t_m)$ is a non-increasing function of
$t_1,\cdots,t_m\geq 0$. Denote by $\nu$ the distribution of
$(T_i)_{1\leq i\leq m}$ on $\mathbb{R}^m $. By independent
assumptions of $(T_i)_{1\leq i\leq m}$, $\nu$ is a product
probability measure on $\mathbb{R}^m $.

To simply notations,  for a set $D$ we denote $I_D$ by $\chi(D)$  in
what below.
 For every $A,B\in \mathcal{C}_{nm}$, we have   by Theorem 1.1 and the  FKG inequality (\ref{FKGH})
\begin{align}
P_X(A\cap B)=&\mathbb{P}\Big( {X}\in A\cap B\Big)
\nonumber\\
=&\mathbb{E}\Big(\mathbb{E}\big( \chi
{((B_{1,t_1},\cdots,B_{m,t_m})\in A\cap B}
) |T_1=t_1,\cdots,T_m=t_m\big)\Big)\nonumber\\
\geq &\mathbb{E}\Big(\mathbb{E} \big( \chi
{((B_{1,t_1},B_{2,t_2},\cdots,B_{m,t_m})\in A }) |T_1=t_1,\cdots,T_m=t_m\big)\nonumber\\
&\ \ \ \ \ \ \ \
\cdot\mathbb{E} \big( \chi {((B_{1,t_1},B_{2,t_2},\cdots,B_{m,t_m})\in B } )|T_1=t_1,\cdots,T_m=t_m\big)\Big)\nonumber\\
=& \int_{\mathbb{R}^m} f_A(t_1,\cdots,t_m)f_B(t_1,\cdots,t_m)d\nu\nonumber\\
\geq &\int_{\mathbb{R}^m} f_A(t_1,\cdots,t_m)d\nu
\int_{\mathbb{R}^m} f_B(t_1,\cdots,t_m)d\nu\nonumber\\
=&P_X(A)P_X(B),\nonumber\end{align}  which completes the proof. \qed
\medskip

A non-decreasing one dimensional L\'{e}vy process is called
subordinator, c.f. \cite{B96}.  A subordinator    $(X_t)$ can be
characterized by its Laplace exponent $\Psi$, i.e.,
$$
\mathbb{E}(e^{-\lambda X_t})=e^{-t\Psi(\lambda)},\ \ \ \ \forall\ t,
\lambda \geq 0.
$$
A function $\Psi$ is the Laplace exponent of a subordinator if and
only if it has the following form:
\begin{align}
\Psi(\lambda)=a \lambda +\int_{(0,\infty)} (1-e^{-\lambda x})
\Pi(dx),
\end{align}
where  $a\geq 0$ and $\Pi(dx)$ is a measure  on $(0,\infty)$ such
that $\int_{(0,\infty)} (x\wedge 1)\Pi(dx)<\infty$.

Let $(B_t) $ be  a Brownian motion on $\mathbb{R}^n$ and   $(X_t)$
be a subordinator with Laplace exponent $\Psi$ which is independent
with $(B_t) $. A subordinate Brownian motion associated with $(X_t)$
is a time changed Brownian motion defined by \begin{align}B^{
\Psi}_t= B_{X_t},\ \ \ \ \ \forall t\geq 0\label{Sb}.\end{align} The
process $(B^{ \Psi}_t)$ is
  a  L\'{e}vy process.  Denote by $D([0,\infty);\mathbb{R}^n)$
the space of c\`{a}dl\`{a}g functions from $[0,\infty)$ to
$\mathbb{R}^n$, i.e., right continuous functions with left limits
defined  on $  [0,\infty)$.  We know that
$D([0,\infty);\mathbb{R}^n)$ is  a Polish space when it is endowed
with Skorohod topology. Under this topology,
  a sequence   $(\omega_{m,t})_{m\geq 1}\in D([0,\infty);\mathbb{R}^n)$
converges to some $(\omega_t)\in D([0,\infty);\mathbb{R}^n)$ if and
only if there exists a sequence of strictly increasing  functions
$(\lambda_m(t))_{m\geq 1}$ from $[0,\infty)$ onto $[0,\infty)$ such
that
\begin{align}\label{dis}
\lim_{m\rightarrow \infty}&\sup \{ |\lambda_m(t)-t|: t\geq 0\}=0\nonumber\\
\lim_{m\rightarrow \infty}&\sup \{|\omega_{m,\lambda_m(t)}-\omega_t|
: 0\leq t\leq T\}=0,\ \ \ \ \forall\    T\in [0,\infty).
\end{align}
We refer to  \cite{B96} and   \cite{BI68} for some  properties of
subordinator and   Skorohod topology, respectively.

\begin{lemma}\label{lim3}  Let $(s_k)$
be a dense set of $[0,\infty)$ and $(\overline{\omega}_t)\in
D([0,\infty);\mathbb{R}^n)$. Assume that $(\omega_{m,t})_{m\geq
1}\in D([0,\infty);\mathbb{R}^n)$ converges to some $(\omega_t)\in
D([0,\infty);\mathbb{R}^n)$ under Skorohod topology. Assume also
that for every $m\geq 1$ and every $k=1,\cdots, m$,
 $\overline{\omega}_{s_k}=\omega_{m,s_k}$.  Then $\overline{\omega}_t=\omega_t$ for
 every  $t\geq 0$.
\end{lemma}
\noindent{\bf Proof}\
 By the assumption of convergence,
  there exists a sequence of strictly increasing functions $(\lambda_m(t))_{m\geq 1}$
from $[0,\infty)$ onto  $[0,\infty)$  such that $(\ref{dis})$ holds.
Let $0\leq t_0 < T$ for some $T>0$. For any $\varepsilon>0$, by
$(\ref{dis})$, there exists $m_0=m_0(\varepsilon)$ such that for
$m\geq m_0$
\begin{align}\label{dis2}
|\omega_{m,\lambda_m(t)}-\omega_t|+|\lambda_m(t)-t|\leq
\varepsilon,\ \ \ \ \ if\  0\leq t\leq T+1.
\end{align}
Since $(s_k)$ is a dense set, for every $\varepsilon\in (0,1/3)$
there exists $s_{k_0}$ for some integer $k_0=k_0(\varepsilon)\geq 1$
such that $t_0+\varepsilon<s_{k_0}<t_0+2\varepsilon$. Therefore, for
$m\geq m_0$, we have by (\ref{dis2}) that
\begin{align} \lambda_m(t_0) < s_{k_0}<\lambda_m(t_0+3\varepsilon) .\nonumber\end{align}
 Then, noticing that
$\overline{\omega}_{s_{k_0}}=\omega_{m_0',s_{k_0}}$ for
$m_0':=m_0\vee k_0 $, we have by (\ref{dis2}) and the right
continuity of $(\omega_t) $ and $(\overline{\omega}_t) $
\begin{align}
\omega_{t_0}=\lim_{ \varepsilon\rightarrow 0
}\omega_{s_{k_0}}=\lim_{ \varepsilon\rightarrow 0
}(\omega_{s_{k_0}}-\omega_{m_0', s_{k_0}})+\lim_{
\varepsilon\rightarrow 0 } \omega_{m_0', s_{k_0}}=\lim_{
\varepsilon\rightarrow 0 } \overline{\omega}_{
s_{k_0}}=\overline{\omega}_{t_0},\nonumber
\end{align}
which completes the proof. \qed \medskip

 Denote  by $\mathcal{F}_{\Psi}$ the
  Borel $\sigma$-algebra of $D([0,\infty);\mathbb{R}^n)$ under   Skorohod topology and  denote the
distribution of $(B^{ \Psi}_t)$ on
$(D([0,\infty);\mathbb{R}^n),\mathcal{F}_{\Psi})$ by $\mu_{\Psi}$.

\begin{corollary}\label{Ba2}
Let $\big(D([0,\infty);\mathbb{R}^n),\mathcal{F}_\Psi,\mu_\Psi\big)$
be the probability space  described  as above.  Then for    any
measurable symmetric convex sets $A$ and $B$ of $\mathcal{F}_\Psi$,
 \begin{align}\label{hi2}
\mu_\Psi(A\cap B)\geq \mu_\Psi(A)\mu_\Psi(B).
 \end{align}

\end{corollary}
\noindent{\bf Proof}\ Since $\mu_\Psi$ is a Radon measure, we can
assume that $A$ and $B$ are both symmetric compact sets in
$\mathcal{F}_\Psi$ as  in the proof of Theorem \ref{C}. Let $(s_k)$
be a dense subset of $[0,\infty)$. For any  symmetric compact set
$D\in \mathcal{F}_\Psi$, denote for every  $m\geq 1$
\begin{align}
D_m=\{ (\omega_t)_{t\geq 0}: \omega_{s_k}=\overline{\omega}_{s_k},\
1\leq k\leq m,\ \ for\ some \ (\overline{\omega}_t)_{t\geq 0}\in
D\}. \nonumber\end{align}
 Notice   that
$D\subseteq D_m$ and $D_m$ is a symmetric convex set for every
$m\geq1$.  Next  we show that \begin{align}\label{m} D=\cap_{m\geq
1} D_m.\end{align} Since      $ D\subseteq\cap_{m\geq 1} D_m$, we
only need to show the converse relation. Suppose that
$(\overline{\omega}_t)\in \cap_{m\geq 1}D_m$, then there exists $(
{\omega}_{m,t})\in D_m, m\geq 1,$ such that
$\overline{\omega}_{s_k}= {\omega}_{m,s_k}$ for $1\leq k \leq m$.
From the assumption that  $D$ is compact, there exists $(
{\omega}_t)\in D$ and a subsequence $(m_k)$ such that $(
 {\omega}_{m_k,t})$ converges to $( {\omega}_t)$ under
Skorohod topology. Therefore, by Lemma \ref{lim3} we have
$(\overline{\omega}_t)=(\omega_t)\in D$.

Set for every  $m\geq 1$
\begin{align}
D_m'=\{  (x^{(k)})_{1\leq k\leq m}:
x^{(k)}=\overline{\omega}_{s_k},\ \ for\ some \
(\overline{\omega}_t)_{t\geq 0}\in D\}. \nonumber\end{align}Let
$m\ge1$. We see that  $D_m'$ is a symmetric convex set of
$\mathbb{R}^{mn}$ and
\begin{align}
\{(B^\Psi_{s_k})_{1\leq k\leq m}\in D_m'\}=\{(B^\Psi_{t})_{t\geq
0}\in D_m\}.\label{=}
\end{align}
Denote by $F_m$  the following   transformation from $\mathbb{R}^{mn
}$ to itself
\begin{align}
F_m((x^{(k)})_{1\leq k \leq m})=(x^{(1)},x^{(2)}-x^{(1)},\cdots,
x^{(m)}-x^{(m-1)}), \ \ \ x^{(k)}\in \mathbb{R}^n,\ 1\leq k \leq
m.\nonumber
\end{align}
Set $D_m''= F_m(D_m')$.  Since $\Psi$ is linear,  $D''$ is also  a
symmetric convex set of $\mathbb{R}^{mn }$. We also have
\begin{align}
\{(B^\Psi_{s_1},B^\Psi_{s_2}-B^\Psi_{s_1},\cdots,B^\Psi_{s_m}-B^\Psi_{s_{m-1}})\in
D_m''\}= \{(B^\Psi_{s_1},B^\Psi_{s_2},\cdots,B^\Psi_{s_m})\in
D_m'\},\ \ \ \ \forall\  m\geq 1.\label{==}
\end{align}
Since $(B^\Psi_{t})$ is a L\'{e}vy process, the distribution of
$(B^\Psi_{s_1},B^\Psi_{s_2}-B^\Psi_{s_1},\cdots,B^\Psi_{s_m}-B^\Psi_{s_{m-1}})$
is the same as the distribution of
$X(s_1,\cdots,s_m):=(B^\Psi_{1,s_1},B^\Psi_{2,s_2-s_1},\cdots,B^\Psi_{m,s_m-s_{m-1}})$,
where    $(B^\Psi_{i,t}) ,1\leq i\leq m,$ are independent processes
with the same   distribution  as  $(B^\Psi_{t})$.  Then,  we have by
Lemma \ref{FKF}, (\ref{Sb}), (\ref{=}) and (\ref{==})
\begin{align}
&\mu_\Psi((B^\Psi_{t})_{t\geq 0}\in A_m,\ (B^\Psi_{t})_{t\geq
0}\in B_m)\nonumber\\
=&\mathbb{P}_X(X(s_1,\cdots,s_m)\in A_m'',\ X(s_1,\cdots,s_m)\in B_m'')\nonumber\\
\geq &\mathbb{P}_X(X(s_1,\cdots,s_m)\in A_m'')\mathbb{P}_X(
X(s_1,\cdots,s_m)\in B_m'')
\nonumber\\
= &\mu_\Psi((B^\Psi_{t})_{t\geq 0}\in
A_m)\mu_\Psi((B^\Psi_{t})_{t\geq 0}\in B_m), \nonumber
\end{align}
where $\mathbb{P}_X$ is the distribution of $X(s_1,\cdots,s_m)$.  By
(\ref{m}), we get the conclusion by letting $m\rightarrow \infty$ in
the estimate above.
 \qed
\medskip

\noindent$\mathbf{ Acknowledgement}$\ The main results of this paper
were  reported with partial proofs at   a workshop  on stochastic
analysis in BeiJing supported by the  research  group of  AMSS, CAS
on complex and high dimensional data and structure,
 and also reported subsequently  in    a seminar talk at  the school of mathematical sciences of  DaLian
University of Technology, both  in November 2011.

 %The paper would not be finished
%without the supports of my family.

\

\

\

\noindent   Address:

\

\noindent Qingyang Guan

\

\noindent Institute of Applied Mathematics

\noindent Academy of Mathematics and Systems Science

\noindent Chinese Academy of Sciences

\noindent BeiJing 100190

\noindent  China

\

\noindent Email address: guanqy@amt.ac.cn

\end{document}